\documentclass[12pt, makeidx]{amsart}
\usepackage[dvips]{graphics}
\usepackage{amsthm}
\usepackage{amsfonts}
\usepackage{amssymb}
\newtheorem{thm}{Theorem}[section]
\newtheorem{cor}[thm]{Corollary}
\newtheorem{lem}[thm]{Lemma}
\newtheorem{prop}[thm]{Proposition}
\theoremstyle{definition}
\newtheorem{defn}[thm]{Definition}
\theoremstyle{remark}
\newtheorem{rem}[thm]{Remark}
\newtheorem{eg}[thm]{Example}
\theoremstyle{construction}
\newtheorem{con}[thm]{Construction}
\theoremstyle{notation}
\newtheorem{notation}[thm]{Notation}
\theoremstyle{addendum}

\numberwithin{equation}{section} \numberwithin{figure}{section}
\newcommand{\ab}{\rm{ab}}
\newcommand{\Aut}{\rm{Aut}}

\newcommand{\expgrowth}{\rm{exp}}
\newcommand{\Orient}{\mathcal{O}}
\newcommand{\Out}{\rm{Out}}
\newcommand{\Integer}{\mathbb {Z}}
\newcommand{\modulo}{\rm {mod}}
\newcommand{\Nat}{\mathbb {N}}
\newcommand{\PG}{{\rm PG}}

\newcommand{\Real}{\mathbb{R}}

\newcommand{\abs}[1]{\left\vert#1\right\vert}

\newcommand{\norm}[1]{\left\Vert#1\right\Vert}
\newcommand{\rev}[1]{\overline{#1}}

\newcommand{\combine}[1]{\vee(#1)}
\newcommand{\combineandfold}[1]{\vee[#1]}
\newcommand{\combineC}[1]{\vee^{\circlearrowleft}(#1)}
\newcommand{\combineandfoldC}[1]{\vee^{\circlearrowleft}[#1]}

\newcommand{\diam}{{\rm Diam}}

\newcommand{\period}{{{\rm period}_f}}
\newcommand{\edgedegree}{{\rm degree}}
\newcommand{\structure}{{\rm str}}
\newcommand{\lcm}{{\rm lcm}}

\newcommand{\length}[1]{l(#1)}

\newcommand{\separate}{\diamond}

\makeindex
\begin{document}


 \baselineskip = 16pt 


\title[Preliminary Draft - Not for Circulation]{Detecting the Growth of Free Group Automorphisms by their Action on the Homology of Subgroups of Finite Index}%
\author{Adam Piggott}%
\address{Mathematical Institute, University of Oxford, 24-29 St Giles, Oxford, OX1 3LB, UK}%
\email{piggott@maths.ox.ac.uk}


\date{28 July, 2004}%
\begin{abstract}
In this paper we prove that if $F$ is a finitely generated free
group and $\phi \in \Aut(F)$ is a polynomially growing
automorphism then there exists a characteristic subgroup $S \leq
F$ of finite index such that the automorphism of $S^{\rm \ab}$
induced by $\phi$ grows polynomially of the same degree as $\phi$.
The proof is geometric in nature and makes use of Improved
Relative Train Track representatives of free group automorphisms.
\end{abstract}
\maketitle

The study of automorphisms of non-abelian free groups has been
reinvigorated in recent years by a program to understand free
group automorphisms as homotopy equivalences of finite graphs,
called topological representatives (see, for example,
\cite{ModuliOfGraphs, TrainTracks, BFH1, BFH2}). This programme is
driven largely by analogy with the study of surface automorphisms
and has led to significant progress in the field. In a series of
papers, Bestvina, Feighn and Handel have developed powerful normal
forms for topological representatives, called Improved Relative
Train Track (IRTT) representatives (in analogy with train track
representatives of surface automorphisms) \cite{TrainTracks, BFH1,
BFH2}. This technology has allowed them to prove a number of
important results, most notably the Scott Conjecture
\cite{TrainTracks} and the Tits Alternative for $\Out(F)$
\cite{BFH1, BFH2}.  In many applications, such as our Main
Theorem, the detailed structure inherent in IRTT representatives
allows one to use geometric intuition to evade difficult and
unsightly cancellation arguments.

Let $F$\index{sF@$F$} be a finitely generated non-abelian free
group and $\phi \in \Aut(F)$ an automorphism.  The growth function
$\mathcal{G}_{\phi}:\Nat \to \Nat$ of $\phi$ quantifies the rate
at which repeated application of the automorphism changes the
`size' of a basis of $F$ (see $\S$\ref{GrowthChapter}). The
asymptotic behaviour of $\mathcal{G}_{\phi}$ does not depend on
the choice of basis for $F$ and is robust when passing to
subgroups of finite index. We write $F^{\ab}$ for the
abelianisation of $F$; $\phi^{\ab}$ for the image of $\phi$ under
the natural map $\Aut(F) \to \Aut(F^{\ab})$;
$\mathcal{G}^{\ab}_{\phi}$ for the growth function of
$\phi^{\ab}$; $\simeq$ for an equivalence relation on maps $\Nat
\to \Nat$ which respects asymptotic behaviour (see
$\S$\ref{GrowthChapter}); $p_{\expgrowth}$ for the
$\simeq$-equivalence class containing $k \mapsto 2^k$; $p_d$ for
the $\simeq$-equivalence class containing $k \mapsto k^d$, for $d
\in \Nat$; and $\PG(F)$ for the subset of $\Aut(F)$ consisting of
all non-exponentially growing automorphisms.

\begin{thm}[Main Theorem\index{wMainTheorem@Main Theorem}]\label{FirstMainTheorem}
Let $F$ be a finitely generated free group and $\phi \in \PG(F)$
an automorphism with polynomial growth. There exists a
characteristic subgroup $S \leq F$ of finite index such that, for
$\theta := \phi|_S$,
\begin{equation*}
\mathcal{G}^{\ab}_{\theta} \simeq \mathcal{G}_{\theta} \simeq
\mathcal{G}_{\phi}.
\end{equation*}
\end{thm}

In the following remarks we offer some context in which to
consider the Main Theorem: automorphisms of free abelian groups
may, of course, be understood as elements of ${\rm SL}(n,
\Integer)$. The following theorem follows from the Jordan
Canonical Form Theorem for ${\rm GL}(n, \Real)$.

\begin{thm}[Traditional]\label{AbelianAutomorphismGrowthTheorem}
Let $F^{\ab}$ be a finitely generated free abelian group of rank
$n \geq 1$ and $\phi^{\ab} \in \Aut(F^{\ab})$ an automorphism.
Either $\mathcal{G}^{\ab}_{\phi} \in p_{\expgrowth}$ or there
exists an integer $\eta$ such that $1 \leq \eta < n$ and
$\mathcal{G}^{\ab}_{\phi} \in p_\eta$.
\end{thm}

\noindent In $\S$\ref{IRTTGrowthSection} below we prove, as a
corollary to the IRTT Theorem, that an analogous statement may be
made about the automorphisms of finitely generated (non-abelian)
free groups.

\begin{thm}[Bestvina, Feighn,
Handel]\label{AutomorphismGrowthTheorem} Let $F$ be a finitely
generated free group of rank $n \geq 2$ and $\phi \in \Aut(F)$ an
automorphism.  Either $\mathcal{G}_{\phi} \in p_{\expgrowth}$ or
there exists an integer $\eta$ such that $1 \leq \eta < n$ and
$\mathcal{G}_{\phi} \in p_\eta$.
\end{thm}

\noindent The Main Theorem elucidates the equality of the growth
spectra of $\PG(F)$ and $\PG(F^{\ab})$. It may also be considered
an extension of the following theorem of Grossmann
\cite{Grossman}: for each automorphism $\phi \in \Aut(F)$ there
exists a characteristic subgroup $S \leq F$ of finite index such
that $(\phi|_S)^{\ab}$ is non-trivial.  Further, the Main Theorem
shares a theme with Lubotzky's \cite{Lubotzky} characterisation of
the inner automorphisms of $F$ as those which act trivially on the
set of normal subgroups of $F$ of prime-power index.

Our proof of the Main Theorem is inspired by the following simple
observations: Let $f:G \to G$ be a topological representative of
an automorphism $\phi \in \Aut(F)$.  Each element $w \in F$
corresponds to a closed path $\rho$ in $G$. The length of
$w^{\ab}$ (the element of $F^{\ab}$ induced by $w$) in the
word-metric (with respect to some generating set $A$) is less than
the length of $w$ if and only if a generator and its inverse both
appear in the unique reduced word in $A^{\pm}$ equal to $w$. In
the topological representative this corresponds to some subpath
$\mu$ of $\rho$ being traversed first in one direction and later
in the reverse direction as $\rho$ is traversed.  We call this
`winding' and `unwinding' $\mu$.  A key observation is that for
any path $\rho$ we may construct a finite cover of $G$ such that
the winding and unwinding in $\rho$ lift to different sheets of
the cover. Now, $\mathcal{G}^{\ab}_{\phi} \not \simeq
\mathcal{G}_{\phi}$ only if for each fastest growing closed path
$\rho$ in $G$ the $f$-iterates of $\rho$ contain significant
amounts of winding and unwinding. Our strategy for proving the
Main Theorem is to construct a covering graph $\tilde{G}$ of $G$
such that, for some fastest growing closed path $\rho$, large
amounts of this winding and unwinding lift to different sheets of
$\tilde{G}$.

Incidental to the proof of the Main Theorem, we prove the
following corollary to the IRTT Theorem (see
$\S$\ref{Proofsofgrowthpropertiessection}), the analogue of which
is unknown to the Author in the case of an arbitrary finitely
generated group. We also indicate how this result may be
considered a corollary to the Main Theorem.

\begin{thm}\label{InverseHasSameGrowth}
Let $F$ be a finitely generated free group and $\phi \in \Aut(F)$
an automorphism.  Then $\mathcal{G}_{\phi} \simeq
\mathcal{G}_{\phi^{-1}}$.
\end{thm}

\begin{rem}[Algorithmic properties of the Main Theorem]\label{AlgorithmicProperties}
It is shown in \cite{TrainTracks} that there exists an algorithm
to determine whether or not an automorphism $\phi \in \Aut(F)$ is
polynomially growing. Given $\phi \in \PG(F)$ and an IRTT
representative $f:G \to G$ of some iterate of $\phi$, our proof of
the Main Theorem shows how to construct a characteristic subgroup
$S \leq F$ of finite index with the property that
$\mathcal{G}^{\ab}_{\phi|_S} \simeq \mathcal{G}_{\phi}$.
Alternatively, without necessarily knowing an IRTT representative
of an iterate of $\phi$, we may find a characteristic subgroup $S$
with the desired property by performing two partial algorithms as
follows: let $S_1, S_2, \dots$, be an enumeration of the subgroups
of $F$ of finite-index; for each $i \in \Nat$, let $k_i \in \Nat$
be such that $\phi^{k_i}(S_i) = S_i$, let $d_i \in \Nat$ be such
that $\mathcal{G}^{\ab}_{\phi^{k_i}|_{S_i}} \in p_{d_i}$ and let
$D_i := \max\{d_1, d_2, \dots, d_i\}$. Now, $\{D_i\}$ is a
non-decreasing sequence and, by the Main Theorem, there exists
$i_0 \in \Nat$ such that $i \geq i_0$ implies $\mathcal{G}_\phi
\in p_{D_i}$. In Remark \ref{UpperBoundForAlgorithm} we show how
to enumerate a non-increasing sequence $\{U_i\}$ of natural
numbers such that there exists $i_1 \in \Nat$ for which $i \geq
i_1$ implies $\mathcal{G}_\phi \in p_{U_i}$.  We may enumerate
$\{D_i\}$ and $\{U_i\}$ until $D_i = U_i$, then
$\mathcal{G}_{\phi^{k_i}|_{S_i}} \in p_{D_i}$ and the intersection
$S$ of all subgroups of $F$ of index $[F: S_i]$ is a
characteristic subgroup of finite index such that
$\mathcal{G}^{\ab}_{\phi|_S} \simeq \mathcal{G}_{\phi}$.
\end{rem}

We now outline the organisation of this paper: In
$\S$\ref{GrowthChapter} we formally introduce the growth function
of an automorphism and some basic properties. In
$\S$\ref{GraphChapter} we remind the reader of Stallings' notation
for directed graphs, Stallings' Folding Operation and Stallings'
Algorithm \cite{Stallings} for extending a graph immersion to a
graph covering. We also introduce notation for end-pointed and
base-pointed graphs, which are the building blocks of the
constructions we use to prove the Main Theorem.  In
$\S$\ref{IRTTChapter} we give a brief exposition of IRTT
representatives of automorphisms in $\PG(F)$. The obvious links
are developed between, on the one hand, the growth of paths and
circuits in an IRTT representative $f:G \to G$ and, on the other,
the growth of the automorphism $\phi \in \PG(F)$ represented by
$f$. The important notion of the reverse $\overline{f}:G \to G$ of
an IRTT representative $f:G \to G$ is introduced and we prove
theorems \ref{AutomorphismGrowthTheorem} and
\ref{InverseHasSameGrowth} and complete the discussion of Remark
\ref{AlgorithmicProperties}. In
$\S$\ref{TranslatingTheoremChapter} we use the theory developed in
$\S$\ref{IRTTChapter} to translate the Main Theorem into a theorem
stated in the language of topological representatives, the Apt
Immersion Theorem (Theorem \ref{FriendlyImmersionTheorem}). For a
path $\rho$ in an IRTT representative $f:G \to G$, the Apt
Immersion Theorem asserts the existence of a covering graph in
which large amounts of any winding and unwinding which occurs in
the iterates of $\rho$ lift to different sheets.  We then proceed
to prove the Apt Immersion Theorem in the linear case
($\S$\ref{LinearCaseChapter}) and the non-linear case
($\S$\ref{PathUnitsChapter}).  An index of notation and
terminology is included at the back of the paper for the
convenience of the reader.

\section{The growth of an automorphism} \label{GrowthChapter}

\begin{defn}\label{DefnEquiv}\index{saaequivalence@$\simeq$}
We define a relation $\preceq$ on the set of all functions $\Nat
\to \Nat$ by writing $f \preceq g$ if there exist constants $A, B,
D, E
> 0$ and $C \geq 0$ such that \begin{equation*}f(n) \leq Ag(Bn+ C) + Dn + E,\end{equation*} for all  $n \in
\Nat$. Two functions $f,g:\Nat \to \Nat$ are said to be $\simeq$
equivalent if $f \preceq g$ and $g \preceq f$.
\end{defn}

It is easily verified that $\simeq$ is an equivalence relation.
Denote by $p_{\expgrowth}$\index{spexponential@$p_{\expgrowth}$}
the $\simeq$-equivalence class which contains all functions
bounded below by a function $k \mapsto c^k$, for some constant $c
> 1$; denote by $p_1$\index{spi@$p_i$} the $\simeq$-equivalence class which
contains all functions bounded above by a polynomial function of
degree 1; and, for each integer $d \geq 2$, denote by $p_d$ the
$\simeq$-equivalence class which contains the function $k \mapsto
k^d$. The classes $p_{\expgrowth}, p_1, p_2, \dots$ are pairwise
disjoint. For each $d \in \Nat$, we write $f \preceq p_d$ if $f
\preceq (k \mapsto k^d)$. We say that a function $f$ has
\emph{degree}\index{wdegree@degree!of a function} $d$ if $f \in
p_d$ for some $d \in \Nat$, and we say that $f$ is
\emph{linear}\index{wlinearfunction@linear function} if $f \in
p_1$ and $f$ is unbounded.

\begin{notation}
For a group $G$ generated by a finite subset $A \subset G$ and for
each element $w \in G$, write $\abs{w}_A$ for the distance from
the identity of $G$ to $w$ in the word-metric on $G$ with respect
to $A$.
\end{notation}

\begin{defn}[Growth of an automorphism]
\index{wgrowthofanautomorphism@growth of an automorphism}
\index{sgrowthofanautomorphism@$\mathcal{G}_{\phi}$} For each
$\phi \in \Aut(G)$, define $\norm{\phi}_A := \max
\{\abs{\phi(a)}_A \hbox{ } | \hbox{ } a \in A \}$; and define a
function $\mathcal{G}_{\phi, \, A}: \Nat \to \Nat$, called the
\emph{growth (function) of $\phi$ (with respect to $A$)}, by
$\mathcal{G}_{\phi, \, A}(n) = \norm{\phi^n}_A$.
\end{defn}

The following elementary properties of the growth function are
easily verified.

\begin{prop}[Properties of the growth function]\label{PropertiesOfGrowthOfAutomorphism}
Let $G$ be a finitely generated group, let $A \subset G$ be a
finite generating set and let $\phi \in \Aut(G)$ be an
automorphism. The following properties hold:
\begin{enumerate}
 \item [(G1)]
 for each finite generating set $B \subset G$, $\mathcal{G}_{\phi, \, A} \simeq \mathcal{G}_{\phi, \,
 B}$;
 \item [(G2)]
 for each $k \in \Nat$, $\mathcal{G}_{\phi, \, A} \simeq \mathcal{G}_{\phi^k, \,
 A}$;
 \item [(G3)]
 for each $\phi$-invariant subgroup $S \leq G$ of finite index with finite generating set $B$, $\mathcal{G}_{\phi, \, A} \simeq \mathcal{G}_{\phi|_S, \, B}$.
\end{enumerate}
\end{prop}


\begin{notation}\index{sgrowthofanautomorphismabelian@$\mathcal{G}^{\ab}_{\phi}$}
Empowered by Property (G1), we usually omit mention of $A$ from
the notation, writing simply $\mathcal{G}_{\phi}$. Also, as
mentioned in the introduction, we usually write
$\mathcal{G}^{\ab}_{\phi}$ for $\mathcal{G}_{\phi^{\ab}}$.
\end{notation}

\begin{rem}[Other notions of automorphism growth]\label{OtherNotionsOfGrowth}
For alternative notions of the growth of an automorphism the
reader is referred to \cite{PolynomialDehnFunctions}, where
Bridson lists four distinct notions of the growth of an
automorphism and sketches the relationship between them in the
case that $G$ is a finitely generated abelian or non-abelian free
group.  In Bridson's notation, the growth function $G_\phi$ is
written $g_{0, \, \phi}$.
\end{rem}

\section{Graphs and covering graphs}\label{GraphChapter}

The notation for undirected graphs that we use is mostly that of
Stallings \cite{Stallings}.  For the convenience of the reader we
describe this notation below, with some additions, before
introducing the simple notion of an end-pointed graph and some
related constructions in $\S$\ref{EndpointedGraphsSection} and
homotopy equivalences of graphs in $\S$\ref{HomotopyEquivSection}.

\subsection{Graphs}\label{GraphsSection}

\begin{defn}[Graph]\index{wgraph@graph} \index{siota@$\iota$}
A \emph{graph} $G$ consists of sets $\mathcal{E}_G$ and
$\mathcal{V}_G$ and functions $r_G: \mathcal{E} \to \mathcal{E}$
and $\iota_G:\mathcal{E} \to \mathcal{V}$ subject to the
conditions that $r_G \circ r_G (e) = e$ and $r_G(e) \not = e$ for
each $e \in \mathcal{E}$.  We write $G = (\mathcal{V}_G,
\mathcal{E}_G, r_G, \iota_G)$.
\end{defn}
For brevity, we often omit the sets and functions from the
notation, stating simply that $G$ is a graph; we write
$\rev{e}$\index{saareverse@$\rev{\;}\;$}\index{seoverline@$\rev{e}\;$}
for $r(e)$; we define a third map $\tau: \mathcal{E} \to
\mathcal{V}$\index{stau@$\tau$} such that, for each $e \in
\mathcal{E}$, $\tau(e) = \iota(\rev{e})$; and we omit the
subscript from $\iota_G$ unless it is necessary to avoid
ambiguity. We call $\mathcal{V}_G$ the set of \emph{vertices} and
$\mathcal{E}_G$ the set of \emph{(directed)
edges}\index{wdirectededge@directed edge}. For an edge $e \in
\mathcal{E}_G$, we call $\iota(e)$ the \emph{initial point of
$e$}\index{winitialpoint@initial point!winitialpointofanedge@of an
edge}, $\tau(e)$ the \emph{terminal point of
$e$}\index{wterminalpoint@terminal point!wterminalpointofanedge@of
an edge} and $\rev{e}$ the \emph{reverse of
$e$}\index{wreverse@reverse!wrofanedge@of an edge}. A pair $\{e,
\rev{e}\}$ is called a \emph{geometric (or undirected)
edge}\index{wgeometricedge@geometric edge}. An
\emph{orientation}\index{worientation@orientation} $\Orient$ of
$G$ is a set containing exactly one directed edge from each
geometric edge.

\begin{notation}
In general, directed edges will be denoted by lower case letters
and the geometric edge containing a particular directed edge will
be denoted by the corresponding upper case letter.
\end{notation}

Of course, graphs may be considered to be topological objects as
well as combinatorial ones. In general, we will not distinguish
between a graph $G$ and the following geometric realisation:
Realise $G$ as a CW-complex with one 0-cell for each element of
$\mathcal{V}_G$, one 1-cell for each geometric edge and attaching
maps as specified by $\iota$. Define a path-metric on $G$ by
assigning unit length to each 1-cell.  Unless otherwise specified,
we will consider only connected graphs.

If $H = (\mathcal{V}_{H}, \mathcal{E}_{H})$ and $G =
(\mathcal{V}_{G}, \mathcal{E}_{G})$ are graphs, a morphism of
graphs\index{wmorphismofgraphs@morphism of graphs} $p : H \to G$
consists of a pair of functions, $p_\mathcal{V}: \mathcal{V}_{H}
\to \mathcal{V}_{G}$ and $p_\mathcal{E}: \mathcal{E}_{H} \to
\mathcal{E}_{G}$ subject to the conditions that, for each $e \in
\mathcal{E}_H$, $p_\mathcal{V} \circ \iota(e) = \iota \circ
p_\mathcal{E}(e)$ and $p_\mathcal{E}(\rev{e}) =
\rev{p_\mathcal{E}(e)}$. We write $p = (p_\mathcal{V},
p_\mathcal{E})$ and often abuse notation by writing $p$ for both
$p_\mathcal{V}$ and $p_\mathcal{E}$. If $p: H \to G$ is a morphism
of graphs we say that $(H, p)$ is a \emph{$G$-labelled
graph}\index{wGalabelledgraph@$G$-labelled-graph}. The morphism
$p$ is called the \emph{labelling
map}\index{wlabellingmap@labelling map} and, for each edge $e \in
\mathcal{E}_H$, $p(e)$ is called the
\emph{label}\index{wlabel@label on an edge} on $e$. We often omit
mention of the map $p$ if it may be understood from the context;
we say simply that $H$ is a $G$-labelled graph and write
$\hat{e}$\index{saahat@$\hat{\,}$} for $p(e)$. We say that two
$G$-labelled graphs $(H_1, p_1)$ and $(H_2, p_2)$ are
\emph{$G$-labelled-graph-isomorphic}\index{wGalabelledgraphisomorphic@$G$-labelled-graph-isomorphic}
if there is a graph isomorphism $f:H_1 \to H_2$ such that $p_1 =
p_2 \circ f$.

\begin{rem}[`Drawing' $G$-labelled graphs]
Given a graph $G$, we may describe a $G$-labelled graph $(H, p)$
in the following way: consider a third graph $\Sigma$, an
orientation $\Orient_\Sigma$ and a set of paths in $G$ (see
$\S$\ref{PathsSection}) which label the edges of $\Orient_\Sigma$
subject to the condition that, if $e$ and $e'$ are directed edges
in $\Orient_\Sigma$ with labels $\alpha$ and $\alpha'$
respectively and such that $\iota(e) = \iota(e')$, then
$\iota(\alpha) = \iota(\alpha')$. The graph $H$ is the subdivision
of $\Sigma$ such that each directed edge $e \in \Orient_\Sigma$
labelled by a path $\rho = d_1 d_2 \dots d_n$ in $G$ corresponds
to a sequence of $n$ distinct directed edges $e_1, e_2, \dots,
e_n$ in $\mathcal{E}_H$; define $p(e_i) = d_i$. This completely
determines the map $p$.
\end{rem}

The \emph{star of $v$ (in $G$)} is $St(v, G) := \{e \in
\mathcal{E}_{G} \hbox{ } | \hbox{ } \iota(e) =
v\}$\index{sstar@$St$}; thus $\abs{St(v, G)}$ is the \emph{valence
of $v$ (in $G$)}.  A graph for which each vertex has valence at
least two is said to be \emph{minimal}\index{wminimalgraph@minimal
graph}. A morphism of graphs $p: H \to G$ induces a map $p_v:
St(v, H) \to St(p(v), G)$ for each $v \in \mathcal{V}_{H}$. If
$p_v$ is injective for each $v \in \mathcal{V}_{H}$, we say that
$p$ is an \emph{immersion} and that $(H, p)$ is a
$G$-\emph{immersion}\index{wGaimmersion@$G$-immersion}.
If $p_v$ is bijective for each $v \in \mathcal{V}_{H}$, we say
that $p$ is a covering map and that $(H, p)$ is a
$G$-\emph{cover}\index{wGaCover@$G$-cover}\footnote{The definition
of a $G$-cover above is equivalent to the usual topological
definition of a covering of a graph $G$.}. For brevity, we usually
omit mention of the map from $G$-immersions and $G$-coverings,
that is, we say that $H$ is a $G$-immersion or a $G$-covering. For
a $G$-covering $(H, p)$ and vertices $v, w \in G$, the sets
$p^{-1}(v)$ and $p^{-1}(w)$ have the same cardinality $s$; we say
that $H$ is an \emph{$s$-sheeted $G$-cover}.

\begin{rem}\label{PerspectiveOnImmersions}
Let $G = (\mathcal{V}_G, \mathcal{E}_G)$ be a finite graph.
Stallings \cite{Stallings} observed that a finite $G$-immersion
may be identified with a graph $J$ constructed as follows:  For
each $v \in \mathcal{V}_G$, choose an integer $s_v \geq 0$ and
define
\begin{equation*}\mathcal{V}_J := \underset{v \in \mathcal{V}_G} {\amalg}  \{(v, i) \hbox{ } |
\hbox{ } 1 \leq i \leq s_v \}.\end{equation*}  For each edge $e
\in \mathcal{E}_G$, choose an integer $s_e$ such that $0 \leq s_e
\leq \min\{s_{\iota(e)}, s_{\tau(e)}\}$ and $s_e = s_{\rev{e}}$.
Define
\begin{equation*}\mathcal{E}_J := \underset {e \in \mathcal{E}_G} {\amalg} \{(e, i) \hbox{ } |
\hbox{ } 1 \leq i \leq s_e \}.\end{equation*} Choose a map
$\iota:\mathcal{E}_J \to \mathcal{V}_J$ such that the restriction
of $\iota$ to each set $\{(e, i) \hbox{ } | \hbox{ } 1 \leq i \leq
s_e\}$ is an injection into $\{(\iota(e), i) \hbox{ } | \hbox{ } 1
\leq i \leq s_{\iota(e)} \}$. Finally, define $\rev{(e, i)} =
(\rev{e}, i)$ for each $e \in \mathcal{E}_J$, define $p:J \to G$
such that $p((v, i)) = v$ and $p((e, i)) = e$ for each $v \in
\mathcal{V}_G$ and each $e \in \mathcal{E}_G$. Such a
$G$-immersion $J$ is an $s$-sheeted $G$-covering if and only if
$s_v = s_e = s$ for each $v \in \mathcal{V}_G$ and each $e \in
\mathcal{E}_G$.
\end{rem}

\begin{defn}
Let $G$ be a graph.  A \emph{handle}\index{whandle@handle} in $G$
is a maximal subgraph $H$ such that $H$ is a non-trivial
line-segment, the ends of $H$ have valence at least three in $G$
and the remaining vertices of $H$ have valence two in $G$.
\end{defn}

\begin{notation}\label{SetminusNotation}
For a graph $G$ and a subgraph $S$, we write $G \setminus
S$\index{saasetminus@$\setminus$} for the subgraph of $G$ which is
the topological closure of the vertex set $\mathcal{V}_G \setminus
\mathcal{V}_S$ and edge set $\mathcal{E}_G \setminus
\mathcal{E}_S$.
\end{notation}

\subsection{Paths and circuits in graphs}\label{PathsSection}

Let $G$ be a graph.  A \emph{path}\index{wpath@path} $\rho$ in $G$
is either a vertex $v \in \mathcal{V}$ (the \emph{trivial path at
$v$}) or a non-empty finite ordered list of (directed) edges $d_1,
d_2, \dots, d_s \in \mathcal{E}_G$ such that $\tau(d_i) =
\iota(d_{i+1})$ for $1 \leq i < s$ (we usually omit commas in the
list of edges). If $\rho$ is the trivial path at $v$ we write
$\iota(\rho) = \tau(\rho) = v$ and $l(\rho) = 0$, otherwise, we
write $\iota(\rho)$ for $\iota(d_1)$, $\tau(\rho)$ for $\tau(d_s)$
and $l(\rho)$\index{slength@$l$} for $s$ (the \emph{length} of
$\rho$). A \emph{tight
path}\index{wtight@tight!w_tightpath@path}\index{wpath@path!wtightsub@tight}
in $G$ is either a trivial path or a non-trivial path for which
the corresponding finite list of edges is reduced (that is,
$d_{i+1} \neq \rev{d}_i$ for each $i = 1, 2, \dots, s-1$). We say
that a path $\rho$ is \emph{closed (at $v$)} if $\iota(\rho) =
\tau(\rho) = v$. If $\rho$ is the trivial path at $v$ and $n \in
\Nat$, we write $\rho^n$ for trivial path at $v$. If $\rho$ is a
non-trivial path in $G$ and $n \in \Nat$, say $\rho = d_1 d_2
\dots d_s$, we write $\rho^n$ for the closed path with directed
edge list $d_1 d_2 \dots d_n$ repeated $n$ times. Also, for a
closed path $\rho$ in $G$ we write
$l^{\circlearrowleft}(\rho)$\index{slengthcircle@$l^{\circlearrowleft}$}
for the length of the cyclically reduced path corresponding to
$\rho$. A \emph{circuit}\index{wcircuit@circuit} in $G$ is an
equivalence class of closed paths in $G$ under the relation of
cyclic permutation of the list of edges. A \emph{tight circuit}
\index{wtight@tight!w_tightcircuit@circuit} in $G$ is an
equivalence class of closed tight and cyclically reduced paths in
$G$ under the same relation. The map $l$ extends naturally to
circuits.  For a path $\rho$, we write
$[\rho]$\index{saasquarebracket@$[\,]$} for the tight path
obtained by reducing $\rho$;
for a circuit $\sigma$, we write $[\sigma]$ for the circuit
obtained by reducing and cyclically reducing $\sigma$.

\begin{rem}
Switching to the topological perspective, the map $[ \, ]$ from
paths in $G$ to tight paths in $G$ corresponds to tightening
relative to the end-points. Similarly, the map $[ \, ]$ from
circuits in $G$ to tight circuits in $G$ corresponds to
tightening.
\end{rem}

For a non-trivial path $\rho = d_1 d_2 \dots d_s$ in $G$ and a
geometric edge $E$ in $G$, we say that $\rho$
\emph{crosses}\index{wcrosses@crosses} $E$ if $\rho$ is
non-trivial and either $e$ or $\rev{e}$ appear in the list of
edges defining $\rho$. Let $l^{\ab}$ denote the $l^1$ norm on the
cellular chain complex of $G$; equivalently, for an orientation
$\Orient = \{e_i \, | \, 1 \leq i \leq r \}$ of $G$,
\begin{equation*}l^{\ab}(\rho) = \sum_{1 \leq i \leq r} \abs{c_i},\end{equation*}\index{slengthabelian@$l^{\ab}$}
where $c_i := \abs{\{j \, | \, d_j = e_i\}} - \abs{\{j \, | \, d_j
= \rev{e}_i\}}.$

We say that a closed tight path $\rho$ in $G$ is
\emph{primitive}\index{wprimitivepath@primitive
path}\index{wpath@path!wprimitivesub@primitive} if there is no
closed tight path $\mu$ in $G$ and integer $m \geq 2$ such that
$\rho = \mu^m$. If $\rho$ is a primitive closed tight path and
$\delta = \rho^n$ for some positive integer $n$, we say that
$\rho$ is a \emph{primitive closed tight path corresponding to
$\delta$}. The following lemma is easily verified.

\begin{lem}\label{UniqueMinimalPaths}
For each closed tight path $\rho$ in $G$ there is a unique
primitive closed tight path corresponding to $\rho$.
\end{lem}

\begin{defn}\label{EndsAndEndPaths}
For a graph $G$, an \emph{end}\index{wend@end} is a vertex $v \in
\mathcal{V}$ with valence one.  An
\emph{end-path}\index{wendpath@end-path} is a non-trivial tight
path $\rho = d_1 d_2 \dots d_s$ in $G$ such that $\iota(d_1)$ is
an end, each $\tau(d_i)$ has valence two for $i = 1, 2, \dots,
s-1$, and $\tau(d_s)$ has valence not equal to two.
\end{defn}

\subsection{Stallings' folding operation and Stallings'
algorithm}\label{FoldingSection}


\begin{defn}[Stallings' Folding Operation]\index{wStallings@Stallings!wStallingsFoldingOperation@Folding Operation}
Let $G$ be a graph and let $H$ be a $G$-labelled graph which is
{\bf not} a $G$-immersion.  There exist a vertex $v_0 \in
\mathcal{V}_H$ and distinct edges $d_1, d_2 \in \mathcal{E}_H$
such that $\iota(d_1) = \iota(d_2) = v_0$ and $\hat{d}_1 =
\hat{d}_2$ (that is, the label on $d_1$ and $d_2$ is the same).
Define a $G$-labelled graph $H'$ as follows: $\mathcal{V}_{H'}$ is
defined from $\mathcal{V}_H$ by identifying $\tau(d_1)$ and
$\tau(d_2)$ (unless they are already equal); $\mathcal{E}_{H'}$ is
defined from $\mathcal{E}_H$ by identifying $d_1$ and $d_2$ and
identifying $\rev{d}_1$ and $\rev{d}_2$; let $f_V:\mathcal{V}_H
\to \mathcal{V}_{H'}$ and $f_E:\mathcal{E}_H \to \mathcal{E}_{H'}$
denote the natural maps and define $r$ and $\iota$ to be the
unique maps such that $f = (f_\mathcal{V}, f_{\mathcal{E}}): H \to
H'$ is a morphism of graphs.
The morphism $f$ is said to be a \emph{folding
morphism}\index{wfolding@folding |see{Stallings}}. Let $H_0$ and
$H_n$ be $G$-labelled graphs for some $n \in \Nat$. We say that
$H_0$ \emph{folds to} $H_n$ if there exist $G$-labelled graphs
$H_1, H_2, \dots H_{n-1}$ and folding morphisms $f_1, f_2, \dots
f_n$ such that $f_i : H_{i-1} \to H_i$.
\end{defn}

\begin{thm}[Stallings]\label{TheFoldingTheorem}
Let $G$ be a graph.  For each finite $G$-labelled graph $H$ there
is a unique $G$-immersion $H'$ (called the \emph{$G$-immersion
determined by $H$}) such that $H$ folds to $H'$.
\end{thm}

\begin{rem}\label{OrderOfFoldingUnimportant}
We may find the unique $G$-immersion $H'$ by following a simple
algorithm: Define $H_0 = H$.  Inductively, for each integer $i
\geq 0$, if $H_i$ is a $G$-immersion then set $H' = H_i$ and
terminate the algorithm, otherwise there exists a $G$-labelled
graph $H_{i+1}$ such that $H_i$ folds to $H_{i+1}$. Because $H$ is
finite and $\abs{\mathcal{E}_{H_{i+1}}} =
\abs{\mathcal{E}_{H_i}}-1$, the algorithm terminates in a most
$\abs{\mathcal{E}_H}-1$ steps. Theorem \ref{TheFoldingTheorem}
informs us that our choice of $H_{i+1}$ at each stage is
unimportant.
\end{rem}

Let $G$ be a graph, let $(H, p)$ be a $G$-labelled graph, fix
vertices $v \in \mathcal{V}_G$ and $w \in \mathcal{V}_H$ such that
$p(w) = v$ and consider the induced homomorphism $p_{\ast}:
\pi_1(H, w) \to \pi_1(G, v)$.

\begin{thm}[Stallings]
If $H$ is a $G$-immersion then $p_*: \pi_1(H, w) \to \pi_1(G, v)$
is injective.
\end{thm}

\begin{thm}[Stallings]
Let $G$ be a graph and $v \in G$ a vertex, let $(H_1, p_1)$ and
$(H_2, p_2)$ be $G$-labelled graphs such that $H_1$ folds to
$H_2$, let $v_1 \in H_1$ be a vertex and $v_2$ the corresponding
vertex in $H_2$. Then ${p_1}_*(\pi_1(H_1, v_1)) =
{p_2}_*(\pi_1(H_2, v_2))$.
\end{thm}

The following theorem is a slight generalisation of Theorem 6.1
\cite{Stallings} (because we allow $\abs{\mathcal{V}_G} > 1$). The
proof below is that of Stallings, which we include because of its
fundamental importance to this paper.

\begin{thm}[Stallings'
Algorithm]\label{ExtendingImmersionsToCoverings}\index{wStallings@Stallings!wStallingsAlgorithm@
Algorithm} Let $G$ be a finite graph and $H$ a finite
$G$-immersion. There exists a finite $G$-covering $\tilde{H}$ such
that $H$ is $G$-labelled graph isomorphic to a subgraph of
$\tilde{H}$.
\end{thm}

\begin{proof}
By relabelling the vertices and edges of $H$ (if necessary) we may
assume that $H$ is constructed as in Remark
\ref{PerspectiveOnImmersions} and we assume the corresponding
notation. Define $s := \max\{s_v \hbox{ } | \hbox{ } v \in
\mathcal{V}_G\}$, $\mathcal{V}_{\tilde{H}} := \mathcal{V}_G \times
\{1, 2, \dots, s\}$ and $\mathcal{E}_{\tilde{H}} := \mathcal{E}_G
\times \{1, 2, \dots, s\}$. Define $\iota_{\tilde{H}}$ such that
the restriction of $\iota_{\tilde{H}}$ to each set $\{(e, i)
\hbox{ } | \hbox{ } 1 \leq i \leq s\}$ is a bijection into
$\{(\iota_G(e), i) \hbox{ } | \hbox{ } 1 \leq i \leq s\}$ and, for
each $e \in \mathcal{E}_{\tilde{H}}$, the restriction of
$\iota_{\tilde{H}}$ to the set $\{(e, i) \hbox{ } | \hbox{ } 1
\leq i \leq s_e\}$ corresponds to $\iota_H$ (this is possible
since the restriction of $\iota_H$ to the set $\{(e, i) \hbox{ } |
\hbox{ } 1 \leq i \leq s_e\}$ is injective). Finally, define
$\rev{(e, i)} = (\rev{e}, i)$ for each $e \in
\mathcal{E}_{\tilde{H}}$, define $p:\tilde{H} \to G$ such that
$p((v, i)) = v$ and $p((e, i)) = e$ for each $v \in \mathcal{V}_G$
and each $e \in \mathcal{E}_G$.  It follows from Remark
\ref{PerspectiveOnImmersions} that $\tilde{H}$ satisfies the
conclusions of theorem.
\end{proof}

\begin{rem}
It is clear from the proof above that
$\abs{\mathcal{V}_{\tilde{H}}} \leq \abs{\mathcal{V}_H} .
\abs{\mathcal{V}_G}$.
\end{rem}

\subsection{End-pointed graphs}\label{EndpointedGraphsSection}

Assigning end-points and base-points to graphs allows us to
discuss combining graphs and the movement of a path through a
graph in a natural way.

\begin{defn}[A vocabulary for end-pointed graphs]\index{wendpointedgraph@end-pointed
graph}\index{wbasepointedgraph@base-pointed
graph}\index{winitialpoint@initial point!winitialpointofagraph@of
a graph}\index{wterminalpoint@terminal
point!wterminalpointofagraph@of a graph} An \emph{end-pointed
graph} is simply a graph $H$ with two distinguished vertices
called the \emph{initial point of $H$}, denoted $\iota(H)$, and
the \emph{terminal point of $H$}, denoted $\tau(H)$.  We refer to
$\iota(H)$ and $\tau(H)$ collectively as the \emph{end-points of
$H$}.  A path \emph{across} $H$ is a non-trivial path $\rho$ such
that $\iota(\rho) = \iota(H)$ and $\tau(\rho) = \tau(H)$.  A
\emph{base-pointed graph} is an end-pointed graph $H$ for which
$\iota(H) = \tau(H)$, in which case we call $\iota(H)$ the
\emph{base-point} of $H$.
\end{defn}

\begin{rem}
The end-points of a graph are not necessarily ends in the sense of
Definition \ref{EndsAndEndPaths}.
\end{rem}

\begin{notation}
When depicting an end-pointed graph it will be our convention to
denote the initial point by a square, the terminal point by an
asterisk and all other vertices by circles (see, for example,
Figure \ref{EndpointedGraphsFigure}).
\end{notation}

\begin{defn}[More vocabulary for end-pointed graphs]
Let $H$ be an end-pointed $G$-labelled graph for some graph $G$.
The \emph{end-pointed $G$-immersion determined by
$H$}\index{wendpointeddeterminedby@end-pointed $G$-immersion
determined by $H$} is the $G$-immersion determined by $H$ with
initial point the natural image of $\iota(H)$ and terminal point
the natural image of $\tau(H)$; the \emph{base-pointed
$G$-labelled graph determined by
$H$}\index{wbasepointedlabelledgraphdeterminedby@base-pointed
$G$-labelled graph determined by $H$} is the $G$-labelled graph
$H'$ obtained from $H$ by identifying $\iota(H)$ and $\tau(H)$ and
defining $\iota(H')$ to be the natural image of $\iota(H)$;
finally, the \emph{base-pointed $G$-immersion determined by
$H$}\index{wbasepointedimmersiondeterminedby@base-pointed
$G$-immersion determined by $H$} is the $G$-immersion $H''$
determined by $H'$ with $\iota(H'')$ defined to be the natural
image of $\iota(H)$.
\end{defn}

\begin{con}[Lines and circles]\label{CirclesAndLines}
Let $G$ be a graph and let $\rho$ be a path in $G$. Define an
end-pointed $G$-labelled graph, denoted $L(\rho)$\index{sL@$L$},
in the following way: $L(\rho)$ is an interval subdivided into
$l(\rho)$ edges; specify one end of the graph as the initial
point, the other end as the terminal point and assign labels such
that the unique tight path across $L(\rho)$ is labelled by $\rho$.
Further, denote by $C(\rho)$\index{sC@$C$} the base-pointed
$G$-labelled graph determined by $L(\rho)$ (see Figure
\ref{EndpointedGraphsFigure}). Observe that if $\rho$ is a
cyclically reduced path then $C(\rho)$ is a base-pointed
$G$-immersion.
\end{con}

\setcounter{figure}{\value{thm}} \stepcounter{thm}
\begin{figure}
\includegraphics{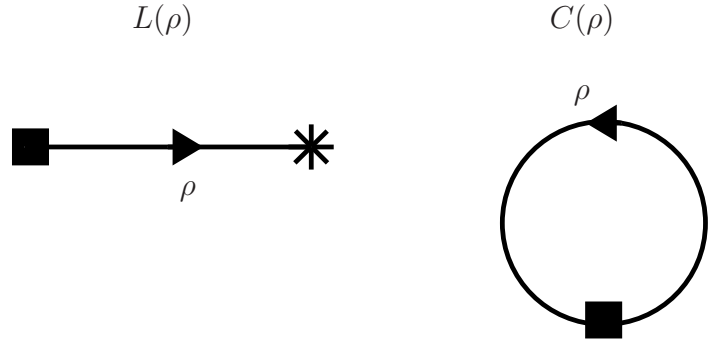}
\caption{$L(\rho)$ and $C(\rho)$. \label{EndpointedGraphsFigure}}
\end{figure}

\begin{con}[Combining end-pointed graphs]\index{saaVeeopen@$\combine{\,}$}\index{saaVeesquare@$\combineandfold{\,}$}\index{saaVeesquarewithC@$\combineandfoldC{\,}$} Let $H_1, H_2,
\dots, H_s$ be end-pointed $G$-labelled graphs for some graph $G$.
Define an end-pointed $G$-labelled graph,
\begin{equation*}\combine{H_1, H_2, \dots, H_s} := H_1 \amalg H_2 \amalg
\dots \amalg H_s / \sim,\end{equation*} where $\sim$ identifies
$\tau(H_i)$ and $\iota(H_{i+1})$ for each $i = 1, 2, \dots, s-1$.
Define the initial point of $\combine{H_1, H_2, \dots, H_s}$ to be
the natural image of $\iota(H_1)$ and the terminal point to be the
natural image of $\tau(H_s)$ (see Figure
\ref{ConstructionVee1Figure}). Further, let $\combineandfold{H_1,
H_2, \dots, H_s}$ (respectively $\combineC{H_1, H_2, \dots, H_s}$,
$\combineandfoldC{H_1, H_2, \dots, H_s}$) denote the end-pointed
$G$-immersion (respectively base-pointed $G$-labelled graph,
base-pointed $G$-immersion) determined by $\combine{H_1, H_2,
\dots, H_s}$  (see Figure \ref{ConstructionVee2Figure}).
\end{con}

\setcounter{figure}{\value{thm}} \stepcounter{thm}
\begin{figure}
\includegraphics{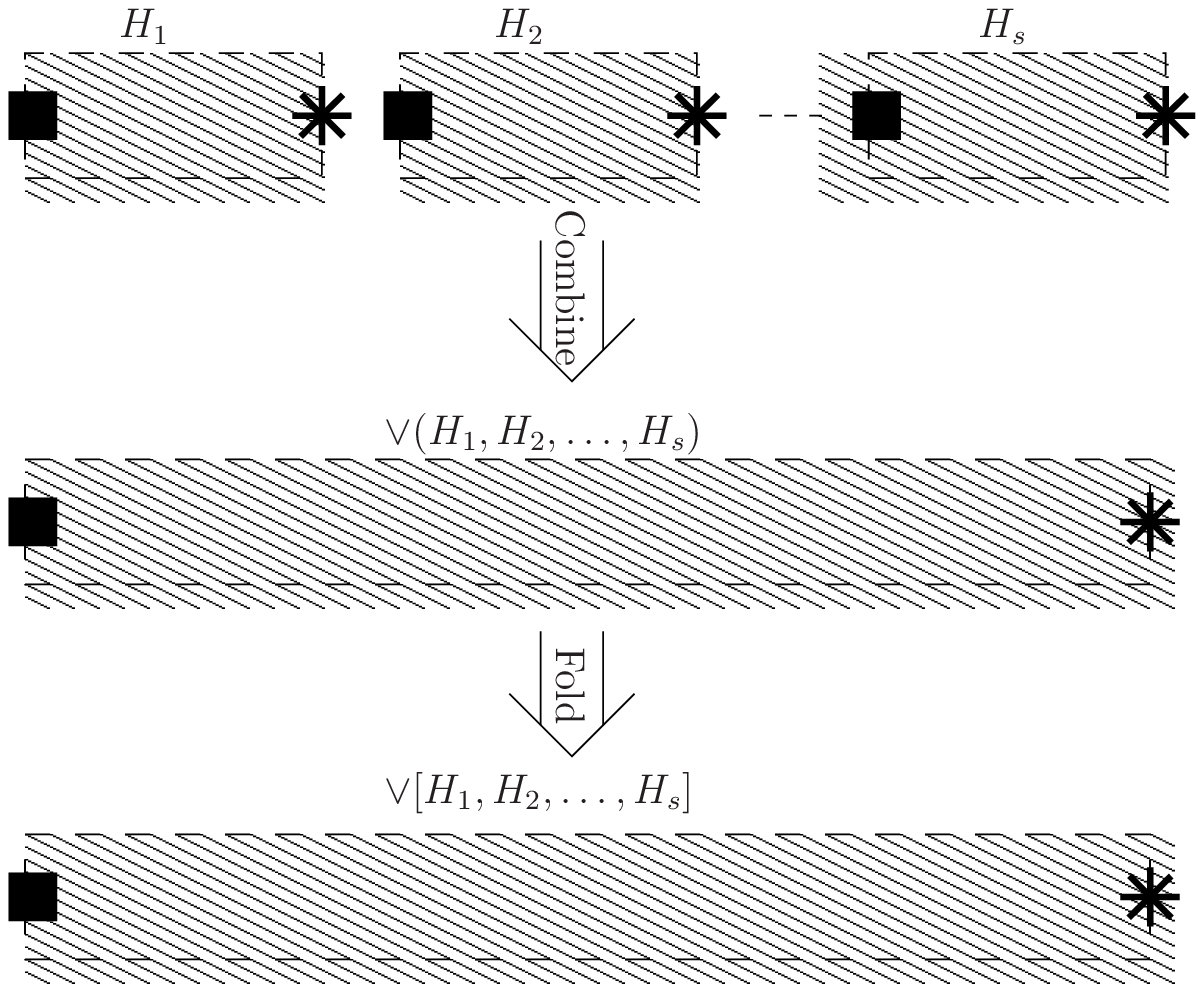}
\caption{A schematic depiction of the construction of
$\combine{H_1, H_2, \dots, H_s}$ and $\combineandfold{H_1, H_2,
\dots, H_s}$. \label{ConstructionVee1Figure}}
\end{figure}

\setcounter{figure}{\value{thm}} \stepcounter{thm}
\begin{figure}
\includegraphics{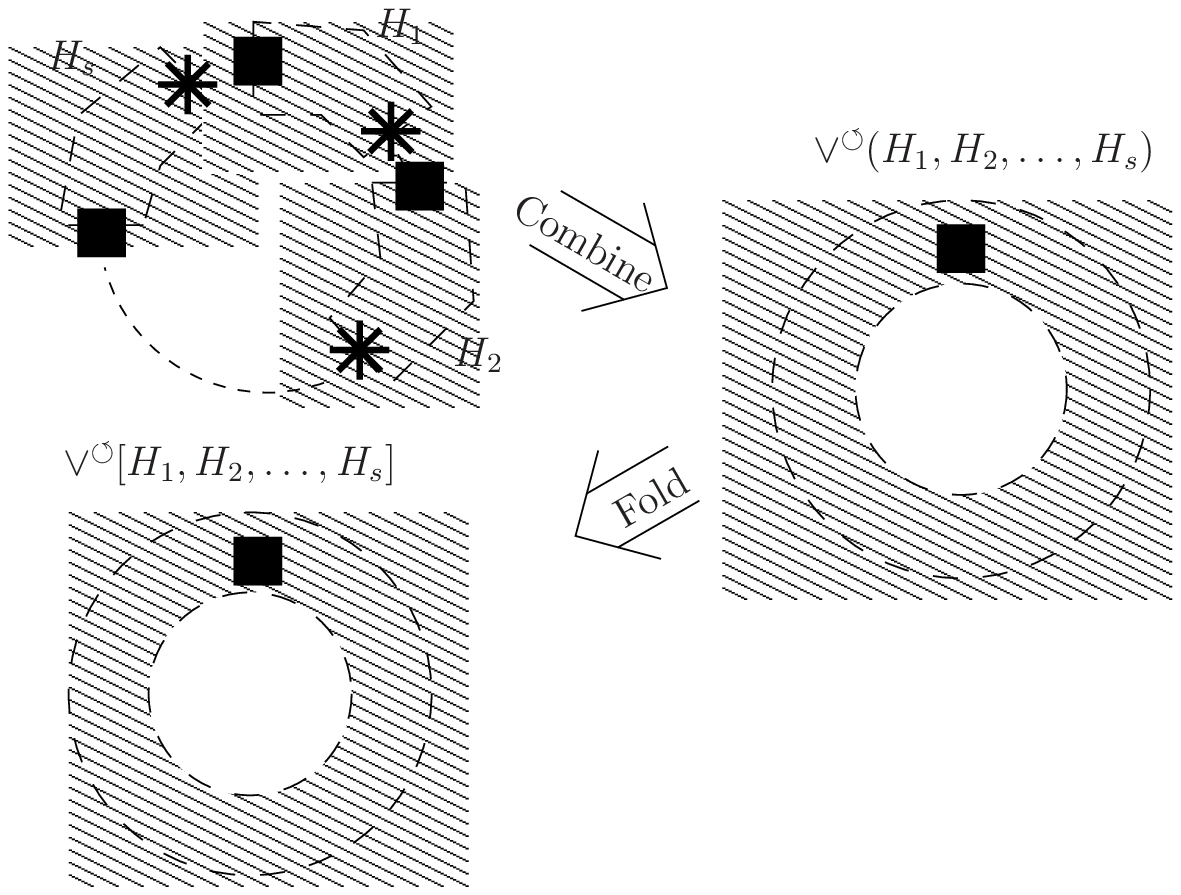}
\caption{A schematic depiction of the construction of
$\combineC{H_1, H_2, \dots, H_s}$ and $\combineandfoldC{H_1, H_2,
\dots, H_s}$. \label{ConstructionVee2Figure}}
\end{figure}

\subsection{Homotopy equivalences of graphs}\label{HomotopyEquivSection}

Since graphs may be thought of as topological objects, we may
consider homotopy equivalences of graphs.  Let $G$ and $H$ be
graphs.  For technical reasons, we consider only those homotopy
equivalences\index{whomotopyequivalence@homotopy equivalence of
graphs} $f: H \to G$ with the properties that $f:\mathcal{V}_H \to
\mathcal{V}_G$ and, for each $e \in \mathcal{E}_H$, $f(e)$ is a
tight path in $G$. Such a homotopy equivalence induces a map (also
denoted by $f$) from the set of paths in $H$ to the set of paths
in $G$. Denote by $f_{\#}$ the map from the set of (tight) paths
in $H$ to the set of tight paths in $G$ defined by $\rho \mapsto
[f(\rho)]$.

\section{Improved relative train track representatives of automorphisms in $\PG (F)$}\label{IRTTChapter}

In the first two sections of this chapter we give an exposition of
those parts of the theory of improved relative train track
representatives necessary for the work that follows. It is
included for the convenience of the reader and is, by necessity,
brief and far from comprehensive. In particular, we discuss only
the $\PG(F)$ case of Bestvina, Feighn and Handel's Improved
Relative Train Track Theorem (Theorem 5.1.5 \cite{BFH1}). Although
this reduces the scope of the theorem significantly, it simplifies
the statement and allows us to state the $\PG(F)$ case with a
combinatorial flavour rather than a topological one.  The reader
is referred to the following series of papers for a full
exposition of this powerful theory: \cite{TrainTracks},
\cite{BFH1}, \cite{BFH2}. Most of the notation used below is that
introduced by Bestvina, Feighn and Handel.  In
$\S$\ref{IRTTGrowthSection} the relationship between the growth of
an automorphism and the growth of tight paths and tight circuits
in an IRTT representative is developed.  Finally, in
$\S$\ref{Proofsofgrowthpropertiessection} we provide proofs of
theorems \ref{AutomorphismGrowthTheorem} and
\ref{InverseHasSameGrowth} and we complete the discussion of
Remark \ref{AlgorithmicProperties}. The important notion of the
reverse of an IRTT representative is also defined.

\subsection{Topological representatives of free group automorphisms}

Fix $n \in \Nat$, let $F$ denote the free group of rank
$n$\index{sn@$n$}, let $R$\index{sR@$R$} denote the graph with one
vertex $b_R$ and $n$ geometric edges, and fix an identification of
$F$ with $\pi_1(R, b_R)$ by identifying a generating set of $F$
with a generating set of $\pi_1(R, b_R)$.

\begin{defn}\label{MarkedGraphDefn}
A \emph{marked graph}\footnote{Culler and Vogtmann
\cite{ModuliOfGraphs}  call this a \emph{marking} on $G$, and save
the term \emph{marked graph} for an equivalence class of markings
under a suitable equivalence relation.  For our purposes it is
enough to consider individual markings, and we will follow the
notation of \cite{TrainTracks} by using the definition of a marked
graph given in the text.}\index{wmarkedgraph@marked graph} is a
pair $(G, m)$ for which:
\begin{enumerate}
  \item $G$ is a non-trivial minimal graph with fundamental group isomorphic to $F$;
  \item $m$ is a homotopy equivalence $R \to G$.
\end{enumerate}
The map $m$ is called a \emph{marking}\index{wmarking@marking} on
$G$ and $b = m(b_R) \in G$ is called the \emph{base-point} of the
marked graph.
\end{defn}

A marking $m$ determines an identification between $\pi_1(G, b)$
and $F$.  Tight circuits in $G$ are in a one-to-one correspondence
with the conjugacy classes of $\pi_1(G, b)$ and hence with the
conjugacy classes of $F$. A marked graph $(G, m)$ and a homotopy
equivalence $f : G \to G$ which fixes $b$ determine an
automorphism of $F$.

\begin{defn}\label{TopRepDefn} \index{wtopologicalrepresentative@topological representative}
Let $(G, m)$ be a marked graph, $f:G \to G$ a homotopy equivalence
which fixes $b$ and $\phi$ the automorphism of $F$ determined by
$f$. We say that the triple $(G, m, f)$ is a \emph{topological
representative} of $\phi$.  More usually, we omit mention of $m$
from the notation, and say simply that $f:G \to G$ is a
topological representative of $\phi$.
\end{defn}


\begin{rem}\label{TopRepsAndCovers}
Let $(G, m, f)$ be a topological representative of $\phi \in
\Aut(F)$. Each finite $G$-cover $(\tilde{G}, p)$ and choice of
point $\tilde{b} \in p^{-1}(b)$ corresponds to a subgroup $S \leq
F$ of finite index.  Let $\tilde{m}: R' \to \tilde{G}$ be the
corresponding marking of $\tilde{G}$, let $k$ be such that $S$ is
$\phi^k$-invariant and let $\tilde{f}: \tilde{G} \to \tilde{G}$ be
the lift of $f^k$ which fixes $\tilde{b}$. Then $(\tilde{G},
\tilde{m}, \tilde{f})$ is a topological representative of
$\phi^k|_S$.
\end{rem}

\subsection{IRTT representatives}\label{IRTTSubsection}

\begin{defn}[IRTT vocabulary]
A \emph{filtration}\index{wfiltration@filtration} for a
topological representative $f:G \to G$ is an increasing sequence
of (not necessarily connected) $f$-invariant subgraphs,
\begin{equation*}\emptyset = G_0 \subset G_1 \subset \dots \subset
G_m = G.\end{equation*} Each set $H_r := G_r\backslash G_{r-1}$ is
called a \emph{stratum} (recall Notation \ref{SetminusNotation}).
A \emph{complete
filtration}\index{wcompletefiltration@complete|see{filtration}}
\index{wfiltration@filtration!wcompletefiltrationsub@complete} is
a filtration such that each $G_i$ is obtained from $G_{i-1}$ by
adding a single geometric edge. A pair consisting of a marked
graph and a (complete) filtration is called a \emph{(completely)
filtered marked graph}. For a complete filtration we usually label
the directed edges of $G$ by $e_1, \rev{e}_1, e_2, \rev{e}_2,
\dots, e_m, \rev{e}_m$ so that the $H_i = \{e_i, \rev{e}_i\}$ for
each $i$. We define the \emph{height}\index{wheight@height of a
path} of a path $\rho \subset G$ (with respect to a filtration),
denoted $h(\rho)$\index{sh@$h$}, to be the maximum value of $i$
for which $\rho$ crosses an edge in $H_i$.  A path $\rho \subset
G$ is a \emph{periodic Nielsen path (for $f$)} if $f^k_{\#}(\rho)
= \rho$ for some $k \geq 1$. If $k = 1$ then we say that $\rho$ is
a \emph{Nielsen path}\index{wNielsenpath@Nielsen path}. A periodic
Nielsen path is said to be
\emph{indivisible}\index{windivisible@indivisible} if it cannot be
written as a concatenation of non-trivial periodic Nielsen paths.
For a tight path $\rho \subset G$, we say that $\rho = \rho_1
\rho_2 \dots \rho_s$ is an
\emph{$f$-splitting}\index{wfsplitting@$f$-splitting}\index{wsplitting@splitting|see{$f$-splitting}}
if $f^k_{\#}(\rho) = f^k_{\#}(\rho_1)f^k_{\#}(\rho_2)\dots
f^k_{\#}(\rho_s)$ for each $k \geq 0$.  We use
$\cdot_{f}$\index{saacdot@$\cdot_{f}$ or $\cdot$} to concatenate
subpaths only if the concatenation is a $f$-splitting, although we
usually omit the map $f$ from the notation if it is clear from the
context. Assume now that we have a complete filtration and an
orientation $\Orient$ where for each $i$ we have $H_i \cap \Orient
= \{e_i\}$ and $f(e_i) = e_i u_i$ for some closed tight path $u_i
\subset G_{i-1}$. A \emph{basic path of height
$i$}\index{wbasicpath@basic path} is a tight path $\rho$ of the
form $e_i \gamma \rev{e}_i$, $e_i \gamma$ or $\gamma \rev{e}_i$
where $e_i \in \Orient$ and $\gamma \subset G_{i-1}$. An
\emph{exceptional path}\index{wexceptionalpath@exceptional path}
is a tight path $\rho$ of the form $e_i
\alpha^k \overline{e}_j$, where $k \in \Integer$, $\alpha$ is a 
closed Nielsen path in $G_{i-1}$, $f(e_i) = e_i \alpha^l$ for some
$l \in \Nat$ and $f(e_j) = e_j \alpha^m$ for some $m \in \Nat$.
\end{defn}

A topological representative of an automorphism in $\PG(F)$ which
satisfies the conclusions of the following theorem (a restriction
of Theorem 5.1.5 \cite{BFH1} to $\PG(F)$) is said to be an
\emph{improved relative train track (IRTT) representative}
\index{wimprovedrtt@improved relative train track representative}.

\begin{thm}[Bestvina, Feighn, Handel - The (\PG) IRTT Theorem]
\label{IRTTThm}\index{wIRTTTheorem@IRTT Theorem} Let $F$ be a
finitely generated free group. For every automorphism $\phi \in
\PG(F)$ there exist a topological representative $f:G \to G$ of an
iterate of $\phi$, a complete filtration
\begin{equation*}\emptyset = G_0 \subset G_1 \subset \dots \subset G_m = G,\end{equation*}
and an orientation $\Orient$ of $G$ such that, if we label the
edges of $\Orient$ such that $H_i \cap \Orient = \{e_i\}$, the
following properties hold:
\begin{enumerate}
\item [(TT1)] each vertex $v \in G$ is fixed by $f$;
\item [(TT2)] each periodic Nielsen path has period one;
\item [(TT3)] for each $i$, either $f(e_i) = e_i$ or
$f(e_i) = e_i \cdot u_i$ for some non-trivial closed tight path
$u_i \subset G_{i-1}$  ($u_i$ is called the \emph{$f$-suffix} of
$e_i$);\index{wfsuffix@$f$-suffix}\index{wsuffix@suffix|see{$f$-suffix}}
\item [(TT4)] if $\sigma \subseteq G_i$ is a basic path of
height $i$ at least one of the following occurs:
\begin{enumerate}
    \item $\sigma$ $f$-splits as a concatenation of two basic paths of
    height $i$;
    \item $\sigma$ $f$-splits as a concatenation of a basic path of
    height $i$ with a tight path contained in $G_{i-1}$;
    \item Some $f^k_\#(\sigma)$ $f$-splits into pieces, one of which
    equals $e_i$ or $\rev{e}_i$;
    \item $u_i$ is a Nielsen path and $\sigma$ is an
    exceptional path of height $i$.
\end{enumerate}
\end{enumerate}
\end{thm}

\begin{cor}\label{CorollaryToIRTT}
Let $f:G \to G$ be an IRTT representative of $\phi \in \PG(F)$,
with the notation of the IRTT Theorem.  The following properties
hold:
\begin{enumerate}
\item \label{CanSplitIntoBasicPaths}
each tight path $\rho \subset G$ may be $f$-split into pieces
which are either basic paths of height $h(\rho)$ or paths of
height less than $h(\rho)$.
\item \label{PathsSplitIntoSingleEdgesAndStuff}
for each tight path $\rho$ in $G$ there exists an integer $M =
M(\rho)$ such that, for each $m \geq M$, $f^m_\#(\rho)$ $f$-splits
into subpaths, each of which is either a single edge, the $k$-th
iterate of an $f$-suffix (or its reverse) for some $k \in
\Integer_+$, or an exceptional path.
\item \label{IRTTAndFiniteCovers}
let $\tilde{G}$ be a finite-sheeted $G$-cover, let $\tilde{b} \in
p^{-1}(b)$ and let $S \leq F$ be the subgroup of finite index
corresponding to $\pi_1 (\tilde{G}, \tilde{b})$. There exists $k
\in \Nat$ such that the following properties hold:
\begin{enumerate}
\item $S$ is $\phi^k$-invariant;
\item $\tilde{f}: \tilde{G} \to
\tilde{G}$ is an IRTT representative of $(\phi^k)|_S$, where
$\tilde{f}: \tilde{G} \to \tilde{G}$ denotes the lift of $f^k$
which fixes $\tilde{b}$.
\end{enumerate}
\end{enumerate}
\end{cor}

\begin{proof}
It follows from (TT3) that we may $f$-split any tight path $\rho$
in $G$ immediately before an occurrence of $e_{h(\rho)}$ and
immediately after an occurrence of $\rev{e}_{h(\rho)}$. If we
$f$-split $\rho$ at each such point, we write $\rho$ as a
concatenation of basic paths of height $h(\rho)$ and paths of
height less than $h(\rho)$. Thus Property
(\ref{CanSplitIntoBasicPaths}) holds.

We prove Property (\ref{PathsSplitIntoSingleEdgesAndStuff}) by
induction on $h(\rho)$. If $h(\rho) = 1$ then $\rho = e_1 \cdot
e_1 \cdot \ldots \cdot e_1$ or $\rho = \rev{e}_1 \cdot \rev{e}_1
\cdot \ldots \cdot \rev{e}_1$. Suppose that, for some integer $k$
such that $2 \leq k < h(G)$, the conclusions of Property
(\ref{PathsSplitIntoSingleEdgesAndStuff}) hold for each tight path
in $G$ of height less than $k$.  By Property
(\ref{CanSplitIntoBasicPaths}) and the inductive hypothesis, to
complete the inductive step it is enough to consider only the case
that $\rho$ is a basic path of height $k$.  We use a second
induction on the length of $\rho$. If $l(\rho) = 1$ there is
nothing more to prove. Suppose the conclusions of Property
(\ref{PathsSplitIntoSingleEdgesAndStuff}) hold for each basic path
of height $i$ and length at most $j \geq 1$.  Suppose that $\rho$
has length $j+1$.  That $\rho$ satisfies the conclusions of
Property (\ref{PathsSplitIntoSingleEdgesAndStuff}) follows
immediately from (TT4) and the two inductive hypotheses.  This
completes the proof of Property
(\ref{PathsSplitIntoSingleEdgesAndStuff}).

Now consider Property (\ref{IRTTAndFiniteCovers}).  For some $k_0
\in \Nat$, $\phi^{k_0}$ leaves $S$ invariant and it follows that
there exists a lift $\tilde{f}':\tilde{G} \to \tilde{G}$ of
$f^{k_0}$ which fixes a particular vertex $\tilde{v} \in
\tilde{G}$. Let $s$ be the number of sheets in the covering
$\tilde{G}$ and let $k_1 = s \abs{\mathcal{V}_G}$.  Clearly,
$\tilde{f} = (\tilde{f}')^{k_1}$ fixes each vertex of $\tilde{G}$,
that is, $\tilde{f}$ has property (TT1).  It is also clear that
the orientation $\Orient$ of $G$ induces an orientation
$\tilde{\Orient}$ of $G$ and, by choosing an order on the elements
of the set $p^{-1}(e)$ for each $e \in \Orient$, we may choose a
complete filtration of $\tilde{G}$ which corresponds to the
complete filtration of $G$. Properties (TT2), (TT3) and (TT4)
follow easily from the corresponding properties of $f$. Hence
Property (\ref{IRTTAndFiniteCovers}) holds with $k = k_0 k_1$.
\end{proof}

\begin{rem}
Corollary \ref{CorollaryToIRTT} (\ref{IRTTAndFiniteCovers}) may be
used to construct examples of IRTT representatives of
automorphisms $\phi \in \PG(F)$ where $F$ has large rank and the
growth of $\phi$ is either exponential or polynomial of small
degree.
\end{rem}

It is convenient to make the following additional definitions.

\begin{defn}[Further IRTT vocabulary]
A complete filtration and an orientation of $G$ which satisfy the
conditions of the IRTT Theorem are said to be
\emph{compatible}\index{wcompatible@compatible} with $f$.  We say
that a closed tight path $\rho$ in $G$ is a \emph{well-chosen
closed tight path}\index{wwellchosen@well-chosen}
\index{wpath@path!wwellchosensub@well-chosen} if either the
initial edge of $\rho$ is $e_{h(\rho)}$ or the terminal edge of
$\rho$ is $\rev{e}_{h(\rho)}$ but not both. A finite tight path
$\alpha$ in $G$ is said to be \emph{essentially
unbounded}\index{wessentiallyunbounded@essentially unbounded}
\index{wpath@path!wessentiallyunboundedsub@essentially unbounded}
if it is not a subpath of any Nielsen path. A $G$-immersion $H$ is
said to be
\emph{$f$-stable}\index{wfstable@$f$-stable}\index{wstable@stable|see{$f$-stable}}
if there exists $q \in \Nat$ such that, for each edge $d \in H$,
$f^q_{\#}(\hat{d})$ labels a path from $\iota(d)$ to $\tau(d)$ (in
which case, the minimum such $q$ is denoted $\period(H)$).
\end{defn}

%

We record some simple properties of the above definitions.

\begin{rem}[A property of well-chosen closed tight paths]
If $\sigma$ is a circuit in $G$ and $\rho$ is a well-chosen closed
tight path representing $\sigma$, then, for each non-negative
integer $k$, $f^{k}_{\#}(\rho)$ is a well-chosen closed tight path
representing $f^k_{\#}(\sigma)$.
\end{rem}

%

\begin{rem}[Properties of essentially unbounded paths]\label{PropertiesOfEssentiallyUnbounded}
Observe the following:
\begin{enumerate}
\item by definition, Nielsen paths contain no essentially unbounded
subpaths;
\item if $\alpha$ is essentially unbounded then
$\rev{\alpha}$ is essentially unbounded;
\item if $\alpha$ is essentially unbounded and $\alpha$ is a
subpath of $\beta$ then $\beta$ is an essentially unbounded
subpath.
\end{enumerate}
\end{rem}



\begin{lem}[Properties of $f$-stable $G$-immersions]
Let $f:G \to G$ be an IRTT representative (of some automorphism
$\phi \in \PG(F)$) and assume the notation of the IRTT Theorem,
let $H$ be a $G$-immersion and $\Orient_H$ the orientation of $H$
induced by $\Orient$ (the orientation of $G$).  The following
properties hold:
\begin{enumerate}
\item $H$ is $f$-stable if and only if, for
each edge $d \in \Orient_H$, there exists a path $\rho = \rho(d)$
in $H$ with $\iota(\rho) = \iota(d)$ and $\hat{\rho} =
f_{\#}(\hat{d})$;
\item if $H$ is $f$-stable with $\period(H) = q$, then for each path
$\alpha \subset H$ and for each $k \in \Integer$,
$f^{kq}_{\#}(\hat{\alpha})$ labels a path from $\iota(\alpha)$ to
$\tau(\alpha)$.
\end{enumerate}
\end{lem}

\subsection{The growth of paths in IRTT representatives}\label{IRTTGrowthSection}

Topological representatives allow us to think of closed tight
paths in $G$ rather than elements of $F$, of tight circuits in $G$
rather than conjugacy classes of $F$ and of homotopy equivalences
of $G$ rather than automorphisms of $F$. Our interest is in the
growth of a basis under repeated application of an automorphism
$\phi \in \Aut(F)$. The aim of this section is to prove Corollary
\ref{FindCircuitToFinishJob} below, which informs us that we may
understand much about the growth of $\phi$ if we understand the
growth of tight circuits in $G$ under the map $f_\#$.

Let $(G, m)$ be a marked graph and let $A$ be a generating set of
$F$. For each $a \in A$, let $\rho_a$ be the closed tight path at
$b$ corresponding to $a$.

\begin{defn}[Growth of a homotopy equivalence]
\index{sgrowthofahomotopy@$\mathcal{G}_{f}$}
\index{sgrowthofahomotopyabelian@$\mathcal{G}^{\ab}_{f}$} For each
homotopy equivalence $f: G \to G$ which fixes $b$, define
$\norm{f}_A := \max \{l(f_{\#}(\rho_a)) \hbox{ } | \hbox{ } a \in
A \}$; define $\mathcal{G}_{f, \, A}: \Nat \to \Nat$ by
$\mathcal{G}_{f, \, A}(k) = \norm{f^k}_A$; define
$\norm{f}^{\ab}_{A} := \max \{l^{\ab}(f_{\#}(\rho_a)) \hbox{ } |
\hbox{ } a \in A \}$; and define $\mathcal{G}^{\ab}_{f, \, A}:
\Nat \to \Nat$ by $\mathcal{G}^{\ab}_{f, \, A}(k) =
\norm{f^k}^{\ab}_{A}$.
\end{defn}

\begin{rem}
As in Proposition \ref{PropertiesOfGrowthOfAutomorphism} (G1), it
is easily verified that each of the functions defined above is
$\simeq$-independent of $A$, and we usually omit mention of the
generating set in our notation.
\end{rem}

We record the following obvious but important consequence of the
above definition.

\begin{lem} \label{PointOfTopReps}
Let $\phi \in \Aut(F)$ be an automorphism and let $f:G \to G$ be a
topological representative of $\phi$. Then $\mathcal{G}_f \simeq
\mathcal{G}_\phi$ and $\mathcal{G}^{\ab}_{f} \simeq
\mathcal{G}^{ab}_{\phi}$.
\end{lem}

\begin{defn} [Growth of a tight path or tight circuit]
\index{wgrowthofapath@growth of a path}
\index{sgrowthofapath@$\mathcal{G}_{f, \, \rho}$}
\index{sgrowthofapathabelian@$\mathcal{G}^{\ab}_{f,\, \rho}$} For
each tight path $\rho$ in $G$, define $\mathcal{G}_{f, \, \rho}:
\Nat \to \Nat$ by $\mathcal{G}_{f, \, \rho}(k) =
l(f^{k}_{\#}(\rho))$ and $\mathcal{G}^{\ab}_{f, \, \rho}: \Nat \to
\Nat$ by $\mathcal{G}^{\ab}_{f, \, \rho}(k) =
l^{\ab}(f^{k}_{\#}(\rho))$. For each tight circuit $\sigma$
represented by a closed tight path $\rho$, define $\mathcal{G}_{f,
\, \sigma}: \Nat \to \Nat$ by $\mathcal{G}_{f, \, \sigma}(k) =
l^{\circlearrowleft}(f^{k}_{\#}(\rho))$ and $\mathcal{G}^{\ab}_{f,
\, \sigma}: \Nat \to \Nat$ by $\mathcal{G}^{\ab}_{f, \, \sigma}(k)
= l^{\ab}(f^{k}_{\#}(\rho))$.
\end{defn}

\begin{cor}\label{IRTTRepsAndGrowth}
Let $f:G \to G$ be an IRTT representative of $\phi \in \PG(F)$,
with the notation of the IRTT Theorem. The following statements
hold:
\begin{enumerate}
\item \label{SplittingRespectsGrowth}
if $\rho = \mu \cdot \nu$ is a tight path in $G$ then
$\mathcal{G}_{f,\,\rho}(k) = \mathcal{G}_{f,\,\mu}(k) +
\mathcal{G}_{f,\,\nu}(k)$ for each $k \in \Nat$;
\item \label{CircuitsPathsAndGrowth}
let $\sigma$ be a tight circuit in $G$ and $\rho$ a closed tight
path representing $\sigma$.  Then $\mathcal{G}_{f, \, \sigma} \leq
\mathcal{G}_{f, \, \rho}$ and $\mathcal{G}^{\ab}_{f, \, \sigma} =
\mathcal{G}^{\ab}_{f, \, \rho}$.  Further, in the case that $\rho$
is a well-chosen closed tight path, the function $G_{f, \,
\sigma}$ is unbounded if and only if $\rho$ is an essentially
unbounded path;
\item \label{EdgesGrowOfIntegerDegree}
for each $i = 1, 2, \dots, h(G)$, there exist $c_i \in \Nat$ such
that $\mathcal{G}_{f, \, e_i} \in p_{c_i}$. Further, if $e_i$ is
not fixed by $f$ (so $f_{\#}(e_i) = e_i \cdot u_i$), there exists
$d_i \in \Nat$ such that $\mathcal{G}_{f, \, u_i} \in p_{d_i}$ and
the following properties hold: $c_i = 1$ and $\mathcal{G}_{f, \,
e_i}$ is linear if and only if $\mathcal{G}_{f, \, u_i}$ is
constant (that is, $u_i$ is a Nielsen path); $c_i = 2$ if and only
if $d_i = 1$ and $\mathcal{G}_{f, \, u_i}$ is linear; and $c_i
\geq 3$ if and only if $d_i = c_i -1 \geq 2$;
\item \label{ExceptionalPathsLinearEdges}
if an exceptional path $\rho$ crosses an geometric edge $E = \{e,
\rev{e}\}$, then $\mathcal{G}_{f, \, e} \in p_1$.
\end{enumerate}
\end{cor}

\begin{proof}
Properties (\ref{SplittingRespectsGrowth}) and
(\ref{CircuitsPathsAndGrowth}) are immediate from the definitions.
Property (\ref{EdgesGrowOfIntegerDegree}) is proved by induction
using the observation that, by (TT3),
\begin{eqnarray}
\mathcal{G}_{f, \, e_i}(n) & = &
l(f^{n}_{\#}(e_i)) \nonumber \\
& = & l(e_i \cdot u_i \cdot \ldots \cdot f^{n-1}_{\#}(u_i)) \nonumber \\
 & = & 1 + l(u_i) + l(f_{\#}(u_i)) + \dots + l(f_{\#}^{n-1}(u_i))
 \nonumber \\
 & = & 1 + \sum_{i = 1}^{n-1} \mathcal{G}_{f, \, u_i}(i). \nonumber
\end{eqnarray}
Property (\ref{ExceptionalPathsLinearEdges}) follows immediately
from Remark \ref{EfficientFiltrationsExist} below and the
observation that the initial and terminal edges of an exceptional
path $\rho$ have linear growth function, and each other edge
crossed by $\rho$ is crossed by the suffix of the initial edge.
\end{proof}

\begin{rem} [Efficient filtration]\label{EfficientFiltrationsExist}
It follows from Property (\ref{EdgesGrowOfIntegerDegree}) that we
may choose a compatible filtration of $G$ and integers $L_1, L_2,
\dots, L_{\eta+1}$ such that $0 < L_1 < L_2 < \dots < L_{\eta+1} =
h(G)+1$ and the following properties hold: if $i < L_1$ then
$f(e_i) = e_i$; if $L_1 \leq i < L_2$ then $\mathcal{G}_{f, \,
e_i}$ is linear and $u_i \subset G_{i-1}$; for $2 \leq j \leq
\eta$, if $L_j \leq i < L_{j+1}$ then $\mathcal{G}_{f, \, e_i} \in
p_j$ and $u_i \subset G_{L_j-1}$. Such a filtration is called
\emph{efficient (with respect to
$f$)}\cite{Macura}\index{wefficientfiltration@efficient|see{filtration}}\index{wfiltration@filtration!wefficientcompletefiltrationsub@efficient}.
In the case of an efficient filtration, define a map $\edgedegree:
\{1, 2, \dots, h(G)\} \to \{0, 1, \dots, d\}$ such that
$L_{\edgedegree (i)} \leq i < L_{\edgedegree(i)+1}$.
\end{rem}

\begin{rem}
Note an important difference between the linear and non-linear
cases in the above: in the case that $L_1 \leq i < L_2$, $u_i
\subset G_{\mathbf{i-1}}$, while in the case that $L_j \leq i <
L_{j+1}$, $u_i \subset G_{\mathbf{L_j-1}}$. This subtlety has a
profound effect on the structure of our proof of the Main Theorem
(see Remark \ref{WhyDifferentCases}).
\end{rem}

\begin{cor}\label{GrowthOfPathOrCircuitLikeHighestEdge}
Let $f:G \to G$ be an IRTT representative of $\phi \in \PG(F)$
with the notation of the IRTT Theorem.  If $\Orient$ is an
efficient filtration then $\mathcal{G}_{f} \simeq \mathcal{G}_{f,
\, e_{h(G)}}$ and, for each tight path (or circuit) $\rho$ in $G$,
$\mathcal{G}_{f, \, \rho} \simeq \mathcal{G}_{f, \, e_{h(\rho)}}$.
\end{cor}

\begin{proof}
Observe that the first part of the conclusion is implied by the
second part of the conclusion.  Assume the notation of Remark
\ref{EfficientFiltrationsExist}. Let $\rho$ be a tight path in
$G$. We use induction on $h(\rho)$.  If $h(\rho) < L_2$ then
$\mathcal{G}_{f, \, \rho}$ is bounded above by the linear function
$k \mapsto k M l^{\circlearrowleft}(\rho)$, where $M := \max
\{l(u_i) \hbox{ } | \hbox{ } L_1 \leq i < L_2 \}$. Hence
$\mathcal{G}_{f, \, \rho}$ and $\mathcal{G}_{f, \, e_{h(\rho)}}$
are both elements of $p_1$ and are $\simeq$-equivalent. Suppose
$h(\rho) \geq L_2$. It follows from Corollary
\ref{CorollaryToIRTT} (\ref{PathsSplitIntoSingleEdgesAndStuff})
and Corollary (\ref{IRTTRepsAndGrowth})
(\ref{ExceptionalPathsLinearEdges}) that, for some $m \in \Nat$,
we may $f$-split $f^m_{\#}(\rho)$ into subpaths $\rho_1 \cdot
\rho_2 \cdot \ldots \cdot \rho_s$, one of which is $e_h(\rho)$ or
$\rev{e}_{h(\rho)}$. Hence $\mathcal{G}_{f, \, \rho} \preceq
\mathcal{G}_{f, \, e_{h(\rho)}}$.  By the inductive hypothesis and
the definition of an efficient filtration, $\mathcal{G}_{f, \,
\rho_i} \preceq \mathcal{G}_{f, \, e_{h(\rho)}}$ for each $ i = 1,
2, \dots, s$, and hence by Corollary \ref{IRTTRepsAndGrowth}
(\ref{SplittingRespectsGrowth}),
\begin{equation*}
\mathcal{G}_{f, \, \rho} = \underset{i = 1}{\overset{s}{\sum}}
\mathcal{G}_{f, \, \rho_i} \preceq \mathcal{G}_{f, \,
e_{h(\rho)}}.
\end{equation*}
Hence $\mathcal{G}_{f, \, \rho} \simeq \mathcal{G}_{f, \,
e_{h(\rho)}}$ and the result holds.
\end{proof}

\begin{cor}\label{FindCircuitToFinishJob}
There exists a circuit $\sigma$ in $G$ such that $\mathcal{G}_{f,
\, \sigma} \simeq \mathcal{G}_f$.  Further, $\mathcal{G}^{\ab}_{f}
\simeq \mathcal{G}_f$ if and only if there exists a circuit
$\sigma$ in $G$ such that $\mathcal{G}^{\ab}_{f, \, \sigma} \simeq
\mathcal{G}_f$.
\end{cor}

\begin{proof}
Assume that $\Orient$ is an efficient filtration of $G$.  Since
$G$ is a minimal graph there exists an circuit $\sigma$ which
crosses $e_{h(G)}$. By Corollary
\ref{GrowthOfPathOrCircuitLikeHighestEdge}, $\mathcal{G}^{\ab}_{f,
\, \sigma} \simeq \mathcal{G}_f$. The second part of the corollary
is then immediate by the observation that, for each circuit
$\delta$ in $G$, $\mathcal{G}^{\ab}_{f, \, \delta} \preceq
\mathcal{G}^{\ab}_f \preceq \mathcal{G}_f$.
\end{proof}

\subsection{Proofs of growth properties for elements of
Aut(F)}\label{Proofsofgrowthpropertiessection}

We begin this subsection with a proof of Theorem
\ref{AutomorphismGrowthTheorem}.

\begin{proof}[{\bf Proof of Theorem \ref{AutomorphismGrowthTheorem}:}]
It is immediate from the definition of \PG(F) that, for each
automorphism $\phi \in \Aut(F)$, $\mathcal{G}_{\phi} \in
p_{\expgrowth}$ or $\phi \in \PG(F)$. Suppose $\phi \in \PG(F)$.
By the IRTT Theorem there exists $k \in \Nat$ such that $\phi^k$
has an IRTT representative. By Proposition
\ref{PropertiesOfGrowthOfAutomorphism} (G2), $\mathcal{G}_{\phi}
\simeq \mathcal{G}_{\phi^k}$, thus we may assume that $\phi$ has
an IRTT representative $f:G \to G$.  By Lemma
\ref{PointOfTopReps}, it is enough to show that $\mathcal{G}_f$
satisfies the conclusions of the theorem.

It is clear from the definitions, Corollary
\ref{IRTTRepsAndGrowth} (\ref{EdgesGrowOfIntegerDegree}) and
Corollary \ref{GrowthOfPathOrCircuitLikeHighestEdge} that
$\mathcal{G}_{f} \in p_\eta$ for some integer $\eta \geq 1$. It
remains to show that $\eta < n$. If $G_{f} \in p_1$ there is
nothing to prove, so we may suppose that $\eta \geq 2$. For each
$1 \leq i \leq h(G)$, let $S_i$ denote the connected component of
$G_{L_{\edgedegree(i)}-1}$ which contains $u_i$. Since $G$ is a
minimal graph, we may choose a maximal subtree $T \subset G$ such
that $T$ does not contain $E_{h(G)}$.  Let $T_i$ denote $T \cap
S_i$. Recall, the number of geometric edges in $G \setminus T$ is
$n$ (the rank of $F$). Thus it suffices to prove the following
claim by induction on $\edgedegree(i)$: $S_i \setminus T_i$
contains at least $\edgedegree(i)$ edges.

Let $i$ be an integer such that $\edgedegree(i) = 2$. The subgraph
$S_i$ contains the non-trivial closed tight path $u_i$ and hence
$S_i \setminus T_i$ contains at least one edge. Suppose that $S_i
\setminus T_i$ contains exactly one edge. Since $u_i$ is a linear
tight path, $S_i$ contains at least one linear edge. Let $j$ be
minimal such that $E_j \subset S_i$ and $e_j$ is linear. Then
$u_j$ is a closed Nielsen path which crosses only fixed edges, and
the unique edge in $S_i \setminus T_i$ must be fixed. It follows
that $S_i$ contracts onto a circle of fixed edges, and hence that
each closed tight path in $S_i$ is a Nielsen path. This
contradicts the fact that $\mathcal{G}_{f, \, u_i}$ is linear,
hence $S_i \setminus T_i$ contains at least two edges. Now suppose
that, for some $k \geq 2$, we have, for each integer $j$ such that
$\edgedegree(j) = k$, $S_j \setminus T_j$ contains at least $k$
edges. Let $i$ be an integer such that $\edgedegree(i) = k+1$.  By
Corollary (\ref{IRTTRepsAndGrowth})
(\ref{EdgesGrowOfIntegerDegree}), $u_i$ crosses an edge $E_j$ such
that $\edgedegree(j) = k$.  By the inductive hypothesis, $S_j
\setminus T_j$ contains at least $k$ edges, and since $S_j \subset
S_i$, we have that $S_i \setminus T_i$ contains at least $k$
edges. Suppose that $S_i \setminus T_i$ contains exactly $k$
edges.  Then $S_i$ contracts onto $S_j$ and it follows that $u_i =
v_1 v_2 \rev{v}_1$ for some paths $v_1$ in $S_i \setminus S_j$ and
$v_2$ in $S_j$. It follows from (TT3) of the IRTT Theorem that
$v_1 = e_{h(u_i)}$. This contradicts (TT3) (since $u_j
f_{\#}(u_j)$ is not a tight path), hence $S_i \setminus T_i$ must
contain at least $k+1$ edges, and the induction is complete.
\end{proof}



\begin{defn}[The reverse of a homotopy equivalence]\index{saareverse@$\rev{\;}\;$}\index{wreverse@reverse!wrofamap@of
a map}
Let $f:G \to G$ be an IRTT representative (of some
automorphism $\phi \in \PG(F)$). It follows from (TT1) that, for
each tight path $\rho$ in $G$, there is a unique tight path $\mu$
in $G$ such that $f_{\#}(\mu) = \rho$. Thus we may (inductively)
define a map $\overline{f}:G \to
G$\index{sfinverse@$\overline{f}$}, called the \emph{reverse of
$f$}, by sending $e_i$ to $e_i$ if $f(e_i) = e_i$, and otherwise
sending $e_i$ to $e_i v_i$, where $v_i$ is the unique tight path
such that $f_{\#}(v_i) = \rev{u}_i$. We define $\overline{f}_{\#}$
from $\overline{f}$ as we defined $f_{\#}$ from $f$ (see
$\S$\ref{HomotopyEquivSection}).
\end{defn}

\begin{rem}
It is clear that $\overline{f}: G \to G$ is a topological
representative of $\phi^{-1}$. However, in the general case,
$\alpha \beta = \alpha \cdot_{f} \beta$ does not imply that
$\alpha \beta = \alpha \cdot_{\overline{f}} \beta$ (cf Lemma
\ref{InverseImageOfPathUnit}).  In particular, it is not
necessarily the case that $\overline{f}$ is an IRTT representative
(of $\phi^{-1}$), as shown by the following example.
\end{rem}

\begin{eg}
Consider an IRTT representative $f:G \to G$ (of an automorphism
$\phi \in \Aut(F)$) where $G$ is the graph with one vertex, three
edges and an orientation $\{e_1, e_2, e_3\}$ and $f$ is defined by
$f(e_1) = e_1$, $f(e_2) = e_2 \cdot e_1$ and $f(e_3) = e_3 \cdot
e_1 e_2$. The map $\overline{f}$ is given by $\overline{f}(e_1) =
e_1$, $\overline{f}(e_2) = e_2 \rev{e}_1$ and $\overline{f}(e_3) =
e_3 e_1 \rev{e}_2 \rev{e}_1$. Let $v_3 = e_1 \rev{e}_2 \rev{e}_1$.
Then $\overline{f}^{2}_{\#}(e_3) = e_3 e_1 \rev{e}_2 e_1 \rev{e}_2
\rev{e_1} \neq e_3 v_3 \overline{f}(v_3)$, hence
$\overline{f}_{\#}(e_3 v_3) \neq \overline{f}_{\#}(e_3)
\overline{f}_{\#}(v_3)$ and $\overline{f}$ fails (TT3).
\end{eg}

Consideration of the map $\overline{f}$ allows us to prove Theorem
\ref{InverseHasSameGrowth}.

\begin{proof}[{\bf Proof of Theorem \ref{InverseHasSameGrowth}}]
It suffices to show the result in the case that $\phi \in \PG(F)$.
Let $d \in \Nat$ be such that $\mathcal{G}_{\phi^{-1}} \in p_d$.
Let $f:G \to G$ be an IRTT representative of some iterate
$\phi^{j_0}$ of $\phi$.  Assume the notation of the IRTT Theorem
and Remark \ref{EfficientFiltrationsExist}.  Recall that
$\overline{f}:G \to G$ is a topological representative of
$\phi^{-j_0}$. We prove (inductively) that
$\mathcal{G}_{\overline{f}} \preceq p_{d}$ and hence
$\mathcal{G}_{\phi^{-1}} \preceq \mathcal{G}_{\phi}$. By an
entirely similar argument, we may show that $\mathcal{G}_{\phi}
\preceq \mathcal{G}_{\phi^{-1}}$ and the result follows.

It is clear that, for each Nielsen path $\mu$ in $G$, $k \mapsto
l(\overline{f}^{k}_{\#}(\mu))$ is the constant function $k \mapsto
l(\mu)$ and hence is an element of $p_1$ as required. Let $i$ be
an integer such that $L_1 \leq i < L_2$. Then
\begin{eqnarray}
l(\overline{f}^{k}_{\#}(e_i)) & \leq & l(e_i v_i \overline{f}_{\#}(v_i) \dots \overline{f}^{(k-1)}_{\#}(v_i)) \nonumber \\
 & = & l(e_i v_i^{k-1}) \nonumber \\
 & = & 1 + \sum_{j = 0}^{k-1} l(v_i). \nonumber
\end{eqnarray}
Hence $\bigl(k \mapsto l(\overline{f}^{k}_{\#}(e_i))\bigr) \preceq
p_1$. It immediately follows that the same is true for a path
$\rho$ such that $L_1 \leq h(\rho) < L_2$. Suppose the following
holds for some integer $d$ such that $1 \leq d \leq \eta$: for
each path $\rho$ with $h(\rho) < L_{d+1}$, we have $\bigl(k
\mapsto l(\overline{f}^{k}_{\#}(\rho))\bigr) \preceq p_d$. Let $i$
be an integer such that $L_{d+1} \leq i < L_{d+2}$ and let $v_i =
g_1 g_2 \dots g_t$, for edges $g_1, g_2, \dots, g_t \in
\mathcal{E}$. Then
\begin{eqnarray}
l(\overline{f}^{k}_{\#}(e_i)) & \leq & l(e_i v_i \overline{f}_{\#}(v_i) \dots \overline{f}^{(k-1)}_{\#}(v_i)) \nonumber \\
 & \leq & l(e_i g_1 \dots g_s \overline{f}_{\#}(g_1) \dots \overline{f}_{\#}(g_s) \dots \overline{f}^{(k-1)}_{\#}(g_1) \dots \overline{f}^{(k-1)}_{\#}(g_s)) \nonumber \\
 & = & 1 + \sum_{l=1}^s \sum_{j = 0}^{k-1} l(\overline{f}^{j}_{\#}(g_l)). \nonumber
\end{eqnarray}
Hence, after applying the inductive hypothesis, we have that
$\bigl(k \mapsto l(\overline{f}^{k}(e_i))\bigr) \preceq p_{d+1}$.
It immediately follows that the same is true for a path $\rho$
such that $L_{d+1} \leq h(\rho) < L_{d+2}$ and the induction is
complete.
\end{proof}

\begin{rem}
Theorem \ref{InverseHasSameGrowth} may also be considered a
corollary to the Main Theorem as follows: again, it suffices to
show the result in the case that $\phi \in \PG(F)$.  Let $f:G \to
G$ be an IRTT representative of some iterate of $\phi$.  It
follows easily from the above definition of $\overline{f}$ and the
IRTT Theorem that the $\rev{f}$-growth of each edge in $G$, and
hence each path in $G$, is bounded above by a polynomial function.
It follows that $\phi^{-1} \in \PG(F)$. By the Main Theorem, there
exists a characteristic subgroup $S \leq F$ of finite index such
that, for $\theta = \phi|_S$, $\mathcal{G}^{\ab}_{\theta} \simeq
\mathcal{G}_{\theta} \simeq \mathcal{G}_{\phi}.$  By the Main
Theorem, there exists a characteristic subgroup $S' \leq S$ of
finite index such that, for $\varphi^{-1} = \theta^{-1}|_{S'}$,
$\mathcal{G}^{\ab}_{\varphi^{-1}} \simeq
\mathcal{G}_{\varphi^{-1}} \simeq \mathcal{G}_{\theta^{-1}}.$ It
is easily verified that $\mathcal{G}^{\ab}_{\varphi} \simeq
\mathcal{G}^{\ab}_{\varphi^{-1}}$.  Combining the above with
Proposition \ref{PropertiesOfGrowthOfAutomorphism} (G3), we have
\begin{equation*} \mathcal{G}_{\phi} \simeq  \mathcal{G}_{\theta} \simeq \mathcal{G}^{\ab}_{\theta} \simeq \mathcal{G}^{\ab}_{\varphi}
  \simeq  \mathcal{G}^{\ab}_{\varphi^{-1}}
 \simeq  \mathcal{G}_{\varphi^{-1}}  \simeq
 \mathcal{G}_{\theta^{-1}} \simeq \mathcal{G}_{\phi^{-1}}.\end{equation*}
\end{rem}

We now demonstrate how to construct the sequence $\{U_i\}$ used in
Remark \ref{AlgorithmicProperties} to put an upper bound on $d$
such that $\mathcal{G}_\phi \simeq p_d$.

\begin{rem}[An upper bound for the degree of $\mathcal{G}_\phi$]\label{UpperBoundForAlgorithm}
Let $k \in \Nat$.  We may enumerate the finite minimal graphs with
fundamental group isomorphic to $F$; for each such graph $G$, we
may enumerate the complete filtrations and orientations of $G$;
for each finite minimal graph $G$ with a complete filtration and
an orientation (assuming the usual notation), we may enumerate the
maps $f:G \to G$ such that, for each integer $i$ such that $1 \leq
i \leq h(G)$, either $f(e_i) = e_i$ or $f(e_i) = e_i u_i$ where
$u_i$ is a closed path in $G_{i-1}$; we may enumerate the
identifications between $F$ and the fundamental group of $G$ (at
each base-point); for each such identification, we may assess
whether $f:G \to G$ is a topological representative of $\phi^k$
and if so, we may determine the minimum integer $\eta'$, called
the \emph{degree bound}, such that the following properties hold:
there exist integers $L_0, L_1, \dots, L_{\eta'+1}$ such that $0 =
L_0 < L_1 < L_2 < \cdots < L_{\eta'+1} = h(G)+1$ and, for each
integer $i$ such that $1 \leq i \leq h(G)$,
\begin{itemize}
\item $f(e_i) = e_i$ if and only if $i \leq L_1$;
\item if $L_1 \leq i < L_2$ then $u_i$ is a Nielsen path;
\item  if $L_j \leq i < L_{j+1}$ for some integer $j$ such that
$2 \leq j \leq \eta'$, then $u_i$ is contained in $G_{L_j-1}$.
\end{itemize}
Thus we may enumerate completely filtered topological
representatives $f:G \to G$ of $\phi^k$ with orientations and we
may calculate the corresponding degree bound; let $\{u^k_i\}$ be
the corresponding sequence of degree bounds. We may enumerate the
set $\{u^j_i \, | \, j \in \Nat\}$ by a diagonal process; let
$\{u_i\}$ be such an enumeration and define $U_i := \min \{u_1,
u_2, \dots, u_i\}$. Clearly, $\{U_i\}$ is a non-increasing
sequence and $\mathcal{G}_f \preceq p_{U_i}$ for each $i \in
\Nat$. By the IRTT theorem and Remark
\ref{EfficientFiltrationsExist}, there exists $i_1 \in \Nat$ such
that $\mathcal{G}_f \simeq p_{U_{i_1}}$.  It follows that $i \geq
i_1$ implies $\mathcal{G}_f \simeq p_{U_{i}}$ (and hence
$\mathcal{G}_{\phi} \simeq p_{U_i}$), as required.
\end{rem}

\section{Translating the Main Theorem}\label{TranslatingTheoremChapter}

In this section we reduce the Main Theorem to a theorem stated in
the language of topological representatives.  We prepare to prove
the latter by fixing some notation for the remainder of the paper.

\subsection{The Apt Immersion Theorem}

\begin{thm}[The Apt Immersion
Theorem]\label{FriendlyImmersionTheorem}\index{wAptImmersionTheorem@Apt
Immersion Theorem} Let $f:G \to G$ be an IRTT representative of
$\phi \in \PG(F)$ with an efficient filtration. Let $\sigma$ be a
circuit in $G$, let $v = \iota(e_{h(\sigma)})$ and let $\rho$ be a
well-chosen closed tight path which represents $\sigma$. There
exist a finite $G$-immersion $\Sigma$, a vertex $\tilde{v} \in
p^{-1}(v) \subset \Sigma$ and a natural number $q \in \Nat$ such
that the following properties hold:
\begin{enumerate}
\item [(AI1)] for each non-negative integer $k$, $f^{kq}_{\#}(\rho)$ labels
a closed path $\tilde{\rho}_{kq} \subset \Sigma$ at $\tilde{v}$;
and
\item [(AI2)] $\bigl(k \mapsto l^{\ab}(\tilde{\rho}_{kq})\bigr) \simeq \bigl(k \mapsto l(\tilde{\rho}_{kq})\bigr)$.
Further, if $k \mapsto l(\tilde{\rho}_{kq})$ is unbounded then $k
\mapsto l^{\ab}(\tilde{\rho}_{kq})$ is unbounded.
\end{enumerate}
\end{thm}

\begin{rem}\label{CarriesAndStretchesRemark}
Informally, (AI1) could be understood to be that $\Sigma$
`carries'\index{wcarries@carries} $\rho$ and its images under
iterates of $f_{\#}^q$, and (AI2) could be understood to be that
$\Sigma$ `stretches'\index{wstretches@stretches} $\rho$.
\end{rem}

\begin{proof}[{\bf Proof that the Apt Immersion Theorem implies the Main Theorem}]
\index{wAptImmersionTheorem@Apt Immersion
Theorem!wimpliesMainTheorem@implies the Main Theorem}
\index{wMainTheorem@Main Theorem!wisa consequenceofApt@is a
consequence of the Apt Immersion Theorem} By the IRTT Theorem
there exists $j \in \Nat$ such that $\phi^j$ has an IRTT
representative $f:G \to G$. By Corollary
\ref{FindCircuitToFinishJob} there exists a circuit $\sigma
\subset G$ such that $\mathcal{G}_{f} \simeq
 \mathcal{G}_{f, \, \sigma}$. By the Apt Immersion Theorem there exist
 a finite $G$-immersion $\Sigma$, a  vertex $\tilde{v} \in \Sigma$ and $q
\in \Nat$ such that (AI1) and (AI2) hold. Extend $\Sigma$ to a
finite $G$-cover $\tilde{G}$ by Stallings' Algorithm and choose a
basepoint $\tilde{b} \in \Sigma$ such that $p(\tilde{b}) = b$.
Choose $i \in \Nat $ such that $S' = p_\ast \pi_1(\tilde{G},
\tilde{b})$ is $f^{iq}_{\#}$-invariant and the lift $\tilde{f}:
\tilde{G} \to \tilde{G}$ of $f^{iq}$ which fixes $\tilde{b}$ fixes
all vertices of $\tilde{G}$. By Remark \ref{TopRepsAndCovers},
(there exists a marking $\tilde{m}$ such that) $\tilde{f} :
\tilde{G} \to \tilde{G}$ is a topological representative of
$\theta' := \phi^{ijq}|_{S'}$. By Proposition
\ref{PropertiesOfGrowthOfAutomorphism} (G3), $\mathcal{G}_\phi
\simeq \mathcal{G}_{\theta'}$. By (AI2) and Corollary
\ref{FindCircuitToFinishJob}, $\mathcal{G}^{\ab}_{\tilde{f}}
\simeq \mathcal{G}_{\tilde{f}}$ and $\mathcal{G}^{\ab}_{\theta'}
\simeq \mathcal{G}_{\theta'}$.  Let $S$ be the intersection of all
subgroups of $F$ with index $[F:S']$. Note that $S$ is a
characteristic subgroup of $F$ and $S$ has finite index in both
$F$ and $S'$. Thus we have, for $\theta = \theta'|_S = \phi|_S$,
$\mathcal{G}^{\ab}_{\theta} \simeq \mathcal{G}_{\theta} \simeq
\mathcal{G}_{\theta'} \simeq \mathcal{G}_{\phi}$ and the Main
Theorem holds.
\end{proof}

\subsection{Some notation}\label{NotationSection}

Our remaining task is to prove the Apt Immersion Theorem. For this
purpose we fix some notation for the remainder of the paper: Let
$F$ be a finitely generated free group of rank $n \geq 2$ and let
$\phi \in \PG(F)$\index{sphi@$\phi$} be an automorphism of $F$
which has an IRTT representative $f:G \to G$\index{sf@$f$}
\index{sG@$G$}. Let $\eta \in \Nat$\index{seta@$\eta$} be such
that $\mathcal{G}_{\phi} \in p_\eta$. Let $h(G)$ be the number of
geometric edges in $G$, let $\Orient$\index{sO@$\Orient$} be an
orientation of $G$ determined by $f$ and let $\emptyset = G_0
\subset G_1 \subset \cdots \subset G_{h(G)} = G$\index{sGi@$G_i$}
be an efficient filtration of $G$ with respect to $f$. Label the
directed edges of $\Orient$ by $\{e_1, e_2, \dots,
e_{h(G)}\}$\index{sei@$e_i$} in such a way that $H_i = G_i
\setminus G_{i-1} = E_i = \{e_i, \rev{e}_i\}$\index{sHi@$H_i$} for
each $i = 1, 2, \dots, {h(G)}$. Let $L_0, L_1, \dots,
L_{\eta+1}$\index{sL_i@$L_i$} be integers such that
\begin{equation*}1 = L_0 < L_1 < \cdots  < L_{\eta+1} = {h(G)}+1\end{equation*} and
$\{e_0, \dots, e_{L_1-1}\}$ is the set of edges in $\Orient$ which
are fixed by $f$, $\{e_{L_1}, \dots, e_{L_2-1}\}$ is the set of
edges in $\Orient$ which grow linearly, and for each $j = 2,
\dots, d$, the set $\{e_k \hbox{ } | \hbox{ } L_j \leq k < L_{j+1}
\}$ is the set of edges in $\Orient$ that grow of degree $j$. For
each integer $i$ such that $L_1 \leq i \leq {h(G)}$, let
$u_i$\index{sui@$u_i$} be the $f$-suffix of $e_i$.

\section{The Apt Immersion Theorem in the linear case}\label{LinearCaseChapter}

Let $\rho$ be a path in $G$ with linear growth.  Recall the
notation of $\S$\ref{CirclesAndLines} for $G$-labelled lines and
circles. Consider an end-pointed $G$-labelled graph $\Lambda'$,
constructed from $L := L(\rho)$ as follows:  for each edge $d \in
\mathcal{E}_L$, if $\hat{d} = e_i \in \Orient$ and there is no
tight path $\nu$ in $L$ such that $\nu$ crosses $d$ and
$\hat{\nu}$ is a Nielsen path in $G$, then adjoin a copy of
$C(\mu_i)$ at $\tau(d)$. Let $\Lambda$ denote the end-pointed
$G$-immersion determined by $\Lambda'$.  It is clear that, for
each $k \in \Integer_+$, $f^k(\rho)$ (note the absence of
tightening) labels a path $\tilde{\rho}'_k$ across $\Lambda'$ and
$k \mapsto l^{\ab}(\tilde{\rho}'_k)$ is unbounded.  It follows
that, for each $k \in \Integer_+$, $f^k_{\#}(\rho)$ labels a path
$\tilde{\rho}_k$ across $\Lambda$. Unfortunately, the second
property of $\Lambda'$ may not be inherited by $\Lambda$ as
folding may ``muddle up'' the images of the embedded circles in
$\Lambda'$. We shall prove the linear case of the Apt Immersion
Theorem by arguing that if we adjoin copies of $C(\mu_i^q)$
instead of $C(\mu_i)$, for $q$ sufficiently large, then we gain
tight control on the amount of folding that is required in the
construction of $\Lambda$. This enables us to prove that $\Lambda$
inherits the desirable second property of $\Lambda'$.

We begin the section by introducing growth units, which allow us
to write a tight linear path as a concatenation of subpaths which
interact in a limited way under iteration of $f$.  By performing
the construction of $\Lambda$ in parts, constructing the subgraphs
appropriate for each growth unit and then combining them using our
standard constructions, the amount of folding that may occur
between circles becomes apparent.

\subsection{Separating linear paths into growth
units}\label{GrowthUnitsSection}\index{wgrowthunits@growth units}

\begin{notation}\index{smui@$\mu_i$}
For each $i$ such that $L_1 \leq i < L_2$, define $\mu_i$ to be
the primitive closed path corresponding to $u_i$ and define $m_i
\in \Nat$ such that $u_i = \mu_i^{m_i}$.
\end{notation}

\begin{defn}\label{PassiveLinearGrowthUnitsDefn}
A \emph{passive (linear) growth unit}\index{wgrowthunits@growth
units!wpassivesub@passive}\index{wpassive@passive|see{growth
units}}\index{wgrowthunits@growth units!wFF@of type (FF)}
\index{wFF@(FF)|see{growth units}}\index{wgrowthunits@growth
units!wFR@of type (FR)}\index{wFF@(FR)|see{growth
units}}\index{wgrowthunits@growth units!wFE@of type
(FE)}\index{wFE@(FE)|see{growth units}} is a tight path $\kappa$
in $G$ which is in one of the following forms:
\begin{enumerate}
\item[(FF)] [Fixed forward edge] \\
$\kappa = e_a$ for some $a \in \Nat$ such that $L_0 \leq a < L_1$;
\item[(FR)] [Fixed reverse edge] \\
$\kappa = \rev{e}_b$ for some $b \in \Nat$ such that $L_0 \leq b <
L_1$;
\item[(FE)] [Fixed exceptional path] \\
$\kappa  =  e_a \mu_a^d \rev{e}_b$ for some $a, b \in \Nat$ and $d
\in \Integer$ such that  $L_1 \leq a, b < L_2$, $\mu_a = \mu_b$
and $m_a = m_b$ (note that $d \neq 0$ if $a = b$).
\end{enumerate}
\end{defn}

\begin{lem}\label{NielsenPathsSplitIntoPassiveGrowthUnits}
Every tight Nielsen path in $G$ can be $f$-split into passive
growth units in exactly one way.  That is, for each tight Nielsen
path $\alpha$ in $G$ there is a unique expression $\alpha =
\kappa_1 \cdot \kappa_2 \cdot \ldots \cdot \kappa_s$ such that
each $\kappa_i$ is a passive growth unit.
\end{lem}

\begin{proof}
This follows immediately from (TT4) of the IRTT theorem.
\end{proof}

\begin{notation}\index{skappaij@$\kappa_{i, j}$}
For each $i$ such that $L_1 \leq i < L_2$, let $\mu_i = \kappa_{i,
\,0} \cdot \kappa_{i,  \,1} \cdot \dots \cdot \kappa_{i, \,
s_i-1}$ be the $f$-splitting of $\mu_i$ into passive growth units.
It is convenient to extend this notation by defining $\kappa_{i,
\, j} := \kappa_{i,  j  \, \modulo  \, s_i}$ for each $j \in
\Integer$.
\end{notation}

We would like to decompose tight linear paths into smaller
subpaths. The interactions between subpaths in a linear path are
more complicated than in the case of a Nielsen path, and a
$f$-splitting is not practical. Instead we introduce the
following:

\begin{defn}\label{ActiveLinearGrowthUnitsDefn}
\index{wgrowthunits@growth units!wactivesub@active}
\index{wactive@active|see{growth units}}
\index{wgrowthunits@growth units!wLF@of type (LF)}
\index{wLF@(LF)|see{growth units}} \index{wgrowthunits@growth
units!wLR@of type (LR)}\index{wLR@(LR)|see{growth units}}
\index{wgrowthunits@growth units!wLE@of type (LE)}
\index{wLE@(LE)|see{growth units}} An \emph{active (linear) growth
unit} is a tight path $\delta$ in $G$ which is in one of the
following forms:
\begin{enumerate}
 \item[(LF)] [Linear forward growth unit] \\
 Either $\delta = e_a$ or
 $\delta = e_a \, \kappa_{a, \, 0} \, \kappa_{a, \, 1} \dots \kappa_{a, \,
 (d-1)}$ or
 $\delta = e_a \, \rev{\kappa}_{a, \, -1} \, \rev{\kappa}_{a, \, -2} \dots \rev{\kappa}_{a, \,
 -d}$, for some $a, d \in \Nat$ such that $L_1 \leq a < L_2$;

 \item[(LR)] [Linear reverse growth unit]  \\
 Either
 $\delta = \rev{e}_b$ or
 $\delta = \rev{\kappa}_{b, \, (d-1)} \, \rev{\kappa}_{b, \, (d-2)} \dots \rev{\kappa}_{b, \, 0} \,
 \rev{e}_b$ or
 $\delta = \kappa_{b, \, -d} \, \kappa_{b, \, -d+1} \dots \kappa_{b, \, -1} \,
 \rev{e}_b$, for some $b, d \in \Nat$ such that $L_1 \leq b < L_2$;

 \item[(LE)] [Linear exceptional path] \\
 Either
 $\delta = e_a \rev{e}_b$ or
 $\delta = e_a \mu_a^d \rev{e}_b$ or
 $\delta = e_a \rev{\mu}_a^d \rev{e}_b$,
 for some$a, b, d \in \Nat$ such that $L_1 \leq a, b < L_2$, $\mu_a = \mu_b$ but $m_a \neq
 m_b$;

 \item[(QE)] [Quasi-exceptional path] \\
 Either $\delta = e_a \rev{e}_b,$ for some $L_1 \leq a, b < L_2$ such that $\mu_a =
 \rev{\mu}_b$;
 or $\delta = e_a \, \kappa_{a, \, 0} \, \kappa_{a, \, 1} \dots \kappa_{a, \, (d-1)} \,
 \rev{e}_b,$ for some $a, b, d \in \Nat$ such that $L_1 \leq a, b <
 L_2$, $\mu_b = \rev{\kappa}_{a, \, (d+s_a-1)} \, \rev{\kappa}_{a, \, (d+s_a-2)} \dots \rev{\kappa}_{a, \, d}$; or
 $\delta = e_a \, \rev{\kappa}_{a, \, -1} \, \rev{\kappa}_{a, \, -2} \dots
 \rev{\kappa}_{a, \, -d} \, \rev{e}_b,$ for some $a, b, d \in \Nat$ such that $L_1 \leq a, b,
 < L_2$, $\mu_b = \kappa_{a, \, -(d+s_a)} \, \kappa_{a, \, -(d+s_a-1)} \dots \kappa_{a, \,
 -(d+1)}$.
\end{enumerate}
\end{defn}

\begin{rem}
One might be struck by an asymmetry in the above definition.
Growth units of type (QE) arise in the situation that there exist
integers $a$ and $b$ such that $L_1 \leq a, b < L_2$ and $\mu_b$
is a non-trivial cyclic permutation of the growth units in the
separation of $\rev{\mu}_a$ (see Definition
\ref{CanonicalSeparation}). For growth units of type (FE) and (LE)
we demand that $\mu_a = \mu_b$. There is no growth unit for the
situation that $\mu_b$ is a non-trivial cyclic permutation of the
growth units in the separation of $\mu_a$. The IRTT Theorem
ensures that this may not occur, since otherwise $e_a \kappa_{a,
\, 0} \dots \kappa_{a, \, d-1} \rev{e}_b$ is a Nielsen path which
does not $f$-split but is not an exceptional path, violating
(TT4).
\end{rem}

\begin{rem}
Observe that if $\delta$ is a growth unit then $\rev{\delta}$ is
also a growth unit, although possibly of a different type.
\end{rem}

Each linear path may be written as a concatenation of growth units
in a trivial way --- simply regard each fixed edge as a growth
unit of type (FF) or (FR) and each linear edge as a growth unit of
type (LF) or (LR) --- and in general there is more than one way to
write a linear path as a concatenation of growth units. Writing a
path as a concatenation of growth units is useful only when the
concatenation distinguishes parts of the path which interact in a
limited way under iteration of $f$.

\begin{defn}[The canonical separation of a linear
path]\label{CanonicalSeparation}\index{wcanonicalseparation@canonical
separation}\index{wseparation@separation|see{canonical
separation}} We describe an algorithm to write a tight linear path
$\rho$ as a concatenation of growth units in a canonical manner;
we call this the \emph{separation of $\rho$ into (linear) growth
units}, or more usually, the \emph{separation} of $\rho$. Write
$\rho = d_1 d_2 \dots d_n$ for some $n \in \Nat$ and edges $d_i
\subset G$. First we define $\delta_1$:
\begin{enumerate}
 \item [Step 1]
 if some initial subpath $\nu$ of $\rho$ is a growth unit of type
 (LR), (LE) or (QE) then define $\delta_1 := \nu$
 (note that $\rho$ has at most one such initial subpath);
 \item [Step 2]
 otherwise, if some initial subpath $\nu$ of $\rho$ is a growth unit of
 type (FE) then define $\delta_1 := \nu$;

 \item [Step 3]
 otherwise, if $d_1$ is a forward linear edge then define $\delta_1$ to be the maximal initial subpath
 of $\rho$ which is a growth unit of type (LF);
 \item [Step 4]
 otherwise, define $\delta_1 := d_1$ (note, $d_1$ is a fixed edge).
\end{enumerate}

Inductively, assume that $\delta_1, \dots, \delta_j$ are defined
but that $\rho \neq \delta_1 \dots \delta_j$.  Let $\rho'$ be the
terminal subpath of $\rho$ such that $\rho = \delta_1 \dots
\delta_j \rho'$.  Define $\delta_{j+1}$ from $\rho'$ in the same
way that $\delta_1$ is defined from $\rho$.
\end{defn}

\begin{notation}\index{saadiamond@$\separate$}
We denote that a concatenation of subpaths $\rho = \delta_1
\delta_2 \dots \delta_2$ is in fact the canonical separation of
$\rho$ by writing the symbol $\separate$ between subpaths, that
is, we write $\rho = \delta_1 \separate \delta_2 \separate \dots
\separate \delta_s$.
\end{notation}

\begin{rem}
The separation of a path is not necessarily a $f$-splitting. A
Nielsen path, however, separates into passive growth units only,
in which case the separation is the unique $f$-splitting of Lemma
\ref{NielsenPathsSplitIntoPassiveGrowthUnits}.
\end{rem}

\begin{rem}\label{SeparationNotSymmetric}
The separation of a path as defined above is not, in general,
symmetric, that is, $\rho = \delta_1 \separate \dots \separate
\delta_s$ does not necessarily imply that $\rev{\rho} =
\rev{\delta}_s \separate \dots \separate \rev{\delta}_1$. Symmetry
may be arranged by replacing Step 3 of the algorithm by a process
which ensures that, for $\delta_i$ a type (LF) growth unit and
$\delta_{i+1}$ a type (LR) growth unit, $\rev{(\delta_i \separate
\delta_{i+1})} = \rev{\delta}_{i+1} \separate \rev{\delta}_i$
(such as applying some order on the set of paths in $G$ to express
a preference for maximising the length of $\delta_i$ or
$\delta_{i+1}$).  Since it is not required in the argument below,
we have opted to sacrifice symmetry for a simpler algorithm.
\end{rem}

We record a simple property of the separation of a linear path.

\begin{lem}\label{PropertiesOfSeparation}
Let $\rho = \delta_1 \separate \delta_2 \separate \dots \separate
\delta_s$ be a linear path in $G$. If $s \geq 3$ and $\delta_2
\dots \delta_{s-1}$ is not a Nielsen path then $\rho$ contains an
essentially unbounded subpath.
\end{lem}

\begin{proof}
The lemma is immediate from the following observations:  if
$\delta_i$ is a growth unit of type (LE) or (QE) then $\delta_i$
is essentially unbounded; if $\delta_i$ is a growth unit of type
(LF) then $\delta_i \delta_{i+1}$ is essentially unbounded; and if
$\delta_i$ is a growth unit of type (LR) then $\delta_{i-1}
\delta_i$ is essentially unbounded.
\end{proof}

\subsection{Diagram units}

\begin{con}[Diagram units]\label{DiagramUnitsDefn}\index{wdiagramunits@diagram
units}\index{sLambdaBig1@$\Lambda(\,)$}\index{sLambdaBig2@$\Lambda[\,]$}
Let $\delta$ be a growth unit. Define an end-pointed $G$-labelled
graph $\Lambda(\delta, q)$ as follows:
\begin{enumerate}
\item
if $\delta$ is a passive growth unit then $\Lambda(\delta, q) :=
L(\delta)$;
\item
if $\delta$ has type (LF), (LE) or (QE) with initial edge $e_a$,
say, then define \begin{equation*}\Lambda(\delta, q) := L(\delta)
\amalg C(\mu_a^q) / \sim,\end{equation*} where $\sim$ identifies
the terminal point of the edge in $L(\delta)$ labelled by $e_a$
with the basepoint of $C(\mu_i^q)$;
\item
if $\delta$ has type (LR) with terminal edge $\rev{e}_b$, say,
then define \begin{equation*}\Lambda(\delta, q) := L(\delta)
\amalg C(\mu_b^q) / \sim,\end{equation*} where $\sim$ identifies
the terminal point of the edge in $L(\delta)$ labelled by $e_b$
with the basepoint of $C(\mu_i^q)$.
\end{enumerate}
In each case define the end-points of $\Lambda(\delta, q)$ to be
the natural images of the end-points of $L(\delta)$ and define
$\Lambda[\delta, q]$ to be the end-pointed $G$-immersion
determined by $\Lambda(\delta, q)$ (see Figure
\ref{DiagramUnitsFigure}).
\end{con}

\setcounter{figure}{\value{thm}} \stepcounter{thm}
\begin{figure}
\includegraphics{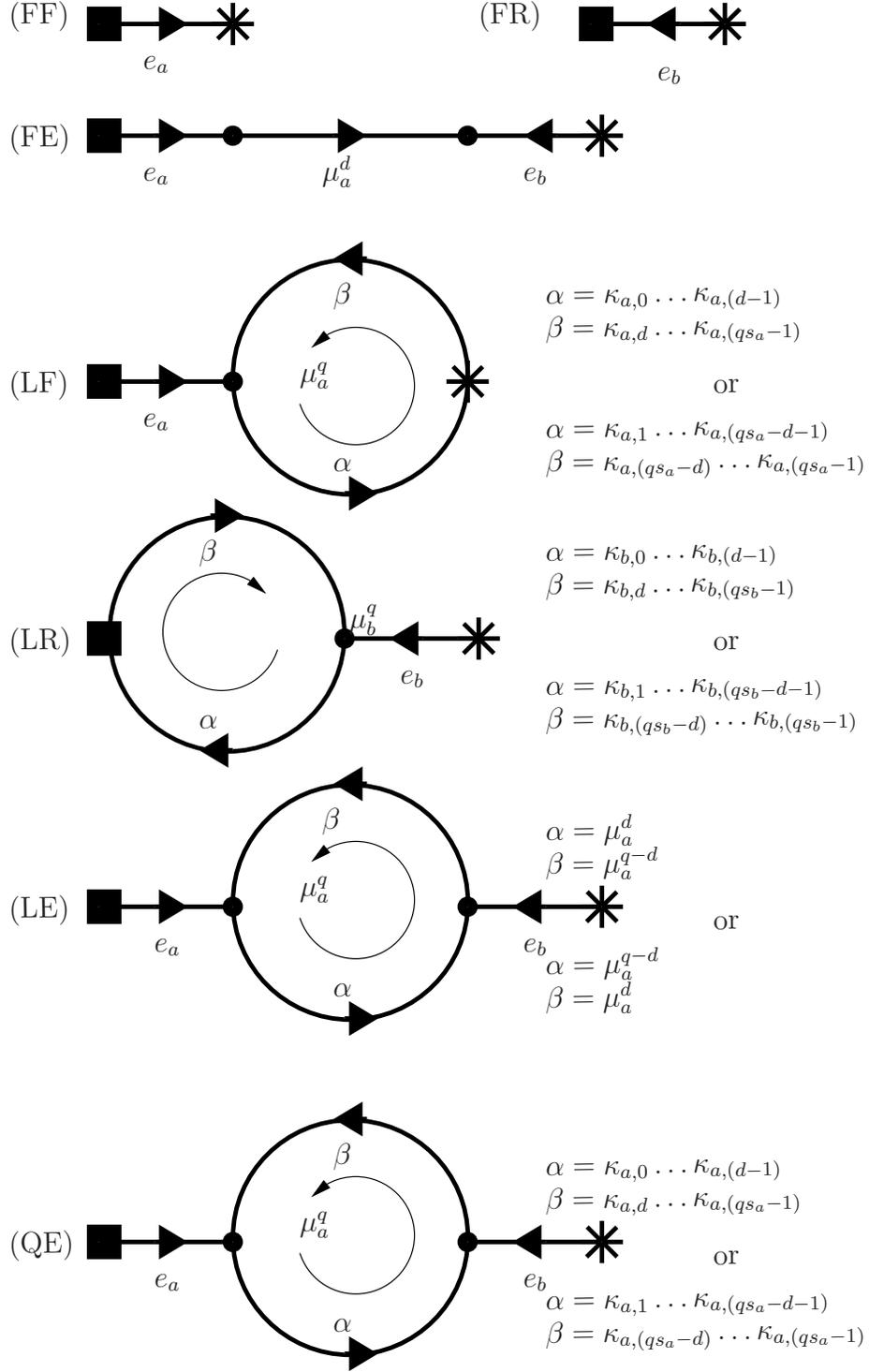}
\caption{The diagram units $\Lambda[\delta, q]$ (assuming the
notation of $\S$\ref{GrowthUnitsSection} for $\delta$).
\label{DiagramUnitsFigure}}
\end{figure}

\begin{notation}\label{LinearBallNotation}\index{sB@$B[\,]$}
For an active growth unit $\delta$, we denote by $B[\delta, q]$
the unique subgraph of $\Lambda[\delta, q]$ which is an embedded
circle.
\end{notation}

We record some elementary properties of diagram units which are
easily verified by inspecting Figure \ref{DiagramUnitsFigure}.

\begin{lem}[Properties of diagram units]\label{PropertiesOfDiagramUnits}
For each growth unit $\delta$ the following statements hold:
 \begin{enumerate}
 \item \label{DiagramUnitsCaptureImagesOfIterates}
 for each non-negative integer $k$, $f^{kq}_{\#}(\delta)$ labels a path $\tilde{\delta}_{kq}$ across
 $\Lambda[\delta, q]$;
 \item \label{DiagramUnitsHaveLinearGrowthOnHomology}
 the function $k \mapsto l^{\ab}(\tilde{\delta}_{kq})$ is linear if and only if $\delta$ is an active growth unit.
 \end{enumerate}
\end{lem}

\begin{con}[The $\Lambda$ and $\Sigma$
constructions]\label{ConstructionOfLambda}\index{sLambdaBig1@$\Lambda(\,)$}\index{sLambdaBig2@$\Lambda[\,]$}
Let $\rho = \delta_1 \separate \delta_2 \separate \dots \separate
\delta_s$ be a path in $G$ with linear growth and let $q \in
\Nat$. Write $\Lambda_i:= \Lambda(\delta_i, q)$ for each $i = 1,
2, \dots, s$, and define an end-pointed $G$-labelled graph
$\Lambda(\rho, q) := \combine{\Lambda_1, \Lambda_2, \dots,
\Lambda_s}$. Further, define $\Lambda[\rho, q]$ to be the
end-pointed $G$-immersion determined by $\Lambda(\rho, q)$ and, if
$\rho$ is a closed path in $G$, define $\Sigma(\rho, q)$
(respectively, $\Sigma[\rho, q]$) to be the base-pointed
$G$-labelled graph (respectively, base-pointed $G$-immersion)
determined by $\Lambda(\rho, q)$.
\end{con}

\begin{rem}\label{AlternativeConstruction}\index{sLambdaBig3@$\Lambda'(\,)$}
For a linear path $\rho = \delta_1 \separate \delta_2 \separate
\dots \separate \delta_s$ in $G$ and $q \in \Nat$, we write
$\Lambda'(\rho, q) := \combine{\Lambda[\delta_1, q],
\Lambda[\delta_2, q], \dots, \Lambda[\delta_s, q]}$.  It follows
from the definitions that $\Lambda(\rho, q)$ folds to
$\Lambda'(\rho, q)$ and we may have defined $\Lambda[\rho, q]$ to
be the end-pointed $G$-immersion determined by $\Lambda'(\rho,
q)$.
\end{rem}

The following property of the $\Lambda$ construction follows
immediately from the definitions and the properties of diagram
units (Lemma \ref{PropertiesOfDiagramUnits}).


\begin{lem}\label{LambdaCapturesIteratesOfAlpha}
Let $q \in \Nat$ and let $\rho$ be a linear path in $G$. For each
non-negative integer $k$ there is a unique path across
$\Lambda[\rho, q]$ which is labelled by $f^{kq}_{\#}(\rho)$.
\end{lem}


\begin{defn}\label{DefnOfPrimaryForm}\index{wprimaryform@primary form}
We say that a linear path $\rho = \delta_1 \separate \delta_2
\separate \dots \separate \delta_t$ in $G$ is \emph{in primary
form with respect to $q$} if each $\delta_i$ takes a minimal
length path across $\Lambda[\delta_i, q]$.  That is, the following
conditions are satisfied for each $i = 1, 2, \dots, s$:
\begin{enumerate}
 \item If $\delta_i$ is of type (LF) then
 $\length{\delta_i} \leq (q \abs{\mu_a})/2 + 1$;
 \item If $\delta_i$ is of type (LR) then
 $\length{\delta_i} \leq (q \abs{\mu_b})/2 + 1$;
 \item If $\delta_i$ is of type (LE) or (QE) then
 $\length{\delta_i} \leq (q \abs{\mu_a})/2 + 2$.
\end{enumerate}
\end{defn}

\begin{notation}\index{slambda0@$\lambda_0$} \index{slambda1@$\lambda_1$}
Define constants $\lambda_0 := lcm \{l(\mu_i) \hbox{ } | \hbox{ }
L_1 \leq i < L_2 \}$ and \begin{equation*}\lambda_1 :=
\max\{l(\mu) \hbox{ } | \hbox{ } \mu \hbox{ a subpath of } \mu_i
\hbox{ for some } L_1 \leq i < L_2 \hbox{ and } \mu \hbox{ an
exceptional path}\}.\end{equation*}
\end{notation}

\begin{lem}\label{C1IsUseful}
Let $q \in \Nat$ be such that $q > \lambda_0$, let $i, j \in \Nat$
be such that $L_1 \leq i, j \leq L_2-1$, write $\mathcal{B}_i :=
C(\mu_i^q)$ and $\mathcal{B}_j := C(\mu_j^q)$, let $v_i \in
\mathcal{B}_i$ and $v_j \in \mathcal{B}_j$ be vertices (not
necessarily the base-points), let $\Delta' := \mathcal{B}_i \amalg
\mathcal{B}_j / \sim$ where $\sim$ equates $v_i$ and $v_j$, let
$\Delta$ be the $G$-immersion determined by $\Delta'$.  The
natural maps $\mathcal{B}_i \to \Delta$, $\mathcal{B}_j \to
\Delta$ are embeddings (we identify $\mathcal{B}_i$ and
$\mathcal{B}_j$ with their respective images under the natural
maps) and exactly one of the following properties holds:
\begin{enumerate}
\item $\Delta$ is $G$-labelled-graph isomorphic to $\mathcal{B}_i$ (and
$\mathcal{B}_j$);
\item $\mathcal{B}_i \cap \mathcal{B}_j$ is a
line-segment of length less than $\lambda_0$.
\end{enumerate}
This statement is illustrated in  Figure
\ref{LinearBalloonsFigure}.
\end{lem}

\begin{proof}
The label on each of $\mathcal{B}_i$ and $\mathcal{B}_j$ is
periodic with period which divides $\lambda_0$. It follows that if
at least $\lambda_0$ edges of $\mathcal{B}_i$ fold with edges of
$\mathcal{B}_j$ then $\mu_i$ is a cyclic permutation of either
$\mu_j$ or $\rev{\mu}_j$.  The result follows.
\end{proof}

\setcounter{figure}{\value{thm}} \stepcounter{thm}
\begin{figure}
\includegraphics{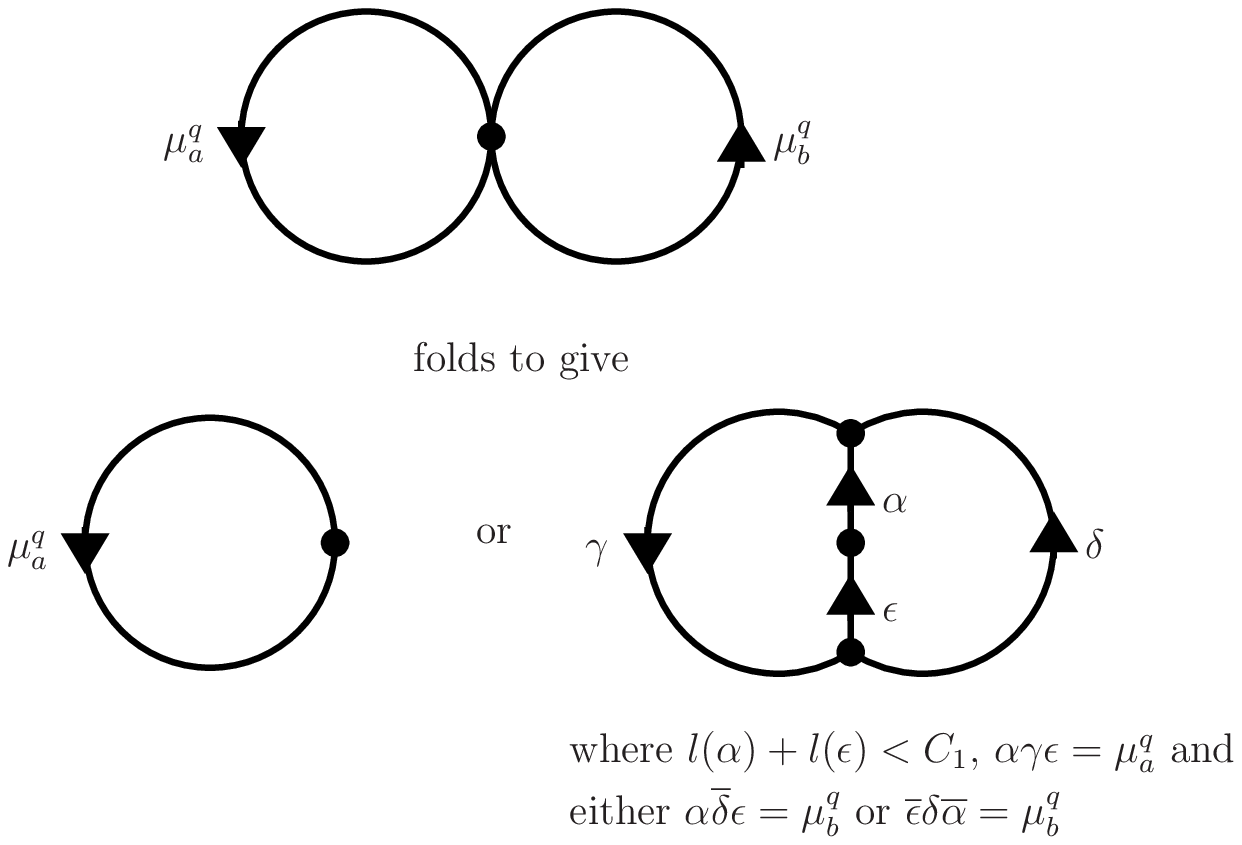}
\caption{Lemma \ref{C1IsUseful}. \label{LinearBalloonsFigure}}
\end{figure}

We are now ready to show that, for sufficiently large $q$,
$\Lambda[\rho, q]$ stretches $\rho$ and the $f^q_{\#}$ iterates of
$\rho$, in the sense of Remark \ref{CarriesAndStretchesRemark}.

\begin{prop}\label{PropertiesOfLambda(PrimaryForm)}
Let $q > 2\max\{\lambda_0, \lambda_1\} + 4\lambda_0$ and let $\rho
= \delta_1 \separate \delta_2 \separate \dots \separate \delta_s$
be a linear path in $G$ which is in primary form with respect to
$q$. The following properties hold:
\begin{enumerate}
 \item \label{PrimaryFormHaveLinearGrowthOnHomology}
 the function $k \mapsto l^{\ab}(\tilde{\rho}_{kq})$ is linear;
 \item \label{BalloonsAreDistinct}
 for each active growth unit $\delta_{i_1}, \delta_{i_2}, \dots,
 \delta_{i_p}$, let $B'_j := B[\delta_{i_j}, q]$ (see Notation
 \ref{LinearBallNotation}) and let $B_j$ be the natural image of
 $B'_j$ in $\Lambda[\rho, q]$.  For each $j = 1, 2, \dots, p$,
 there exists a geometric edge $D_j$ of $\Lambda[\rho, q]$ such that $D_j \subset B_k$ if and
 only if $j = k$.
\end{enumerate}
\end{prop}

\begin{rem}
See Figure \ref{LinearCaseFigure} for a schematic depiction of the
construction of $\Lambda[\rho, q]$ to accompany the argument
below.
\end{rem}

\setcounter{figure}{\value{thm}} \stepcounter{thm}
\begin{figure}
\includegraphics{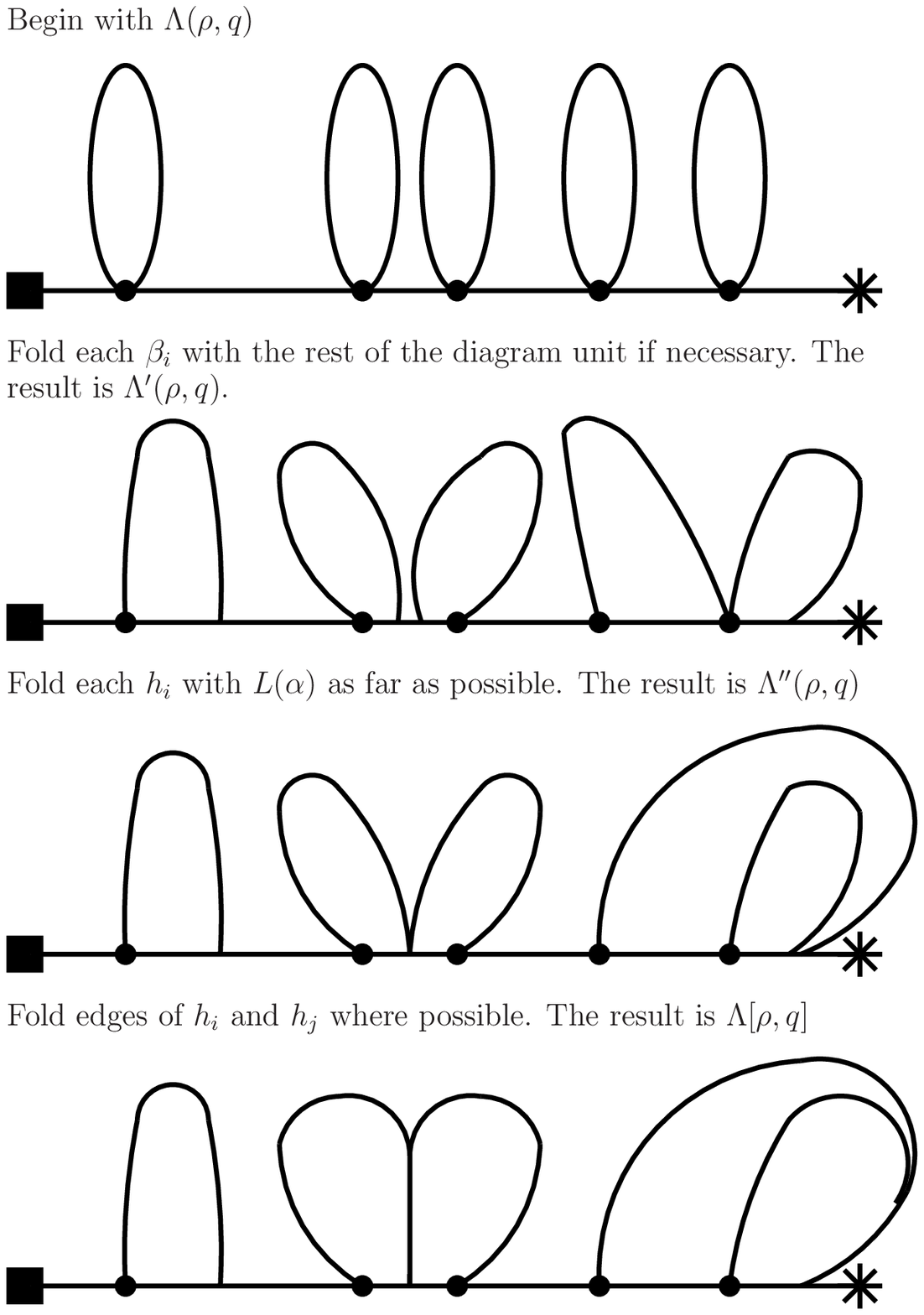}
\caption{A schematic depiction of the construction of
$\Lambda[\rho, q]$. \label{LinearCaseFigure}}
\end{figure}

\begin{proof}
It is clear from the definitions and the properties of diagram
units (in particular Lemma \ref{PropertiesOfDiagramUnits}
(\ref{DiagramUnitsHaveLinearGrowthOnHomology})) that there is a
path $\tilde{\rho}'_{kq}$ across $\Lambda'(\rho, q)$ which is
labelled by $f^{kq}_{\#}(\delta_1) f^{kq}_{\#}(\delta_2) \dots
f^{kq}_{\#}(\delta_s)$ and such that $k \mapsto
l^{\ab}(\tilde{\rho}'_{kq})$ is linear. By Remark
\ref{AlternativeConstruction}, $\Lambda[\rho, q]$ is the immersion
determined by $\Lambda'(\rho, q)$.  It follows that Property
(\ref{BalloonsAreDistinct}) implies Property
(\ref{PrimaryFormHaveLinearGrowthOnHomology}) and it remains only
to show that Property (\ref{BalloonsAreDistinct}) holds.

By the hypothesis that $\rho$ is in primary form with respect to
$q$, we may think of $\Lambda'(\rho, q)$ as consisting of a line
(corresponding to $L(\rho)$) with one handle (the part of
$\Lambda[\delta_i, q]$ not crossed by $\delta_i$) of length at
least $q/2$ attached for each active growth unit. For growth units
of type (LE) and (QE) there is no more folding possible between
the corresponding handle and the line (since $h(\mu_a) < a, b$).
For growth units of type (LF) or (LR) there may be some further
folding possible between the corresponding handle and the line. We
examine the possibilities: suppose that $\delta_i$ is a growth
unit of type (LF) for some integer $i$ such that $1 \leq i \leq s$
and let $h_i$ denote the handle in $\Lambda'(\rho, q)$
corresponding to $\delta_i$. By construction, there is no folding
possible between $h_i$ and the line-segment in $\Lambda'(\rho, q)$
corresponding to $L(\delta_1 \dots \delta_{i-1})$, thus we may
assume that $i < s$. Let $L_{i+1}$ denote the line-segment in
$\Lambda'(\rho, q)$ corresponding to $L(\delta_{i+1} \dots
\delta_s)$. Observe that $\hat{h}_i$ is a Nielsen path. If
$\delta_{i+1}$ is a passive growth unit, less than
$l(\delta_{i+1}) \leq \lambda_1$ edges of $h_i$ may fold with
$L_{i+1}$ (any further folding would violate the maximality
condition in Step 3 of the separation algorithm). If
$\delta_{i+1}$ has type (LE) or (QE), less than $\lambda_1$ edges
of $h_i$ fold with $L_{i+1}$ (otherwise, $\hat{h}_i$ contains an
essentially unbounded subpath $\hat{\delta}_{i+1}$, which is
impossible by Remark \ref{PropertiesOfEssentiallyUnbounded} (1)).
Finally, consider the case that $\delta_{i+1}$ has type (LR).
Since $\hat{h}_i$ is a Nielsen path, less than $l(\delta_{i+1})$
edges of $h_i$ fold with $L_{i+1}$ (otherwise, again, $\hat{h}_i$
contains an essentially unbounded subpath $\hat{\delta}_{i+1}$).
It follows from Lemma \ref{C1IsUseful} and the definition of
growth units that, if $\lambda_0$ edges of $h_i$ fold with
$L_{i+1}$ then $\delta_i \delta_{i+1}$ is a growth unit of type
(FE), (LE) or (QE), contradicting Step 1 of the separation
algorithm. Hence less than $\lambda_0$ edges of $h_i$ fold with
$L_{i+1}$. A similar examination of the possibilities in the case
that $\delta_i$ has type (LR) allows us to conclude the following:
after performing all folding possible between the handles and
$L(\alpha)$ in $\Lambda'(\rho, q)$, the resulting graph
$\Lambda''(\rho, q)$ may be viewed as a line with one handle
attached for each active growth unit, where each such handle has
length at least $q/2 - \max\{\lambda_1, \lambda_0\}$.

There may be some folding possible between handles in
$\Lambda''(\rho, q)$. Suppose at least $\lambda_0$ edges of $h_i$
are identified with edges of $h_j$ in $\Lambda[\rho, q]$, for some
integers $i$ and $j$ such that $1 \leq i < j \leq s$. By Lemma
\ref{C1IsUseful}, $B(\delta_i, q)$ and $B(\delta_j, q)$ are
identified in $\Lambda[\rho, q]$ and hence $\delta_i \delta_{i+1}
\dots \delta_j$ is a growth unit of type (FE),(LE) or (QE)) --- a
contradiction to the separation algorithm. Thus we have that less
than $\lambda_0$ pairs of edges may fold between any two handles
in $\Lambda''(\rho, q)$ and it follows from our hypothesis on $q$
that each handle contains at least one geometric edge which does
not fold with any other handle. Thus Property
(\ref{BalloonsAreDistinct}) holds.
\end{proof}

\begin{thm}\label{AptImmersionTrueForLinear} \index{wAptImmersionTheorem@Apt Immersion Theorem!wproofinlinearcase@proof in the linear case}
The Apt Immersion Theorem holds in the case that $\sigma$ is a
circuit with linear growth.
\end{thm}

\begin{proof}
Let $\rho$ be a well-chosen closed tight path representing
$\sigma$. Choose $q$ such that $q
> 2\max\{\lambda_1, \lambda_0\} + 4\lambda_0$ and $q
> 2 l(\rho)$ so that we may apply Lemma \ref{LambdaCapturesIteratesOfAlpha} and Proposition
\ref{PropertiesOfLambda(PrimaryForm)}
(\ref{PrimaryFormHaveLinearGrowthOnHomology}) to $\Lambda[\rho,
q]$.  It follows immediately that $\Sigma[\rho, q]$ satisfies
(AI1) and (AI2).
\end{proof}

\begin{rem}
Although it would be unjustifiably distracting to develop the
necessary ideas here, it can be shown that there exists a constant
$q_0 \in \Nat$ such that we may replace the hypothesis in
Proposition \ref{PropertiesOfLambda(PrimaryForm)} that $\rho$ is
in primary form with respect to $q$ by the hypothesis that $q \geq
q_0$.
\end{rem}

\section{The Apt Immersion Theorem in the non-linear case}\label{PathUnitsChapter}

In the linear growth case we used growth units to write our path
as a concatenation of subpaths which interact in a limited way
under iteration of $f_{\#}$; in the non-linear case we introduce
the notion of `path units'.   Path units are analogous to basic
paths but more flexible.  A path $\rho$ of degree $d \geq 2$
$f$-splits canonically into path units of degree $d$ and paths of
degree at most $d-1$, and hence we write $\rho$ as a concatenation
of subpaths which do not interact at all under iteration of
$f_{\#}$. We assign to $\rho$ a description, called the `path unit
structure', which summarises this concatenation. Importantly, the
path unit structure of a path is invariant under the action of
$f_{\#}$. Recognising this allows us to construct different
end-pointed $G$-immersions tailored for the different path units
of $\rho$, which can then be combined using our standard
constructions.  By taking care at the neighbourhood of each
end-point of the $G$-immersions constructed, we may ensure that no
folding will be required when the $G$-immersions are combined.  In
this way we reduce the task of proving the Apt Immersion Theorem
to
the task of constructing $G$-immersions which carry and stretch
path units (in the sense of Remark
\ref{CarriesAndStretchesRemark}) and their $f^q_{\#}$ iterates and
for which the neighbourhoods of the end-points are appropriately
simple. We perform the necessary construction inductively,
inducting on the degree $d$ and making repeated use of Stallings'
Algorithm and the structure of the IRTT representative $f:G \to
G$.

\subsection{Path units}

The following definition should be compared to that of a basic
path ($\S$\ref{IRTTSubsection}).

\begin{defn}[Path Units]\label{PathUnitsDefn} \index{wpathunit@path unit}
Let $n \geq 2$ be an integer. A \emph{path unit (of degree
$d$)}\index{wdegree@degree!wofapathunit@of a path unit}
\index{wpathunit@path unit!wdegreeofapathunit@degree} is a path
$\alpha$ in $G$ in one of the following forms:
\begin{enumerate}
\item [(i)] $e_a \gamma \rev{e}_b$;
\item [(ii)] $e_a \gamma$; or
\item [(iii)] $\gamma \rev{e}_b$,
\end{enumerate} where $a, b \in\Nat$ are such that $L_d \leq a, b < L_{d+1}$ and $\gamma$ is a path in $G_{L_d -
1}$.  A path unit which has form $(t)$ for some $t \in \{i, ii,
iii\}$ is said to \emph{have type
(t)}\index{wtype@type!wofapathunit@of a path
unit}\index{wpathunit@path unit!wtypeofapathunit@type}.
\end{defn}


\begin{defn}[Canonical $f$-splitting of a
path]\label{CanonicalSplittingOfPath}\index{wcanonicalsplitting@canonical
splitting } Let $\rho$ in $G$ be a path. We define the
\emph{canonical $f$-splitting of $\rho$} as follows: if $h(\rho) <
L_1$, then the canonical $f$-splitting is simply the path $\rho$;
if $L_1 \leq h(\rho) < L_2$, then $f$-split $\rho$ immediately
before each occurrence of $e_{h(\rho)}$ and after each occurrence
of $\rev{e}_{h(\rho)}$; if $L_d \leq h(\rho) < L_{d+1}$ for some
$d \geq 2$, then, for each integer $i$ such that $L_d \leq i <
L_{d+1}$, $f$-split $\rho$ immediately before each occurrence of
$e_{i}$ and after each occurrence of $\rev{e}_{i}$.
\end{defn}

\begin{notation}\index{saasterisk@$\ast$}
We denote that an $f$-splitting $\rho = \alpha_1 \cdot \alpha_2
\cdot \dots \cdot \alpha_s$ is in fact the canonical $f$-splitting
by using the symbol $\ast$ between subpaths rather than the symbol
$\cdot$, that is, we write $\rho = \alpha_1 \ast \alpha_2 \ast
\dots \ast \alpha_s$.
\end{notation}

\begin{rem}
If $\rho$ has degree $d \geq 2$, the canonical $f$-splitting
writes $\rho$ as a concatenation of maximal path units of degree
$d$ and paths of degree at most $d-1$.  If $\rho$ is linear, the
canonical $f$-splitting writes $\rho$ as a concatenation of
maximal basic paths of height $h(\rho)$ and paths of height at
most $h(\rho)-1$.
\end{rem}


The following lemma is an immediate consequence of Property (TT4)
of the IRTT Theorem.

\begin{lem}\label{InverseImageOfPathUnit}
Let $\alpha \subset G$ be a path unit.  If $\alpha$ has type (i)
then $f_{\#}(\alpha)$ and $\overline{f}_{\#}(\alpha)$ are path
units of type (i) with the same initial and terminal edges as
$\alpha$; if $\alpha$ has type (ii) then $f_{\#}(\alpha)$ and
$\overline{f}_{\#}(\alpha)$ are path units of type (ii) with the
same initial edge as $\alpha$; if $\alpha$ has type (iii) then
$f_{\#}(\alpha)$ and $\overline{f}_{\#}(\alpha)$ are path units of
type (iii) with the same terminal edge as $\alpha$.
\end{lem}

\begin{rem}\label{PathUnitStructure}
Let $\alpha \subset G$ be a path unit.  We define the
\emph{structure}\index{wstructure@structure of a path unit} of
$\alpha$, denoted $\structure(\alpha)$ as follows: if $\alpha =
e_a \gamma \rev{e}_b$ has type (i) then $\structure(\alpha) =
((i), a, b)$; if $\alpha = e_a \gamma$ has type (ii) then
$\structure(\alpha) = ((ii), a)$; if $\alpha = \gamma \rev{e}_b$
has type (iii) then $\structure(\alpha) = ((iii), b)$. For a path
$\rho = \alpha_1 \ast \alpha_2 \ast \dots \ast \alpha_s$ in $G$ of
degree $d$, define the \emph{path unit structure} of $\rho$ to be
a finite list of sets $s_1, s_2, \dots, s_s$, where $s_i =
\emptyset$ if degree $\alpha_i < d$, otherwise $s_i :=
\structure(\alpha_i)$. It is an immediate corollary to Lemma
\ref{InverseImageOfPathUnit} that the path unit structure of a
path is invariant under the action of $f_{\#}$ and
$\overline{f}_{\#}$.
\end{rem}

\subsection{Tails of edges}

Let $a \geq 1$ be an integer.  We now investigate the structure of
the `$f$-tails' of $e_a$, that is, the infinite paths
$S_a^+$\index{sSa+@$S_a^+$} and $S_a^-$\index{sSa-@$S_a^-$} such
that $f^k_{\#}(e_a) \to e_a S_a^+$ and $\overline{f}^{k}_{\#}(e_a)
\to e_a S_a^-$ as $k \to \infty$. An understanding of this
structure is crucial for building $G$-immersions which carry a
path unit which crosses $E_a$.

For each integer $a$ such that $L_1 \leq a < h(G)$, define an
infinite tight path $S^+_a := u_a f_{\#}(u_a) f^2_{\#}(u_a)
\dots$.  We define a second infinite tight path $S^-_a$ in one of
two ways, depending on whether $u_a$ is a well-chosen closed tight
path or not:

\vskip 10pt

\noindent In the case that $u_a$ is a well-chosen closed tight
path, let $u_a = \epsilon_{a, \, 0} \ast \epsilon_{a, \, 1} \ast
\dots \ast \epsilon_{a, s_a-1}$.  Note that
\begin{equation*}S^+_a = \epsilon_{a, \, 0} \ast \epsilon_{a, \,
1} \ast \dots \ast \epsilon_{a, s_a-1} \ast f_{\#}(\epsilon_{a, \,
0}) \ast f_{\#}(\epsilon_{a, \, 1}) \ast \dots \ast
f_{\#}(\epsilon_{a, s_a-1}) \ast \dots.\end{equation*} Define an
infinite tight path,
\begin{equation*}S^-_a := \overline{f}_{\#}(\rev{\epsilon}_{a, \, s_a-1}) \ast \dots
 \ast \overline{f}_{\#}(\rev{\epsilon}_{a, \, 0}) \ast \overline{f}^{2}_{\#}(\rev{\epsilon}_{a, \,
s_a-1}) \ast \dots \ast \overline{f}^{2}_{\#}(\rev{\epsilon}_{a,
\, 0}) \ast \dots.\end{equation*}

\vskip 10pt

\noindent In the case that $u_a$ is {\bf not} a well-chosen closed
tight path, let $u_a = \epsilon'_{a, \, 0} \ast \epsilon_{a, \, 1}
\ast \dots \ast \epsilon_{a, s_a-1} \ast \epsilon'_{a, \, s_a}$
and define $\epsilon_{a, \, 0} := [\overline{f}_{\#}(\epsilon'_{a,
\, s_a}) \epsilon'_{a, \, 0}]$. Note that
\begin{equation*}S^+_a = \epsilon'_{a, \, 0} \ast \epsilon_{a, \, 1} \ast \dots \ast
\epsilon_{a, s_a-1} \ast f_{\#}(\epsilon_{a, \, 0}) \ast
f_{\#}(\epsilon_{a, \, 1}) \ast \dots \ast f_{\#}(\epsilon_{a,
s_a-1}) \ast \dots.\end{equation*} Define an infinite tight path,
\begin{equation*}S^-_a := \overline{f}_{\#}(\rev{\epsilon}'_{a, \, s_a}) \ast \overline{f}_{\#}(\rev{\epsilon}_{a, \, s_a-1}) \ast \dots
 \ast \overline{f}_{\#}(\rev{\epsilon}_{a, \, 0}) \ast \overline{f}^{2}_{\#}(\rev{\epsilon}_{a, \,
s_a-1}) \ast \dots \ast \overline{f}^{2}_{\#}(\rev{\epsilon}_{a,
\, 0}) \ast \dots.\end{equation*}

In either case, observe that
\begin{equation*}
\begin{split}
[\rev{S_a^-} S_a^+] = \dots \ast \overline{f}_{\#}(\epsilon_{a, \,
0}) \ast \dots & \ast \overline{f}_{\#}(\epsilon_{a, \, s_a-1})
\ast \epsilon_{a, \, 0} \ast \dots \\ & \ast \epsilon_{a, \,
s_a-1} \ast f_{\#}(\epsilon_{a, \, 0}) \ast \dots \ast
f_{\#}(\epsilon_{a, \, s_a-1}) \ast \dots.
\end{split}
\end{equation*}
\begin{notation}\index{salphaai@$\alpha_{a, \, i}$}\index{sbetaai@$\beta_{a, \,
i}$} For each integer $a$ such that $L_1 \leq a < h(G)$, relabel
the canonical $f$-splitting $S_a^+ = \alpha_{a, \, 0} \ast
\alpha_{a, \, 1} \ast \alpha_{a, \, 2} \ast \dots$ and $S_a^- =
\beta_{a, \, 0} \ast \beta_{a, \, 1} \ast \beta_{a, \, 2} \ast
\dots$.
\end{notation}

\begin{rem} \label{SemiPeriodicityOfTails}
For each $a \geq L_2$ and each integer $i \geq 1$, $\alpha_{a, \,
i + s_a} = f_{\#}(\alpha_{a, \, i})$ and $\beta_{a, \, i + s_a} =
\overline{f}_{\#}(\beta_{a, \, i})$.
\end{rem}

Recall, in $\S$\ref{NotationSection} we defined $\mu_i$ to be the
primitive closed path corresponding to $u_i$, for each integer $i$
such that $L_1 \leq i < L_2$.

\begin{lem}[The Linear Balloon Lemma]\label{LinearBalloonLemma}\index{wLinearBalloonLemma@Linear Balloon Lemma}
Let $a, b \in \Nat$ be such that $L_1 \leq a, b < L_2$ and define
$K := l(\mu_a)l(\mu_b)$.  The following properties hold:
\begin{enumerate}
\item if there exist finite tight paths $\rho_1,
\rho'_1, \rho_2$ and infinite tight paths $\rho_3$ and $\rho'_3$
such that $S_a^+ = \rho_1 \rho_2 \rho_3$, $S_b^+ = \rho'_1 \rho_2
\rho'_3$ and $l(\rho_2) \geq K$, then $\mu_a = \mu_b$.
\item if there exist finite tight paths $\rho_1,
\rho'_1, \rho_2$ and infinite tight paths $\rho_3$ and $\rho'_3$
such that $S_a^+ = \rho_1 \rho_2 \rho_3$, $S_b^- = \rho'_1 \rho_2
\rho'_3$ and $l(\rho_2) \geq K$, then $\mu_a$ is a cyclic
permutation of $\rev{\mu}_b$.
\end{enumerate}
\end{lem}

\begin{proof}
Assume the hypothesis of Property (2).  The infinite paths $S_a^+,
S_a^-$ are periodic with period $l(\mu_a)$ and the infinite paths
$S_b^+, S_b^-$ are periodic with period $l(\mu_b)$.  The
periodicity of $S_a^+$ and $S_b^-$ imply that $l(\mu_a) =
l(\mu_b)$ and $\mu_a$ is a cyclic permutation of $\rev{\mu}_b$.
That is, Property (2) holds.

Assume the hypothesis of Property (1).  As above, the periodicity
of $S_a^+$ and $S_b^+$ imply that $\mu_a$ is a cyclic permutation
of $\mu_b$. Suppose $\mu_a \neq \mu_b$, say $\mu_a = \epsilon
\mu_b \rev{\epsilon}$. Then $e_a \epsilon \rev{e}_b$ is a tight
path which violates (TT4) of the IRTT theorem.  Thus $\mu_a =
\mu_b$ and Property (1) holds.
\end{proof}

Lemma \ref{LinearBalloonLemma} is a consequence of the periodicity
of $S^+_a$ and $S^-_a$ in the case that $L_1 \leq a < L_2$.  In
the case that $a \geq L_2$, $S^+_a$ and $S^-_a$ are not periodic,
but Remark \ref{SemiPeriodicityOfTails} can be used to mimic the
role of periodicity.

\begin{lem}[The Non-Linear Balloon
Lemma]\label{BalloonsDontFoldMuch}\index{wNonLinearBalloonLemma@Non-Linear
Balloon Lemma} Let $d \geq 2$, let $a, b \in \Nat$ be such that
$L_d \leq a, b < L_{d+1}$ and define $K := s_a s_b + \min \{s_a,
s_b\} + 1$.  The following properties hold:
\begin{enumerate}
\item if there exist finite tight paths $\rho_1, \rho'_1, \rho_2$ and infinite
tight paths $\rho_3$ and $\rho'_3$ such that $S_a^+ = \rho_1
\rho_2 \rho_3$, $S_b^+ = \rho'_1 \rho_2 \rho'_3$ and the
separation of $\rho_2$ contains at least $K$ complete path units
of the canonical $f$-splitting of $S_a^+$, then $a = b$, $\rho'_1
= \rho_1$ and $\rho'_3 = \rho_3$.
\item if there exist finite tight paths $\rho_1, \rho'_1, \rho_2$ and infinite
tight paths $\rho_3$ and $\rho'_3$ such that $S_a^- = \rho_1
\rho_2 \rho_3$, $S_b^- = \rho'_1 \rho_2 \rho'_3$ and the
separation of $\rho_2$ contains at least $K$ complete path units
of the canonical $f$-splitting of $S_a^-$, then $a = b$, $\rho'_1
= \rho_1$ and $\rho'_3 = \rho_3$.
\item if there exist finite tight paths $\rho_1, \rho'_1, \rho_2$ and infinite
tight paths $\rho_3$ and $\rho'_3$ such that $S_a^+ = \rho_1
\rho_2 \rho_3$ and $S_b^- = \rho'_1 \rho_2 \rho'_3$, then the
separation of $\rho_2$ contains less than $K$ complete path units
of the canonical $f$-splitting of $S_a^+$.
\end{enumerate}
\end{lem}

\begin{proof}
We first claim that the hypothesis of Property (1) implies that
$s_a = s_b$.  The hypothesis implies that there exist $i, j \geq
1$ such that $\alpha_{a, \, i} \dots \alpha_{a, \, i+K-1} =
\alpha_{b, \, j} \dots \alpha_{b, \, j+K-1}$. Thus $\alpha_{a, \,
i + s_a s_b + j} = f^{s_b}_{\#}(\alpha_{a, \, i+j}) =
f^{s_a}_{\#}(\alpha_{a, \, i+j})$, for $j = 0$, $1$, $\dots$,
$\min\{s_a, s_b\} -1$ (by Remark \ref{SemiPeriodicityOfTails}).
But at least one of the path units $\alpha_{a, i +j}$ is not a
Nielsen path, say $j = j_0$, and $f^{s_b}_{\#}(\alpha_{a, \,
i+j_0}) = f^{s_a}_{\#}(\alpha_{a, \, i+j_0})$ implies that $s_a =
s_b$.

Now, suppose that the path $\alpha = [e_a \alpha_{a, \, 0} \dots
\alpha_{a, \, i-1} \rev{\alpha}_{b, \, j-1} \dots \rev{\alpha}_{b,
\, 0} \rev{e}_b]$ is not the trivial path.  We may assume that
$\alpha = e_a \alpha_{a, \, 0} \dots [\alpha_{a, \, i-1}
\rev{\alpha}_{b, \, j-1}] \dots \rev{\alpha}_{b, \, 0} \rev{e}_b$.
Then
\begin{eqnarray}
f_{\#}(\alpha) & = & [e_a \alpha_{a, \, 0} \dots \alpha_{a, \,
i-1} \alpha_{a, \, i} \dots \alpha_{a, \, i+s_a-1}
\rev{\alpha}_{b, \, j + s_b -1} \dots \rev{\alpha}_{b, \, j}
\rev{\alpha}_{b, \, j-1} \dots \rev{\alpha}_{b, \, 0} \rev{e}_b]
\nonumber \\
& = & e_a \alpha_{a, \, 0} \dots [\alpha_{a, \, i-1} \alpha_{a, \,
i} \dots \alpha_{a, \, i+s_a-1} \rev{\alpha}_{a, \, i + s_a -1}
\dots \rev{\alpha}_{a, \, i } \rev{\alpha}_{b, \, j-1}] \dots
\rev{\alpha}_{b, \, 0} \rev{e}_b
\nonumber \\
& = & e_a \alpha_{a, \, 0} \dots [\alpha_{a, \, i-1}
\rev{\alpha}_{b, \, j-1}] \dots \rev{\alpha}_{b, \, 0} \rev{e}_b.
\nonumber
\end{eqnarray}
Hence $\alpha$ is a Nielsen path which crosses a non-linear edge
--- a contradiction to the Corollary \ref{IRTTRepsAndGrowth}
(\ref{ExceptionalPathsLinearEdges}). Thus $\alpha$ is the trivial
path, and the conclusions of Property (1) hold.

Property (2) may be proved by a similar argument to the above.

Assume the hypothesis of Property (3).  Suppose that the
separation of $\rho_2$ contains at least $K$ complete path units
of the canonical $f$-splitting of $S_a^+$.  By hypothesis, there
exist $i, j \in \Nat$ such that $\alpha_{a, \, i} \dots \alpha_{a,
\, i+K-1} = \beta_{b, \, j} \dots \beta_{b, \, j+K-1}$. Thus
$\alpha_{a, \, i + s_a s_b + j} = f^{s_b}_{\#}(\alpha_{a, \, i+j})
= \overline{f}^{s_a}_{\#}(\alpha_{a, \, i+j})$, for $j = 0, 1,
\dots, \min\{s_a, s_b\} -1$.  This implies that each $\alpha_{a,
\, i+j}$ is a Nielsen path, a contradiction to the fact that $a
\geq L_d \geq L_2$. Hence Property (2) holds.
\end{proof}

\begin{lem}\label{TailsContainEssentiallyUnboundedSubpaths}
For each $a \in \Nat$ such that $a \geq L_2$ and for each
non-negative integer $k$, the infinite paths $\alpha_{a, \, k}
\alpha_{a, \, k+1} \dots$ and $\beta_{a, \, k} \beta_{a, \, k+1}
\dots$ each contain an essentially unbounded subpath.
\end{lem}

\begin{proof}
Immediate by Lemma \ref{PropertiesOfSeparation}.
\end{proof}

\subsection{A strategy of proof}\index{wAptImmersionTheorem@Apt Immersion Theorem!wstrategyofproof@strategy of proof}

In this section we indicate our strategy for completing the proof
of the Apt Immersion Theorem.

\begin{lem}\label{CanOpenSinglePathUnits}
Let $\alpha \subset G$ be a path unit. There exists an end-pointed
$G$-immersion $\Gamma$ with the following properties:
\begin{enumerate}
\item $\iota(\Gamma) \neq \tau(\Gamma)$;
\item if $\alpha$ has type (i) then both $\iota(\Gamma)$ and
$\tau(\Gamma)$ have valence 1; if $\alpha$ has type (ii) then
$\iota(\Gamma)$ has valence 1;  if $\alpha$ has type (iii) then
$\tau(\Gamma)$ has valence 1;
\item $\Gamma$ is $f$-stable;
\item $\alpha$ labels a path across $\Gamma$.
\end{enumerate}
\end{lem}

\begin{proof}
Consider the case that $\alpha$ has type (i), say, $\alpha = e_a
\gamma \rev{e}_b$. Let $d$ be such that $L_d \leq h(\alpha) <
L_{d+1}$ and let $G'$ be the connected component of $G_{L_d-1}$
containing $\gamma$. Extend $L(\gamma)$ to a $G'$-cover $P$.
Define the initial and terminal points of $P$ to be those
corresponding to the initial and terminal points of $L(\gamma)$
respectively. Define $\Gamma := \combine{(L(e_a), P
,L(\rev{e}_b)}$.  Properties (1), (2) and (4) are immediate by
construction.   Property (3) follows easily from the construction
and the fact that $\mu_a, \mu_b \subset G'$. An upper bound on
$period(\Gamma)$ is given by $\abs{\mathcal{V}_\Gamma}!$. The
cases that $\alpha$ has type (ii) and type (iii) are proved
similarly.
\end{proof}

\begin{proof}[{\bf Sketch of a proof of the Apt Immersion Theorem in the non-linear case}]\label{HowToFinishProof}
Let $\sigma \subset G$ be a circuit of degree $d$ and let $\rho =
\alpha_1 \ast \alpha_2 \ast \dots \ast \alpha_s$ be a well-chosen
closed tight path corresponding to $\sigma$.  For each $i = 1, 2,
\dots, s,$ there exists an end-pointed $G$-immersion $\Gamma_i$
which satisfies the conclusions of Lemma
\ref{CanOpenSinglePathUnits} for $\alpha_i$.  For each $ i = 1, 2,
\dots, s$, let $q_i = \period(\Gamma_i)$.  Define $q = \lcm\{q_i
\hbox{ } | \hbox{ } 1 \leq i \leq s \}$ and $\Gamma :=
\combineC{\Gamma_1, \dots, \Gamma_s}.$ Properties (1) and (2) of
Lemma \ref{CanOpenSinglePathUnits} imply that $\Gamma$ is a
$G$-immersion.  Properties (3) and (4) of Lemma
\ref{CanOpenSinglePathUnits} imply that, for each non-negative
integer $k$, $f^{kq}_{\#}(\rho)$ labels a path across $\Gamma$.
Further, since $f^{kq}_{\#}(\rho) = f^{kq}_{\#}(\alpha_1) \ast
f^{kq}_{\#}(\alpha_2) \ast \dots \ast f^{kq}_{\#}(\alpha_s)$ for
each non-negative integer $k$, $f^{kq}_{\#}(\alpha_i)$ labels a
path across $\Gamma_i \subset \Gamma$. Since $\rho$ is a
well-chosen closed tight path, either the initial edge of $\rho$
is $e_{h(\rho)}$ or the terminal edge of $\rho$ is
$\rev{e}_{h(\rho)}$. Without loss of generality, we may assume
that the former case holds. Let $\tilde{\alpha}_{1, \, k}$ be the
path across $\Gamma_1$ labelled by $f^{k q_1}_{\#}(\alpha_1)$. If
we can ensure that $\Gamma_1$ is such that $k \mapsto
l^{\ab}(\tilde{\alpha}_{1, \, k})$ is an element of $p_d$, then
$\Gamma$ and $q$ will satisfy the conclusions of the Apt Immersion
Theorem (and the proof of the Main Theorem will be complete).
\end{proof}

\subsection{The Periodic Open Immersions Lemma}

Next we use a simple finiteness argument to find, for a path unit
$\alpha$ in $G$, a periodic sequence of end-pointed $G$-immersions
$\{\Gamma_i\}_{i \in \Nat}$ with Properties (1) and (2) of Lemma
\ref{CanOpenSinglePathUnits} and such that $\Gamma_i$ carries
$f^i_{\#}(\alpha)$. The periodicity of the sequence is crucial
because it allows us to join finite strings of such end-pointed
$G$-immersions into base-pointed $G$-immersions such that, for
some integer $a \geq L_2$ and some $q \in \Nat$, the result
carries $S^+_a$ and $S^-_a$.

\begin{lem}\label{TheOpenLemma} Let
$\alpha \subset G$ be a path unit of degree $d \geq 2$. There
exist $m, r \in \Nat$ and a finite set of end-pointed
$G$-immersions $\{\Gamma_0, \dots, \Gamma_{r-1}\}$ such that, for
each $i = 0, 1, \dots, r-1$ and each $j = 0, 1, \dots, m-1$, the
following properties hold:
\begin{enumerate}
\item $\iota(\Gamma_i) \neq \tau(\Gamma_i)$;
\item if $\alpha$ has type (i) then both $\iota(\Gamma_i)$ and
$\tau(\Gamma_i)$ have valence 1; if $\alpha$ has type (ii) then
$\iota(\Gamma_i)$ has valence 1;  if $\alpha$ has type (iii) then
$\tau(\Gamma_i)$ has valence 1;
\item $\Gamma_i$ is $f$-stable and $\period(\Gamma_i)$ divides $m$;
\item there exist $i_j \in \{0, 1, \dots,
r-1\}$ such that $f^j_{\#}(\alpha)$ labels a path across
$\Gamma_{i_j}$.
\end{enumerate}
\end{lem}

\begin{proof}
Let $a = h(\alpha)$, let $d$ be such that $L_d \leq h(\alpha) <
L_{d+1}$ and assume, without loss of generality, that the initial
edge of $\alpha$ is $e_a$ (that is, $\alpha$ has type (i) or type
(ii)). Let $G'$ be the connected component of $G_a$ which contains
$e_a$ and let $H = \pi_1(G', \iota(e_a))$. If $\alpha$ has type
(ii), the result holds with $r =1$, $m = 1$ and $(\Gamma_0, p_0)$
constructed from $G'$ simply by detaching the initial point of
$e_a$. Thus we may assume that $\alpha$ has type (i), that is,
$\alpha = e_a \gamma \rev{e}_b$ for some $b \in \Nat$ such that
$L_d \leq b \leq a$ and for some $\gamma \subset G_{L_d -1}$. If
$a \neq b$, the result holds with $r =1$, $m = 1$ and $(\Gamma_0,
p_0)$ constructed from $G'$ simply by detaching the initial point
of $e_a$ and the initial point $e_b$. Thus we may assume that $a =
b$.

Because free groups are residually finite, there exists a
finite-index subgroup $H_0 \leq H$ such that $\hat{\alpha} \not
\in H_0$; write $I := [H: H_0]$. Now $H$ is $\phi$-invariant and
$\phi|_H$ permutes the set of subgroups of $H$ of index $I$.  Let
$\{H_0, H_1, \dots, H_{s-1}\}$ be the $\phi|_H$-orbit of $H_0$
(indexed such that $(\phi|_H)^k(H_0) = H_{k \, \modulo \, s}$).
For each $i = 0, 1, \dots, s-1$, let $(\Delta'_i, p'_i)$ be a
$G'$-cover and let $\tilde{b}_i \in \Delta_i$ be a vertex such
that ${p'_i}(\tilde{b}_i) = \iota(e_a)$ and ${p'_i}_{\ast}
\pi_1(\Delta'_i, \tilde{b}_i)$ corresponds to $H_i$.  Let
$(\Delta_i, p_i)$ be constructed from $(\Delta'_i, p'_i)$ by
detaching each edge with label $e_a$ at the initial point. Let
$\tilde{c}_i^0, \tilde{c}_i^1, \dots, \tilde{c}_i^{I-1}$ be an
enumeration of the vertices in $(p_i)^{-1}(\iota(e_a))$ such that
$\tilde{c}_i^0$ corresponds to $\tilde{b}_i$. For each $j = 1, 2,
\dots, I-1$, let $(\Delta_i^j, p_i^j)$ be the end-pointed
$G_i$-labelled graph $(\Delta_i, p_i)$ with initial point
$\tilde{c}_i^0$ and terminal point $\tilde{c}_i^j$.

Now, $\alpha \not \in H_0$ implies that, for each non-negative
integer $j$, $f^j_{\#}(\alpha) \not \in H_{j \, \modulo \, s}$.
Hence $f^j_{\#}(\alpha)$ lifts to an open path in $\Delta'_{j \,
\modulo \, s}$ at $\tilde{b}_{j \, \modulo \, s}$. It follows that
$f^k j_{\#}(\alpha)$ labels a path across $\Delta_{j \, \modulo \,
s}^l$ for some integer $l$ such that $1 \leq l \leq I-1$. Since
any lift of $(f|_{G'})^{I!}$ to $\Delta'_i$ fixes each vertex of
$\Delta'_i$, for each non-negative integer $k$, we know that
$f^{j+k(I!)}_{\#}(\alpha)$ labels a path across $\Delta_{j \,
\modulo \, s}^l$. Thus the result holds, with $m = I!$, $r = sI$,
and $\{(\Gamma_0, p_0), \dots, (\Gamma_{r-1}, p_{r-1})\} =
\{(\Delta_i^j, p^i_j) \hbox{ } | \hbox{ } 0 \leq i < s; 1 \leq j <
I\}$.
\end{proof}

The following corollary to Lemma \ref{TheOpenLemma} is obtained by
replacing the {\bf set} of $G$-immersions by an ordered {\bf list}
of $G$-immersions.

\begin{cor}[The Periodic Open Immersions
Lemma]\label{TheUsefulOpenLemma}\index{wPeriodicOpenImmersionsLemma@Periodic
Open Immersions Lemma} Let $\alpha_1, \dots, \alpha_s \subset G$
be a finite ordered list of path units of degree $d \geq 2$. There
exists $q \in \Nat$ and, for each $i = 1, 2, \dots, s$, there
exists a finite list of end-pointed $G$-immersions $\Gamma_i^0,
\dots, \Gamma_i^{q-1},$ such that, for each $j = 0, 1, \dots,
q-1$, the following properties hold:
\begin{enumerate}
\item $\iota(\Gamma_i^j) \neq \tau(\Gamma_i^j)$;
\item if $\alpha_i$ has type (i) then both $\iota(\Gamma_i^j)$ and
$\tau(\Gamma_i^j)$ have valence 1; if $\alpha_i$ has type (ii)
then $\iota(\Gamma_i^j)$ has valence 1;  if $\alpha_i$ has type
(iii) then $\tau(\Gamma_i^j)$ has valence 1;
\item $\Gamma_i^j$ is $f$-stable with $\period(\Gamma_i^j) \, | \, q$; and
\item $f^{j}_{\#}(\alpha_i)$ labels a path across $\Gamma_i^j$.
\end{enumerate}
\end{cor}

\begin{rem}
For each $i = 1, 2, \dots, s$, consider the periodic bi-infinite
sequence of end-pointed $G$-immersions $\{\Gamma_i^k \}_{k \in \,
\Integer}$, where for $k < 0 $ and $k \geq q$, $\Gamma_i^k :=
\Gamma_i^{k \, \modulo \, q}$. Properties (3) and (4) imply that
$\alpha_i$ labels a path across $\Gamma_i^0$ and, for each $k \in
\Nat$, $f^{k}_{\#}(\alpha_i)$ labels a path across $\Gamma_i^k$
and $\overline{f}^{k}_{\#}(\alpha_i)$ labels a path across
$\Gamma_i^{-k}$.
\end{rem}

\subsection{The Apt Immersion Theorem in the quadratic case}

The following completes the proof of the Apt Immersion Theorem in
the quadratic case and hence completes the proof of the Main
Theorem in the case that $\mathcal{G}_{\phi} \in p_2$.

\begin{prop}\label{TheQuadraticPathUnitLemma}\index{wAptImmersionTheorem@Apt Immersion Theorem!wproofinquadraticcase@proof in the quadratic case}
Let $\alpha \subset G$ be a path unit of degree 2.  There exists
an end-pointed $G$-immersion $\Sigma$ such that the following
conditions hold:
\begin{enumerate}
\item [$(\Sigma 1)$] $\iota(\Sigma) \neq \tau(\Sigma)$;
\item [$(\Sigma 2)$] If $\alpha$ has type (i) then both $\iota(\Sigma)$ and
$\tau(\Sigma)$ have valence 1; if $\alpha$ has type (ii) then
$\iota(\Sigma)$ has valence 1;  if $\alpha$ has type (iii) then
$\tau(\Sigma)$ has valence 1;
\item [$(\Sigma 3)$] For each non-negative integer $k$, $f^{kq}_{\#}(\alpha)$
labels a path $\tilde{\alpha}_{kq}$ across $\Sigma$;
\item [$(\Sigma 4)$] $k \mapsto l^{\ab}(\tilde{\alpha}_{kq})$ is an element of $p_2$.
\end{enumerate}
\end{prop}

\begin{proof}
Consider the case that $\alpha$ has type (i) (the hardest case).
Suppose that $\alpha$ has type (i), say $\alpha = e_a \alpha'
\rev{e}_b$, and consider $h(u_a)$, $h(u_b)$ and $h(\alpha')$. We
may assume, without loss of generality, that $h(u_a) \geq h(u_b)$.
The case that $h(u_a) = h(u_b) \geq h(\alpha')$ is Lemma
\ref{TheQuadraticPathUnitLemmaCase1}; the case that $h(u_a)
> h(u_b)$ and $h(u_a) \geq h(\alpha')$ is Lemma
\ref{TheQuadraticPathUnitLemmaCase2}; the case that $h(\alpha')
> h(u_a) = h(u_b)$ is Lemma \ref{TheQuadraticPathUnitLemmaCase3};
the case that $h(\alpha') > h(u_a) > h(u_b)$ is Lemma
\ref{TheQuadraticPathUnitLemmaCase4}.  The proof in the case that
$\alpha$ has type (ii) is performed similarly.  By considering
$\rev{\alpha}$ instead of $\alpha$, it is clear that the case that
$\alpha$ has type (iii) is equivalent to the case that $\alpha$
has type (ii).
\end{proof}

\begin{lem}\label{TheQuadraticPathUnitLemmaCase1}
Proposition \ref{TheQuadraticPathUnitLemma} holds in the case that
$\alpha = e_a \alpha' \rev{e}_b$ and $h(u_a) = h(u_b) \geq
h(\alpha')$.
\end{lem}

\begin{proof} (An example construction of $\Sigma$ as below is illustrated
schematically in Figure \ref{QuadraticCase1Figure}). Define $h :=
h(u_a)$, $K := s_a s_b + \min \{s_a, s_b\} + 1$ and $l :=
l(\alpha') + 2K+1$. Let $U_a^+$ (respectively, $U_a^-, U_b^+,
U_b^-$) be the initial subpath of $S_a^+$ (respectively, $S_a^-,
S_b^+, S_b^-$) consisting of the first $l$ path units in the
canonical $f$-splitting. Let $d_a \in \mathcal{E}_{L(\alpha)}$
(respectively, $d_b \in \mathcal{E}_{L(\alpha)}$) be the unique
edge in $L(\alpha)$ labelled by $e_a$ (respectively, $e_b$).
Define an end-pointed $G$-labelled graph
\begin{equation*}T' := L(\alpha) \amalg L(U_a^+) \amalg L(U_a^-) \amalg L(U_b^+)
\amalg L(U_b^-) / \sim,\end{equation*} where $\sim$ equates
$\tau(d_a)$, $\iota(L(U_a^+))$ and $\iota(L(U_a^-))$ and equates
$\tau(d_b)$, $\iota(L(U_b^+))$ and $\iota(L(U_b^-))$. Define the
end-points of $T'$ to be the image of the end-points of
$L(\alpha)$. Define $T$ to be the end-pointed $G$-immersion
determined by $T'$.  Let $v_a$ (respectively, $v_a^{\pm}, v_b,
v_b^{\pm}$) denote the image of $\tau(d_a)$ (respectively,
$\tau(L(U_a^{\pm}))$, $\tau(d_b)$, $\tau(L(U_b^{\pm}))$) in $T$.

Suppose at least $K$ edges of $U_a^+$ and $U_b^+$ become
identified in $T$.  By Lemma \ref{BalloonsDontFoldMuch} (1),
$\mu_a = \mu_b$ and $\alpha' = \mu_a^k$ for some $k \in \Integer$.
Thus, $\alpha$ is an exceptional path --- a contradiction to
Corollary \ref{IRTTRepsAndGrowth}
(\ref{ExceptionalPathsLinearEdges}). Similarly, the folding
between any pair from the set $\{U_a^+, U_a^-, U_b^+, U_b^-\}$ is
limited by Lemma \ref{BalloonsDontFoldMuch}.  Hence the definition
of $l$ implies that the following properties hold:
\begin{itemize}
\item $T$ is a tree with six distinct ends; and
\item the six end-paths (see Definition \ref{EndsAndEndPaths})
have labels $e_a, e_b, \epsilon_a^+, \epsilon_a^-, \epsilon_b^+,
\epsilon_b^-$ such that each of $\epsilon_a^+, \epsilon_a^-,
\epsilon_b^+, \epsilon_b^-$ crosses $E_{h(u_a)}$.
\end{itemize}

Construct a $G$-immersion $T^{\ast}$ from $T$ by extending each
connected component of $p^{-1}(G_{h-1})$ to a connected
$G_{h-1}$-covering (if an end-path of $T$ has label $\rev{e}_{h}$,
then adjoin a $G_{h-1}$ cover at the corresponding end of $T$). It
follows immediately from (TT3) of the IRTT Theorem that the
following properties hold:
\begin{itemize}
\item [(A)]
$T^{\ast} \setminus \{D_a, D_b\}$ is $f$-stable with period $q_0$,
say.
\item [(B)]
$\alpha_{a, \, 0} \dots \alpha_{a, \, l-1}$ labels a
path from $v_a$ to $v_a^+$; \\
$\rev{\beta}_{a, \, l-1} \dots \rev{\beta}_{a, \, 0}$
labels a path from $v_a^-$ to $v_a$; \\
$\alpha_{b, \, 0} \dots \alpha_{b, \, l-1}$ labels a
path from $v_b$ to $v_b^+$; \\
$\rev{\beta}_{b, \, l-1} \dots \rev{\beta}_{b, \, 0}$
labels a path from $v_b^-$ to $v_b$;\\
$\alpha'$ labels a path from $v_a$ to $v_b$.
\end{itemize}

\noindent By Proposition \ref{PropertiesOfLambda(PrimaryForm)},
for sufficiently large $q_1 \in \Nat$, the following property
holds:
\begin{itemize}
\item [(C)] $\Lambda := \Lambda(\alpha_{a, \, l}, q_1)$ satisfies the
conclusions of Proposition \ref{PropertiesOfLambda(PrimaryForm)}.
\end{itemize}

\noindent By the Periodic Open Immersions Lemma (applied to
$\alpha_{a, \, 1}$, $\dots$, $\alpha_{a, \, s_a}$), there exist
$q_2 \in \Nat$ and a bi-infinite sequence of end-pointed
$G$-immersions $\{\Gamma_{a, \, i}\}_{i \in \, \Integer}$ such
that the following properties hold for each $i \in \Integer$:
\begin{itemize}
\item [(D)] Properties (1) and (2) of the Periodic Open Immersions Lemma are satisfied;
\item [(E)] $\Gamma_{a,
\, i}$ is $f$-stable with $\period(\Gamma_{a, \, i}) \, | \, q_2$;
\item [(F)] $\alpha_{a, \, i}$ labels a path across $\Gamma_{a, \, i}$.
\end{itemize}
Similarly, there exist $q_3 \in \Nat$ and $\{\Gamma_{b, \, i}\}_{i
\in \, \Integer}$ such that properties (D'), (E') and (F'),
analogous to (D), (E) and (F) respectively, hold. Choose $m \in
\Nat$ such that $q := m q_0 q_1 q_2 q_3 s_a s_b > 2l+1$. Define
end-pointed $G$-labelled graphs
$L_a := \combine{\Gamma_{a, \, l+1}, \Gamma_{a, \, l+2}, \dots,
\Gamma_{a, \, q-l-2}}$ and $L_b := \combine{\Gamma_{b, \, l},
\Gamma_{b, \, l+1}, \dots, \Gamma_{b, \, q-l-2}}.$  Properties (D)
and (D') imply that $L_a$ and $L_b$ are end-pointed
$G$-immersions. Properties (A), (B), (E), (F), (E') and (F') imply
that, for each $k \in \Integer$,
\begin{itemize}
\item [(G)]
$f^{kq}_{\#}(\alpha_{a, \, 0} \dots \alpha_{a, \, l-1})$ labels a
path from $v_a$ to $v_a^+$; \\
$f^{kq}_{\#}(\alpha_{a, \, l})$ labels a path $\tilde{\beta}_{kq}$
across $\Lambda$ such that $k \mapsto l^{\ab}(\tilde{\beta}_{kq})
\in p_{d-1}$;\\
$f^{kq}_{\#}(\alpha_{a, \, l} \dots \alpha_{a, \, q - l-2})$
labels a path across $L_a$; \\
$f^{kq}_{\#}(\alpha_{a, \, q-l-1} \dots \alpha_{a, \, q-1})$
labels a path from $v_a^-$ to $v_a$; \\
$f^{kq}_{\#}(\alpha_{b, \, 0} \dots \alpha_{b, \, l-1})$ labels a
path from $v_b$ to $v_b^+$;\\
 $f^{kq}_{\#}(\alpha_{b, \, l} \dots
\alpha_{b, \, q - l-2})$ labels a path across $L_b$; and\\
$f^{kq}_{\#}(\alpha_{b, \, q-l-1} \dots \alpha_{b, \, q-1})$
labels a path from $v_b^-$ to $v_b$.
\end{itemize}
\noindent Define \begin{equation*}\Sigma := T^{\ast} \amalg
\Lambda \amalg L_a \amalg L_b / \sim,\end{equation*} where $\sim$
equates $v_a^+$ (respectively, $\tau(\Lambda), \tau(L_a), v_b^+,
\tau(L_b)$) with $\iota(\Lambda)$ (respectively, $\iota(L_a),
v_a^-, \iota(L_b), v_b^-$). It follows from the construction of
$T^{\ast}$ and Properties (C), (D) and (D') that $\Sigma$ is a
$G$-immersion. Properties $(\Sigma1)$ and $(\Sigma2)$ are easily
verified. Properties $(\Sigma3)$ and $(\Sigma4)$ follow from
Property (G).
\end{proof}

\setcounter{figure}{\value{thm}} \stepcounter{thm}
\begin{figure}
\includegraphics{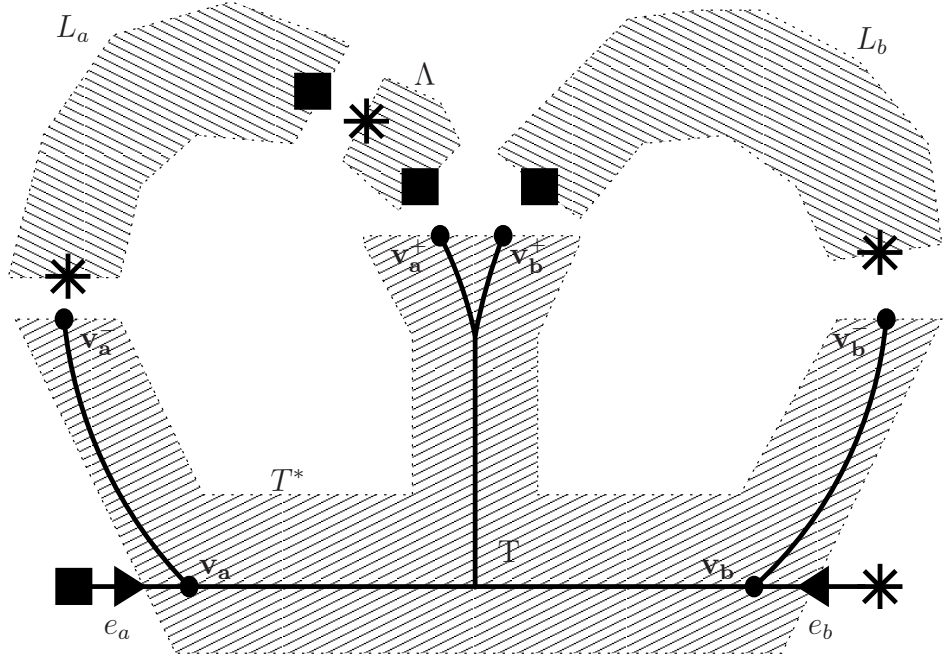}
\caption{A schematic depiction of a construction of $\Sigma$ in
Lemma \ref{TheQuadraticPathUnitLemmaCase1}.
\label{QuadraticCase1Figure}}
\end{figure}

\begin{lem}\label{TheQuadraticPathUnitLemmaCase2}
Proposition \ref{TheQuadraticPathUnitLemma} holds in the case that
$\alpha = e_a \alpha' \rev{e}_b$, $h(u_a) > h(u_b)$ and $h(u_a)
\geq h(\alpha')$.
\end{lem}

\begin{proof}
(An example construction of $\Sigma$ as below is illustrated
schematically in Figure \ref{QuadraticCase2Figure}). The proof is
similar to Lemma \ref{TheQuadraticPathUnitLemmaCase1} except that
we need not consider $L(U^+_b)$ or $L(U^-_b)$ in the construction
of $T$, or $L_b$ in the construction of $\Sigma$.

\end{proof}

\setcounter{figure}{\value{thm}} \stepcounter{thm}
\begin{figure}
\includegraphics{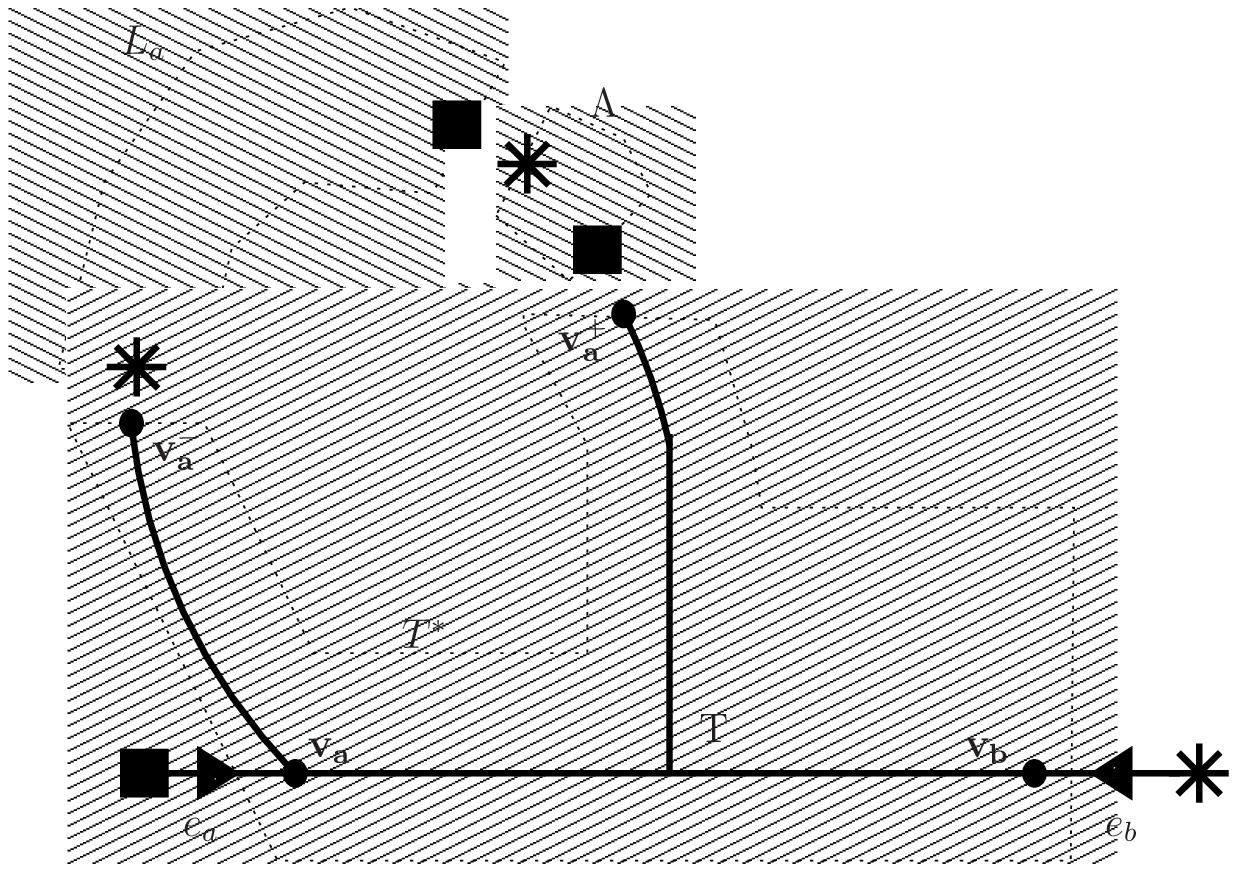}
\caption{A schematic depiction of a construction of $\Sigma$ in
Lemma \ref{TheQuadraticPathUnitLemmaCase2}.
\label{QuadraticCase2Figure}}
\end{figure}

\begin{lem}\label{TheQuadraticPathUnitLemmaCase3}
Proposition \ref{TheQuadraticPathUnitLemma} holds in the case that
$\alpha = e_a \alpha' \rev{e}_b$ and $h(\alpha') > h(u_a) =
h(u_b)$.
\end{lem}

\begin{proof}
(An example construction of $\Sigma$ as below is illustrated
schematically in Figure \ref{QuadraticCase3Figure}).  Define $h :=
h(\alpha')$ and $l := l(\alpha') + \max\{s_a, s_b\} + 1 $. Let
$U_a^+$ (respectively, $U_a^-, U_b^+, U_b^-$) be the initial
subpath of $S_a^+$ (respectively, $S_a^-, S_b^+, S_b^-$)
consisting of the first $l$ path units in the canonical
$f$-splitting. Let $d_a \in \mathcal{E}_{L(\alpha)}$
(respectively, $d_b \in \mathcal{E}_{L(\alpha)}$) be the unique
edge in $L(\alpha)$ labelled by $e_a$ (respectively, $e_b$).
Define an end-pointed $G$-labelled graph
\begin{equation*}T' := L(\alpha) \amalg L(U_a^+) \amalg L(U_a^-) \amalg L(U_b^+)
\amalg L(U_b^-) / \sim,\end{equation*} where $\sim$ equates
$\tau(d_a)$, $\iota(L(U_a^+))$ and $\iota(L(U_a^-))$ and equates
$\tau(d_b)$, $\iota(L(U_b^+))$ and $\iota(L(U_b^-))$. Define the
end-points of $T'$ to be the image of the end-points of
$L(\alpha)$. Define $T$ to be the end-pointed $G$-immersion
determined by $T'$. Let $v_a$ (respectively, $v_a^{\pm}, v_b,
v_b^{\pm}$) denote the image of $\tau(d_a)$ (respectively,
$\tau(L(U_a^{\pm}))$, $\tau(d_b)$, $\tau(L(U_b^{\pm}))$) in $T$.

Since $h(\alpha') > h(u_a) = h(u_b)$, there is at least one edge
of $L(\alpha')$ which acts as a sentinel, ensuring that $\bigl(
L(U_a^+) \cup L(U_a^+) \bigr) \cap \bigl( L(U_b^+) \cup L(U_b^+)
\bigr) = \emptyset$. Combined with the definition of $l$, this
implies that the following properties hold:
\begin{itemize}
\item [(A)]$T$ is a tree with six distinct ends;
\item [(B)] the six end-paths (see Definition \ref{EndsAndEndPaths})
have labels $e_a, e_b, \epsilon_a^+, \epsilon_a^-, \epsilon_b^+,
\epsilon_b^-$ such that each of $\epsilon_a^+, \epsilon_a^-,
\epsilon_b^+, \epsilon_b^-$ contains an essentially unbounded
subpath.
\end{itemize}

We define a finite sequence of $G$-immersions $T = T_{h+1}$, $T_h,
\dots, T_{h(u_a)} = T^+$ inductively as follows:  let $q_i
> \diam(T_{i+1}) + 2l(\mu_i)$, let $d_1, d_2, \dots, d_N$ be a
complete list of the edges in $T_{i+1}$ with label $e_i$, for each
$j = 1, 2, \dots, N$, let $\beta_i^j$ be a copy of
$C(\mu_{i}^{q_i})$ and define
\begin{equation*}T'_i := T_{i+1} \amalg \beta_i^1 \amalg \beta_i^2 \amalg
\dots \amalg \beta_i^N / \sim,\end{equation*} where $\sim$
identifies the basepoint of $\beta_i^j$ with $\tau(d_j)$ for each
$j = 1, 2, \dots, N$.  Define the end-points of $T'_i$ to be those
inherited from $T_{i+1}$.  Let $T_i$ be the end-pointed
$G$-immersion determined by $T'_i$.  We claim that the natural map
$T_{i+1} \to T_i$ is an embedding: let $T''_i$ be the end-pointed
$G$-labelled graph obtained from $T'_i$ by performing all folding
possible where one edge is from $\beta_i^1 \amalg \beta_i^2 \amalg
\dots \amalg \beta_i^N$ and the other from $T_{i+1}$.  By the
definition of $q_i$, for each $j = 1, 2, \dots, N$, at least
$2l(\mu_i) + 1$ edges of $\beta_i^j$ do not fold with the image of
$T_{i+1}$. Thus $T''_i$ consists of a copy of $T_{i+1}$ and $N$
distinct handles $h_1, \dots, h_N$ of length at least $2l(\mu_i) +
1$. Since the label on each $\beta_i^j$ is periodic, it follows
that if two such handles fold for more than $l(\mu_i)$ edges, then
the handles may be identified by folding and the end-points of the
handles are identical (this may only happen if there exists a path
$\gamma$ in $T''_i$ with label $e_i \mu_i^k \rev{e}_i$). It
follows that $T_i$ is obtained from $T''_i$ by folding some parts
of the handles $h_1, \dots, h_N$ (possibly identifying some
handles). Thus the natural map $T_{i+1} \to T_i$ is an embedding.
Hence we have that the following property holds:
\begin{itemize}
\item [(C)] the natural map $T \to T^+$ is an embedding.
\end{itemize}

\noindent We claim that $T^+$ also has the following property:
\begin{itemize}
\item [(D)] if $\epsilon_a^+$ has initial edge $e_{h(u_a)}$ then there
exists a unique edge $d \in \mathcal{E}_{T^+}$ such that $\hat{d}
= e_{h(u_a)}$ and $\iota(d) = v_a^+$; otherwise, there is no such
edge $d \in \mathcal{E}_{T^+}$
\end{itemize}
Suppose that $\epsilon_a^+$ has initial edge $e_{h(u_a)}$. That
there exists at least one edge $d \in \mathcal{E}_{T^+}$ with the
required property is immediate by Property (C).  Now consider the
inductive construction of $T^+$.  By Remark
\ref{PropertiesOfEssentiallyUnbounded}, a Nielsen path contains no
essentially unbounded subpaths.  It follows from Property (B) that
for each $i > h(u_a)$ and each $j$, $\beta_i^j$ does not fold past
the essentially unbounded subpath in $\epsilon_a^+$.  Analogous
properties (D'), (D'') and (D''') hold for $\epsilon_a^-$,
$\epsilon_b^+$ and $\epsilon_b^-$ respectively.

Construct a $G$-immersion $T^{\ast}$ from $T^+$ by extending each
connected component of $p^{-1}(G_{h-1})$ to a $G_{h-1}$-covering
(if an end-path of $T^+$ has label $\rev{e}_{h}$, then adjoin a
$G_{h-1}$ cover at the corresponding end of $T^+$).  It follows
from our construction and (TT3) of the IRTT Theorem that the
following properties hold:
\begin{itemize}
\item [(E)]
$T^{\ast} \setminus \{D_a, D_b\}$ is $f$-stable with period $q_0$,
say.
\item [(F)]
$\alpha_{a, \, 0} \dots \alpha_{a, \, l-1}$ labels a
path from $v_a$ to $v_a^+$; \\
$\rev{\beta}_{a, \, l-1} \dots \rev{\beta}_{a, \, 0}$
labels a path from $v_a^-$ to $v_a$; \\
$\alpha_{b, \, 0} \dots \alpha_{b, \, l-1}$ labels a
path from $v_b$ to $v_b^+$; \\
$\rev{\beta}_{b, \, l-1} \dots \rev{\beta}_{b, \, 0}$
labels a path from $v_b^-$ to $v_b$; and\\
$\alpha'$ labels a path from $v_a$ to $v_b$.
\end{itemize}

The rest of the proof proceeds as in Lemma
\ref{TheQuadraticPathUnitLemmaCase1}, with Properties (D), (D'),
(D'') and (D''') used to show that the construction $\Sigma$ is a
$G$-immersion.
\end{proof}

\setcounter{figure}{\value{thm}} \stepcounter{thm}
\begin{figure}
\includegraphics{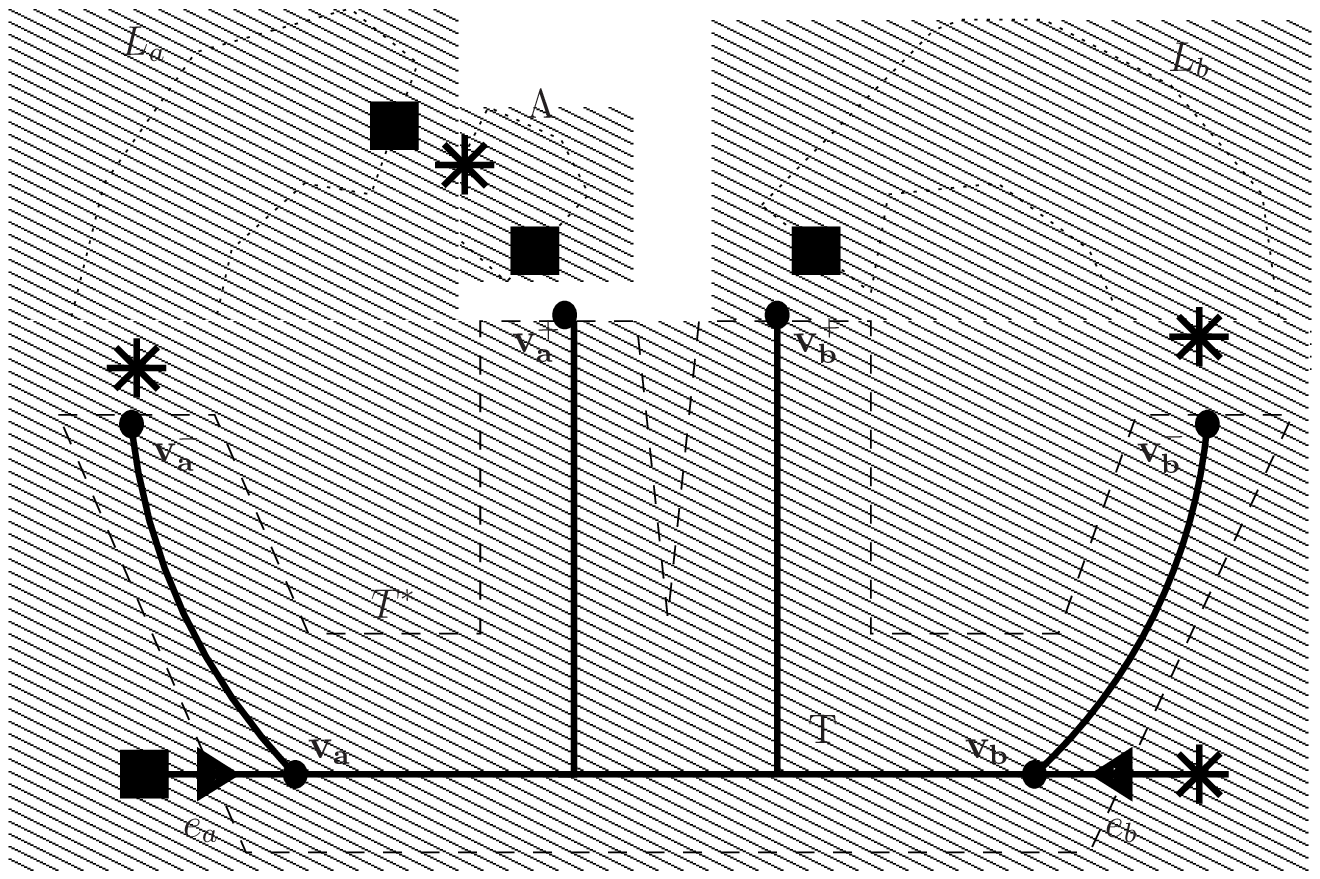}
\caption{A schematic depiction of a construction of $\Sigma$ in
Lemma \ref{TheQuadraticPathUnitLemmaCase3}.
\label{QuadraticCase3Figure}}
\end{figure}

\begin{lem}\label{TheQuadraticPathUnitLemmaCase4}
Proposition \ref{TheQuadraticPathUnitLemma} holds in the case that
$\alpha = e_a \alpha' \rev{e}_b$ and $h(\alpha') > h(u_a) >
h(u_b)$.
\end{lem}

\begin{proof}
The proof is similar to Lemma \ref{TheQuadraticPathUnitLemmaCase3}
except that we need not consider $L(U^+_b)$ or $L(U^-_b)$ in the
construction of $T$, or $L_b$ in the construction of $\Sigma$.
\end{proof}

\subsection{The case $d \geq 3$}

\begin{lem}[The Tree Lemma]\label{TheTreeLemma}\index{wTreeLemma@Tree Lemma}
Let $T$ be a finite $G$-immersion which is a tree and let $d \geq
2$ be such that $L_d \leq h(T) < L_{d+1}$. We may extend $T$ to a
$G$-immersion $T^{\ast}$ such that the following properties hold:
\begin{enumerate}
\item $h(T^{\ast} \setminus T) < L_d$;
\item $T^{\ast}$ is $f$-stable.
\end{enumerate}
\end{lem}

\begin{proof}
Let $D_1, D_2, \dots, D_s$ be a complete list of the geometric
edges in $T$ such that $\hat{D_i} \in \{E_{L_d}, E_{L_d +1},
\dots, E_{L_{d+1}-1}\}$. Construct $T^{\ast}$ from $T$ by
extending each connected component of $T \setminus (\cup_{i = 1}^s
D_i)$ to a $G_{L_d - 1}$-cover by Stallings' Algorithm, and
adjoining a cover of $G_{L_d-1}$ at any end of $T$ for which the
corresponding end-path has initial edge with label in
$\{\rev{e}_{L_d}, \rev{e}_{L_d +1}, \dots, \rev{e}_{L_{d+1}-1}\}$.
\end{proof}

By Remark \ref{HowToFinishProof}, the following lemma completes
the proof of the Apt Immersion Theorem and the Main Theorem.

\begin{prop}[The Path Unit Proposition]\label{ThePathUnitLemma}\index{wAptImmersionTheorem@Apt Immersion Theorem!wproofincasedgeq3@proof in the case $d \geq 3$}
Let $\alpha \subset G$ be a path unit of degree $d \geq 3$. There
exist an end-pointed $G$-immersion $\Sigma$ and $q \in \Nat$ such
that the following conditions hold:
\begin{enumerate}
\item [$(\Sigma 1)$] $\iota(\Sigma) \neq \tau(\Sigma)$;
\item [$(\Sigma 2)$] if $\alpha$ has type (i) then both $\iota(\Sigma)$ and
$\tau(\Sigma)$ have valence 1; if $\alpha$ has type (ii) then
$\iota(\Sigma)$ has valence 1;  if $\alpha$ has type (iii) then
$\tau(\Sigma)$ has valence 1;
\item [$(\Sigma 3)$] for each non-negative integer $k$, $f^{kq}_{\#}(\alpha)$
labels a path $\tilde{\alpha}_{kq}$ across $\Sigma$; and
\item [$(\Sigma 4)$] $k \mapsto l^{\ab}(\tilde{\alpha}_{kq}) \in p_d$.
\end{enumerate}
\end{prop}

\begin{proof}
(Figure \ref{QuadraticCase1Figure} can be reused to illustrate
schematically an example construction of $\Sigma$ as below).  We
use induction on $d$, the degree of the path unit. The case that
$d = 2$ has been completed in Proposition
\ref{TheQuadraticPathUnitLemma}.  Assume the result holds for each
path unit of degree $d-1$, for some $d \geq 3$. Let $\alpha
\subset G$ be a path unit of degree $d$.  We will complete the
inductive step in the case that $\alpha$ has type (i) (the most
difficult case). The case that $\alpha$ has type (ii) is proved by
an argument similar to that executed below. By considering
$\rev{\alpha}$ instead of $\alpha$, it is clear the case that
$\alpha$ has type (iii) is equivalent to the case that $\alpha$
has type (ii).

Assume $\alpha = e_a \alpha' \rev{e}_b$ for some $L_d \leq a, b <
L_{d+1}$ and some $\alpha' \subset G_{d-1}$. Consider $L(\alpha)$.
Let $d_a, d_b \in \mathcal{E}_{L(\alpha)}$ be the edges labelled
by $e_a$ and $e_b$ respectively. Let $g$ be the number of path
units in the canonical $f$-splitting of $\alpha'$ and define $K :=
s_a s_b + \min \{s_a, s_b\} + 1$. Choose $l \in \Nat$ such that $l
\geq g + 2 K$ and $\alpha_{a, \, l}$ is a path unit of degree
$d-1$. Let $U_a^+$ (respectively, $U_a^-, U_b^+, U_b^-$) be the
initial subpath of $S_a^+$ (respectively, $S_a^-, S_b^+, S_b^-$)
consisting of the first $l$ path units in the canonical
$f$-splitting of $S_a^+$ (respectively, $S_a^-, S_b^+, S_b^-$).
Define an end-pointed $G$-labelled graph \begin{equation*}T' :=  T
\amalg L(U_a^+) \amalg L(U_a^-) \amalg L(U_b^+) \amalg L(U_b^-) /
\sim,\end{equation*} where $\sim$ equates the initial point of
$L(U_a^+)$ (respectively, $L(U_a^-)$, $L(U_b^+)$, $L(U_b^-)$) with
$\tau(d_a)$ (respectively, $\tau(d_a)$, $\tau(d_b)$, $\tau(d_b)$).
Let $v_a$ (respectively, $v_a^{\pm}$, $v_b$, $v_b^{\pm}$) denote
the image of $\tau(d_a)$ (respectively, $\tau(L(U_a^{\pm}))$,
$\tau(d_b)$, $\tau(L(U_b^{\pm}))$) in $T$. Let $T$ denote the
end-pointed $G$-immersion determined by $T'$. Let $T''$ be
obtained from $T'$ by performing all folding possible where one
edge is from $L(U_a^-) \cup L(U_a^+) \cup L(U_b^+) \cup L(U_b^-)$
and the other from $L(\alpha)$.  It follows from our hypothesis on
$l$ that $T'$ is a tree with 6 distinct ends and at least $s_a
s_b$ complete path units of $L(U_a^+)$ (respectively, $L(U_a^-),
L(U_b^+), L(U_b^-)$) remain unfolded. It follows from Lemma
\ref{BalloonsDontFoldMuch} that $T$ (which is also the
$G$-immersion determined by $T''$) is a tree with 6 distinct ends.
By the Tree Lemma, we may extend $T$ to a $G$-immersion $T^{\ast}$
such that the following properties hold:
\begin{itemize}
\item [(A)] $h(T^\ast \setminus T) < L_d$;
\item [(B)] $T^{\ast} \setminus \{D_a, D_b\}$ is $f$-stable, with period $q_0$, say.
\end{itemize}

\noindent By the inductive hypothesis the following property
holds:
\begin{itemize}
\item [(C)] there exist an end-pointed $G$-immersion $\Lambda$
and $q_1 \in \Nat$ such that the conclusion of the Path Unit Lemma
hold with $\alpha_{a, \, l}$ in place of $\alpha$, $\Lambda$ in
place of $\Sigma$ and $d-1$ in place of $d$.
\end{itemize}

\noindent By the Periodic Open Immersions Lemma (applied to
$\alpha_{a, \, 1}$, $\dots$, $\alpha_{a, \, s_a}$), there exist
$q_2 \in \Nat$ and a bi-infinite sequence of end-pointed
$G$-immersions $\{\Gamma_{a, \, i}\}_{i \in \, \Integer}$ such
that the following properties hold for each $i \in \Integer$:
\begin{itemize}
\item [(D)] Properties (1) and (2) of the Periodic Open Immersions Lemma are satisfied;
\item [(E)] $\Gamma_{a,
\, i}$ is $f$-stable with $\period(\Gamma_{a, \, i}) \, | \, q_2$;
\item [(F)] $\alpha_{a, \, i}$ labels a path across $\Gamma_{a, \, i}$.
\end{itemize}
Similarly, there exist $q_3 \in \Nat$ and $\{\Gamma_{b, \, i}\}_{i
\in \, \Integer}$ such that Properties (D'), (E') and (F'),
analogous to (D), (E) and (F) respectively, hold. Choose $m \in
\Nat$ such that $q := m q_0 q_1 q_2 q_3 s_a s_b > 2l+1$. Define
end-pointed $G$-labelled graphs $L_a := \combine{\Gamma_{a, \,
l+1}, \Gamma_{a, \, l+2}, \dots, \Gamma_{a, \, q-l-2}}$ and $L_b
:= \combine{\Gamma_{b, \, l}, \Gamma_{b, \, l+1}, \dots,
\Gamma_{b, \, q-l-2}}.$  Properties (D) and (D') imply that $L_a$
and $L_b$ are end-pointed $G$-immersions. Properties (B), (E),
(F), (E') and (F') imply that, for each $k \in \Integer$,
\begin{itemize}
\item [(G)]
$f^{kq}_{\#}(\alpha_{a, \, 0} \dots \alpha_{a, \, l-1})$ labels a
path from $v_a$ to $v_a^+$; \\
$f^{kq}_{\#}(\alpha_{a, \, l})$ labels a path $\tilde{\beta}_{kq}$
across $\Lambda$ such that $k \mapsto l^{\ab}(\tilde{\beta}_{kq})
\in p_{d-1}$;\\
$f^{kq}_{\#}(\alpha_{a, \, l} \dots \alpha_{a, \, q - l-2})$
labels a path across $L_a$; \\
$f^{kq}_{\#}(\alpha_{a, \, q-l-1} \dots \alpha_{a, \, q-1})$
labels a path from $v_a^-$ to $v_a$; \\
$f^{kq}_{\#}(\alpha_{b, \, 0} \dots \alpha_{b, \, l-1})$ labels a
path from $v_b$ to $v_b^+$;\\
 $f^{kq}_{\#}(\alpha_{b, \, l} \dots
\alpha_{b, \, q - l-2})$ labels a path across $L_b$;\\
$f^{kq}_{\#}(\alpha_{b, \, q-l-1} \dots \alpha_{b, \, q-1})$
labels a path from $v_b^-$ to $v_b$.
\end{itemize}
\noindent Define \begin{equation*}\Sigma := T^{\ast} \amalg
\Lambda \amalg L_a \amalg L_b / \sim,\end{equation*} where $\sim$
equates $v_a^+$ (respectively, $\tau(\Lambda), \tau(L_a), v_b^+,
\tau(L_b)$) with $\iota(\Lambda)$ (respectively, $\iota(L_a),
v_a^-, \iota(L_b), v_b^-$). It follows from the construction of
$T^{\ast}$ and Properties (A), (C), (D) and (D') that $\Sigma$ is
a $G$-immersion. Properties $(\Sigma1)$ and $(\Sigma2)$ are easily
verified. Properties $(\Sigma3)$ and $(\Sigma4)$ follow from
Property (G).
\end{proof}

\begin{rem}\label{WhyDifferentCases}\index{wAptImmersionTheorem@Apt Immersion Theorem!wstrategyofproof@strategy of proof}
In this remark we clarify why the proof of the quadratic case of
the Apt Immersion Theorem (Proposition
\ref{TheQuadraticPathUnitLemma}) is separate from the proof in the
case that $d \geq 3$ (Proposition \ref{ThePathUnitLemma}), and why
the former case is further split into four sub-cases (the lemmas
\ref{TheQuadraticPathUnitLemmaCase1},
\ref{TheQuadraticPathUnitLemmaCase2},
\ref{TheQuadraticPathUnitLemmaCase3},
\ref{TheQuadraticPathUnitLemmaCase4}).

In the proof of Proposition \ref{ThePathUnitLemma}, we construct a
$G$-immersion $T$ such that $L_{d-1} \leq h(T) < L_d$. We extend
$T$ to a $G$-immersion $T^{\ast}$ by applying Stallings' Algorithm
to extend certain connected subgraphs of $T$ which are
$G_{h(T)-1}$-immersions to $G_{h(T)-1}$-covers.  The new edges in
the extension (that is, edges in $T^\ast \setminus T$) have height
at most $L_{d-1}-1$. Thus, for each $i \geq L_{d-1}$, each edge
with label $e_i$ in $T$ acts as a sentinel in $T$, limiting the
amount that new edges may fold with edges of $T$. In the quadratic
case, we construct a $G$-immersion $T$ (or $T^+$ in the case of
Lemma \ref{TheQuadraticPathUnitLemmaCase3}) such that $L_1 \leq
h(T) < L_2$. Again, we extend $T$ to a $G$-immersion $T^{\ast}$ by
applying Stallings' Algorithm to extend certain connected
subgraphs of $T$ which are $G_{h(T)-1}$-immersions to
$G_{h(T)-1}$-covers. Since it is possible that, for an integer $i$
such that $L_1 \leq i < L_2$, $h(u_i) \geq L_1$, it is not
necessarily the case that every linearly growing edge in $T$ acts
as a sentinel in the way that edges of degree $d-1$ did in the
case that $d \geq 3$. Thus the quadratic case is more subtle than
the case that $d \geq 3$ and is dealt with separately. In the
quadratic case, observe that each edge with label $e_{h(T)}$ does
act as a sentinel in $T$. The quadratic case is split into four
sub-cases, depending on how such sentinel edges arise, and in the
cases of Lemma \ref{TheQuadraticPathUnitLemmaCase3} and Lemma
\ref{TheQuadraticPathUnitLemmaCase4}, how the edges with label
$e_{h(u_a)}$ arise.
\end{rem}

\nocite{Grossman, ModuliOfGraphs, Stallings, BFH1, BFH2,
TrainTracks, MarshallHall, ThurstonSurfaces}
\bibliographystyle{amsplain}
\bibliography{GrowthPaperBib}

\providecommand{\bysame}{\leavevmode\hbox to3em{\hrulefill}\thinspace}
\providecommand{\MR}{\relax\ifhmode\unskip\space\fi MR }
\providecommand{\MRhref}[2]{%
  \href{http://www.ams.org/mathscinet-getitem?mr=#1}{#2}
}
\providecommand{\href}[2]{#2}
\begin{thebibliography}{10}

\bibitem{BFH2}
M.~Bestvina, M.~Feighn, and M.~Handel, \emph{{The} {Tits} alternative for
  {$Out(F_n)$} {II}: {A} {Kolchin} type theorem}, Preprint, Univerity of Utah,
  1999.

\bibitem{BFH1}
\bysame, \emph{{The} {Tits} alternative for {$Out(F_n)$} {I}: {Dynamics} of
  exponentially growing automorphisms}, Ann. of Math. \textbf{151} (2000),
  517--623.

\bibitem{TrainTracks}
M.~Bestvina and M.~Handel, \emph{{Train} tracks and automorphisms of free
  groups}, Ann. of Math. \textbf{135} (1992), 1--51.

\bibitem{PolynomialDehnFunctions}
M.R. Bridson, \emph{{Polynomial} {Dehn} functions and the length of
  asynchronously automatic structures}, Proc. London Math. Soc. \textbf{3}
  (2002), no.~85, 441--466.

\bibitem{ModuliOfGraphs}
Marc Culler and Karen Vogtmann, \emph{Moduli of graphs and automorphisms of
  free groups}, Invent. math. \textbf{84} (1986), 91--119.

\bibitem{Grossman}
Edna~K. Grossman, \emph{Representations of the automorphism groups of free
  groups}, J. Algebra \textbf{30} (1974), 388--399.

\bibitem{MarshallHall}
Marshall~Hall{,} Jr., \emph{Coset representations in free groups}, Trans. Amer.
  Math. Soc. \textbf{67} (1949), 421--432.

\bibitem{Lubotzky}
Alexander Lubotzky, \emph{Normal automorphisms of free groups}, Journal of
  Algebra \textbf{63} (1980), 494--498.

\bibitem{Macura}
N.~Macura, \emph{Quadratic isoperimetric inequality for mapping tori of
  polynomially growing automorphisms of free groups}, Geom. Funct. Annal.
  \textbf{10} (2000), 874--901.

\bibitem{Stallings}
John~R. Stallings, \emph{{Topology} of finite graphs}, Invent. Math.
  \textbf{71} (1983), 551--565.

\bibitem{ThurstonSurfaces}
W.P. Thurston, \emph{On the geometry and dynamics of diffeomorphisms of
  surfaces}, Bull. A.M.S. \textbf{19} (1988), 417--431.

\end{thebibliography}
\printindex
\end{document}